\documentclass[11pt]{article}
\usepackage{latexsym}
\usepackage{amsmath}
\renewcommand{\epsilon}{\varepsilon}

\usepackage{amsmath,bm}
\usepackage{amssymb}
\usepackage{amsthm}
\usepackage[ansinew]{inputenc}
\usepackage{graphicx}
\usepackage{epsfig}
\usepackage{dsfont}
\usepackage{bbm}
\usepackage{subfig}
\usepackage{url}

\usepackage[authoryear]{natbib}
\newtheorem{satz}{Theorem}[section]

\newtheorem{lem}[satz]{Lemma}

\newtheorem{rem}[satz]{Remark}
\newtheorem{assumption}[satz]{Assumption}

\def\3{\ss}
\def\er{\mathbb{R}}

\def\en{\mathbb{N}}

\newcommand{\E}{\mathbbm{E}}

\newcommand{\bea}{\begin{eqnarray*}}
\newcommand{\eea}{\end{eqnarray*}}
\newcommand{\be}{\begin{eqnarray}}
\newcommand{\ee}{\end{eqnarray}}

\newcommand{\ba}{\begin{array}}
\newcommand{\ea}{\end{array}}
\newcommand{\cum}{\text{\rm cum}}

\newcommand{\Cov}{\text{\rm Cov}}
\newcommand{\Var}{\text{\rm Var}}
\def\3{\ss}
\def\er{\mathbb{R}}

\def\en{\mathbb{N}}

\allowdisplaybreaks

\begin{document}

\title{Detecting long-range dependence in non-stationary time series}

\author{Holger Dette, Philip Preu\ss, Kemal Sen\\
Ruhr-Universit\"at Bochum \\
Fakult\"at f\"ur Mathematik \\
44780 Bochum \\
Germany \\
{\small email: holger.dette@ruhr-uni-bochum.de} \\
{\small email: philip.preuss@ruhr-uni-bochum.de}  \\
{\small email: kemal.sen@ruhr-uni-bochum.de}
}

 \maketitle
\begin{abstract}
An important problem in time series analysis
is the discrimination between non-stationarity and long-range dependence.
Most of the literature considers the problem of testing specific
parametric hypotheses of non-stationarity (such as a change in the mean)
against long-range dependent stationary alternatives. In this paper we suggest a simple
 approach, which can be used to test the null-hypothesis of a general
non-stationary short-memory against the alternative of a non-stationary long-memory  process.
The test procedure works in the spectral domain and uses a sequence of approximating tvFARIMA models
to estimate the time varying long-range dependence parameter. We prove uniform consistency of this estimate
and  asymptotic normality of an averaged version. These results yield a simple test (based on the quantiles of the
standard normal distribution), and  it  is demonstrated in a simulation study  that - despite of its semi-parametric nature - the new
test outperforms the currently available methods, which are constructed to discriminate between  specific
parametric hypotheses of non-stationarity short- and stationarity long-range dependence.

\end{abstract}

AMS subject classification: 62M10, 62M15, 62G10

Keywords and phrases: spectral density, long-memory, non-stationary processes, goodness-of-fit tests, empirical spectral measure,
integrated periodogram, locally stationary process, approximating models

\section{Introduction}
\def\theequation{1.\arabic{equation}}
\setcounter{equation}{0}

\begin{figure}
\begin{center}
\includegraphics[width=0.4\textwidth]{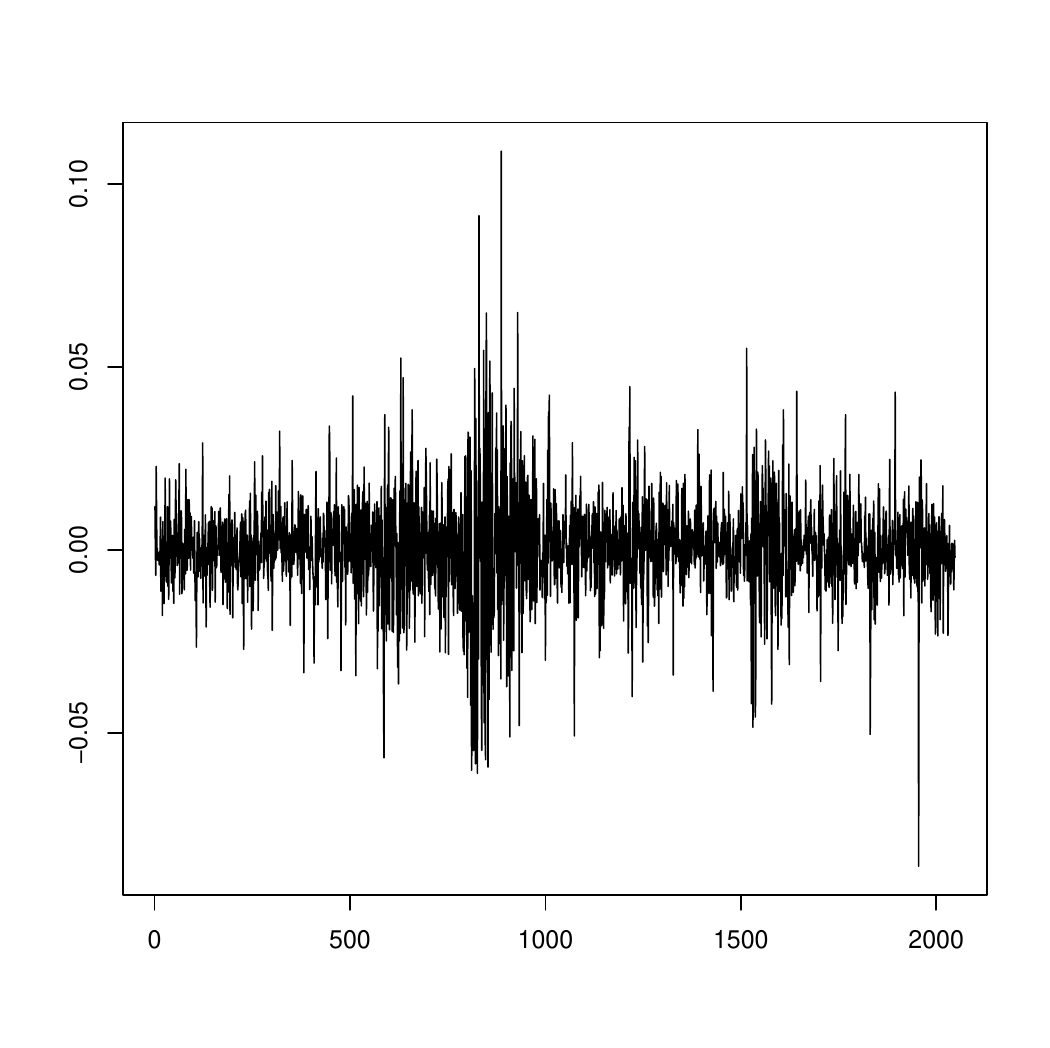}~~
\includegraphics[width=0.4\textwidth]{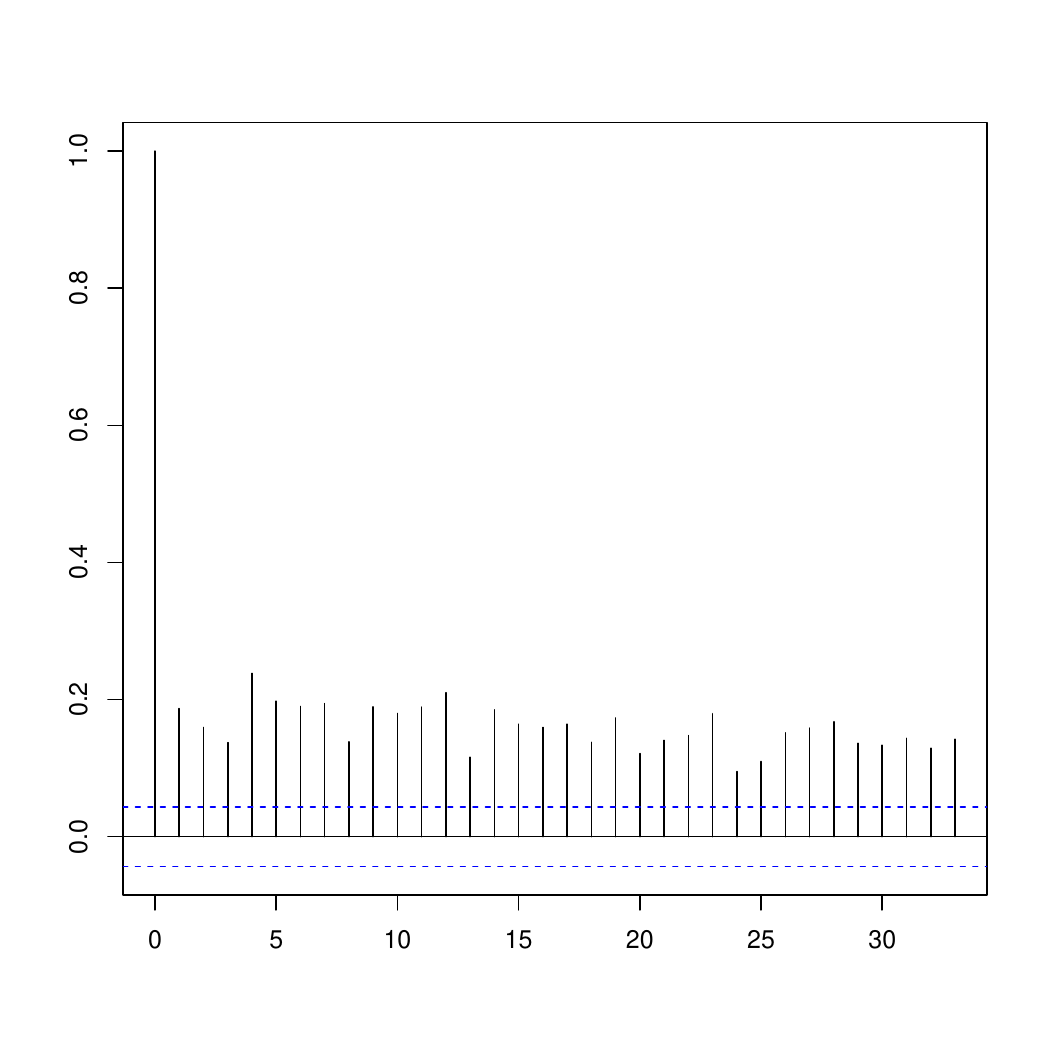}~~
  \caption{ \label{ibm} {\it Left panel: log-returns of the IBM stock between July 15th 2005 and  August 30th 2013;
  right panel: Sample autocovariance function of the squared returns $X_t^2$} }
\end{center}
\end{figure}

\begin{figure}
\begin{center}
\includegraphics[width=0.4\textwidth]{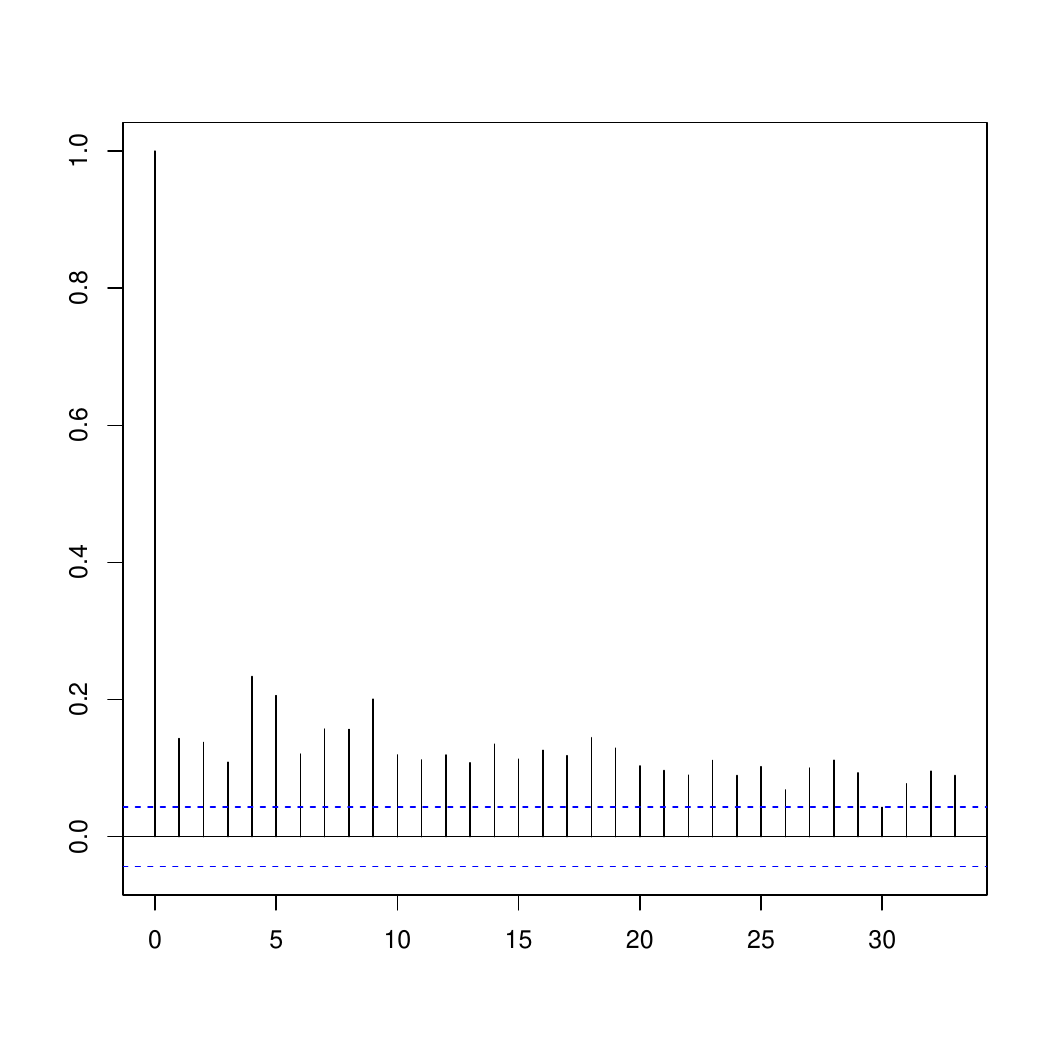}~~
\includegraphics[width=0.4\textwidth]{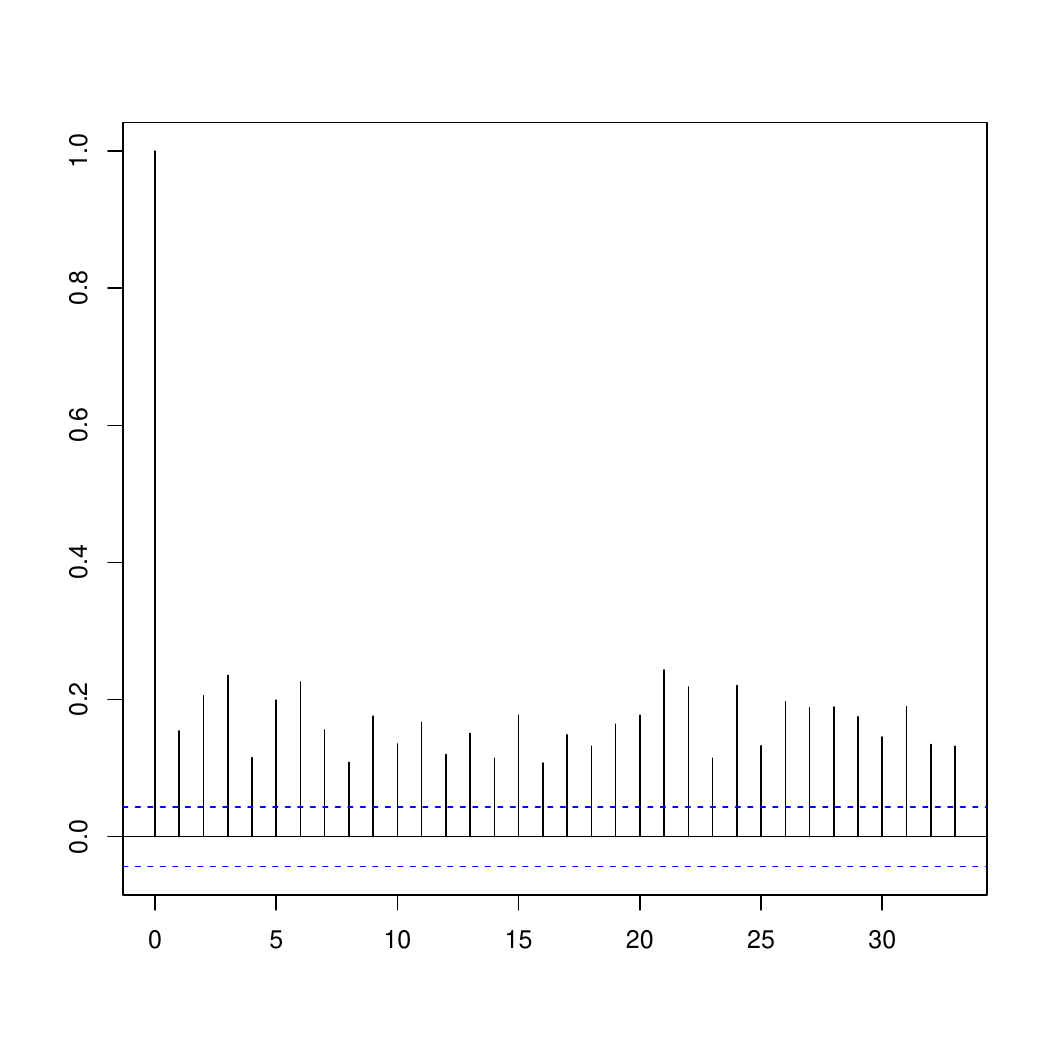}~~
  \caption{ \label{ibmfit} {\it Left panel: Sample autocovariance function of a simulated time series from a FARIMA(3,d,0)-fit to the $2048$ squared IBM-returns $X_t^2$, right panel: Sample autocovariance function of $X_t^2$ for $X_t$ simulated from the model $X_{t,T}=\hat \sigma(t/T) Z_t$ with $\hat \sigma (\cdot)$ estimated by a rolling-window of length $128$.} }
\end{center}
\end{figure}

Many time series  [like asset volatility or regional temperatures]
exhibit a slow decay in the sample
 autocorrelation function and
simple
 stationary short-memory models can not be used to analyze this type of data.
 A typical example is displayed in Figure \ref{ibm}, which shows
 $2048$ log-returns of the IBM stock between July $15$th $2005$ and  August $30$th $2013$,
with estimated autocovariance function of the squared returns $X_t^2$.
In this example  the assumption of stationarity with a summable sequence of autocovariances, say
 $(\gamma(k))_{k \in \mathbb{N}}$,  is hard to justify for the volatility process. Long-range dependent processes have been introduced as an attractive alternative to model features of this type
using an autocovariance function with the property
\begin{eqnarray*}
\gamma(k) \sim C k^{2d-1}
\end{eqnarray*}
as $k \rightarrow \infty$, where $d \in (0,0.5)$  denotes a ``long memory'' parameter. Statistical models (and corresponding theory) for long-range dependent processes are  very well developed [see \cite{book1} or \cite{book2} for recent surveys]
 and have found applications in numerous fields
 [see \cite{breicrali}, \cite{anwendung2} or \cite{anwendung3} for such an approach in the framework of asset volatility, video traffic and wind power modeling]. However, it was pointed out by several authors that the observation of ``long memory'' features in the sample autocovariance function can be as well explained by
non stationarity [see \cite{mikosch2004} or \cite{motivation2} among many others].
This is clearly demonstrated in Figure \ref{ibmfit}, which shows the sample autocovariances
of  the squared returns from a fit of the (non-stationary) model  $X_{t,T}=\sigma(t/T) Z_t$ for the returns [here  $Z_t$ is an i.i.d. sequence and $\sigma(\cdot)$ is piecewise-constant, cf. \cite{starica2005} or  \cite{fryzlewicz2006} for more details], and from a
 stationary FARIMA($3,d,0$)-fit for the squared ones $X_t^2$. Both  models  are able to explain the observed effect of 'long-range dependence' for the volatility process. So, in summary, the same effect can be explained by two completely different modeling approaches.

For this reason several authors have pointed  out the importance to distinguish between long-memory and non-stationarity
 [see \cite{starica2005}, \cite{motivation1} or \cite{motivation2} to mention only a few]. However, there exists a surprisingly small number of statistical procedures which address problems of this type.
 To the best of our knowledge,  \cite{kuensch1986} is the first reference
 investigating the existence of  ``long memory''  if non-stationarities appear in the time series.
 In this article a procedure to discriminate between a long-range dependent model and
 a process with a monotone mean functional and weakly dependent innovations is derived. Later on, \cite{heydedai} developed a method for distinguishing between long-memory and small trends. \cite{sibbertsen2009} tested the null hypothesis of a constant long-memory parameter against a break in the long-memory parameter.  Furthermore, \cite{horvarth2006}, \cite{baekpipi} and \cite{Yau2012} investigated CUSUM and likelihood ratio tests to discriminate between the null hypothesis of no long-range
 and weak dependence with one change point in the mean.

Although the procedures  proposed in these articles are technically mature and work rather well in suitable situations, they are, however,  only designed to discriminate between long-range dependence and a very specific change in the first-order structure, like one structural break and two stationary segments of the series. This is rather restrictive, since the expectation might change in a different way than assumed [there could be, for example, continuous changes or multiple breaks instead of a single one] and the second-order structure could be time-varying as well. However, if these or more general non-stationarities occur, the  discrimination techniques, which have been proposed in the literature so far,  usually fail, and a procedure which is working under less restrictive assumptions is still missing. \\
The objective of this  paper is to fill this gap and to develop a test for the null hypothesis of no long-range dependence in a framework which is flexible enough to deal with different types of
non-stationarity
 in both the first and second-order structure.  The general model is introduced in Section \ref{sec2}.
Our approach uses an estimate of a (possibly time varying) long-range dependence parameter, which
 is derived by a sequence of approximating tvFARIMA models with a slightly enlarged  parameter space.
This statistic estimates a functional which vanishes if and only if the null hypothesis of a short-memory locally stationary process is satisfied.
The method is based on some non-intuitive features of averages of
unconstrained estimators  in models  with a constrained  parameter space, which become
clear from the rather technical proofs given in Section \ref{appendix}. In order to make
these phenomena also visible  to readers which are less familiar with the technical machinery used
for the asymptotic analysis of non-stationary long range dependent processes we provide in Section \ref{sec2A} a motivation of our approach in the context of the classical
nonparametric regression model with repeated observations. \\
In Section \ref{sec3} we return to the locally stationary long range dependent time series model and
prove consistency and asymptotic normality of a corresponding test statistic under the null hypothesis of no long-range dependence.
As a consequence we obtain a nonparametric test, which is based on the quantiles of the standard normal distribution and therefore very easy to implement. The finite sample properties of the new test are investigated in Section \ref{sec4}, which also provides a comparison with the competing procedures with a focus on non-stationarities. We demonstrate the superiority of the new method and also illustrate its application in two data examples.

\section{Locally stationary long-range dependent processes} \label{sec2}
\def\theequation{2.\arabic{equation}}
\setcounter{equation}{0}

In order to develop a test for the presence of long-range dependence which can deal with different kinds of non-stationarity, a set-up is required which includes short-memory processes with a rather general time-varying first and second order structure and a reasonable long-range dependent extension. For this purpose, we  consider a triangular scheme $(\{X_{t,T}\}_{t=1,...,T})_{T \in \mathbb{N}}$ of locally stationary long-memory processes, which have an MA($\infty$) representation of the form
\begin{align}
\label{proc}
X_{t,T}=\mu(t/T)+\sum_{l=0}^{\infty}\psi_{t,T,l}Z_{t-l}, \quad t=1,\ldots,T,
\end{align}
where
\be
\label{summesqendl}
\sup_{T \in \mathbb{N}}\sup_{t \in \{ 1, \ldots, T\}} \sum_{l=0}^\infty \psi_{t,T,l}^2&<&\infty,
\ee
 $\mu:[0,1] \rightarrow \mathbb{R}$  is a  ``smooth'' function and $\{Z_t \}_{t\in  \mathbb{Z} } $ are independent standard normal distributed random variables.
  The assumption of a normal distribution for the innovations
 is made to simplify the technical arguments in the proofs of our results [see Section \ref{appendix}] and can be replaced by  the existence of
  moments of all order of the random variables $Z_{t}$ - see Remark \ref{prozess nicht-gauss}  for more details. Note
  also that the  random variables
  $Z_t$ have been standardized to have variance $1$. Alternatively, one could normalize by
 $\psi_{t,T,1} =1$  and allow for an additional parameter in the variance.
 For the coefficients $\psi_{t,T,l}$ and the function $\mu$ in the expansion (\ref{proc}) we make the following additional assumptions.

\begin{assumption}
\label{ass1}
Let $(\{X_{t,T}\}_{t=1,...,T})_{T \in \mathbb{N}}$ denote a sequence of stochastic processes  which have an MA($\infty$) representation of the form \eqref{proc} satisfying \eqref{summesqendl},
where $\mu$ is twice continuously differentiable. Furthermore, we assume that the following conditions are satisfied:
\begin{itemize}
\item[1)]
There exist  twice continuously differentiable   functions   $\psi_l: [0,1] \rightarrow \mathbb{R}$ ($l\in  \mathbb{Z}$)  such that the conditions
\be
\label{apprbed}
\sup_{t = 1, \ldots,T}\big |\psi_{t,T,l}-\psi_l(t/T) \big | &\leq& C T^{-1} I(l)^{D-1} \quad \forall l \in \mathbb{N} \\
\psi_l(u) &=&a(u)I(l)^{d_0(u)-1}+O(I(l)^{D-2}) \label{apprpsi}
\ee
are satisfied  uniformly with respect to $u \in [0,1]$ as $l \rightarrow \infty$ , where $I(x):=|x| \cdot 1_{\{x\not= 0\}} +1_{\{x=0\}}$  and
$
D=\sup_{u \in [0,1]}d_0(u)<1/2.
$
Moreover, the functions  $a:[0,1] \rightarrow \mathbb{R}$, $d_0:[0,1] \rightarrow [0,1/2) $
in \eqref{apprpsi}  are
twice continuously differentiable.
\item[2)] The time varying spectral density $f:[0,1]\times [-\pi,\pi]\rightarrow \mathbb{R}_{0}^{+}$
\be
\label{tvspectraldensity}
f(u,\lambda):=\frac{1}{2\pi}\Bigl|\sum_{l=0}^\infty \psi_l(u)\exp(-i \lambda l)\Bigr|^2
\ee
can be represented as
\begin{eqnarray}
f(u,\lambda)=  |1-e^{i\lambda}|^{-2d_0(u)} g(u,\lambda),
\label{spec}
\end{eqnarray}
where the function $g$ defined by
\begin{eqnarray} \label{gfkt}
g(u,\lambda):=\frac{1}{2\pi}\big|1+\sum_{j=1}^\infty a_{j,0}(u)\exp(-i \lambda j)\big|^{-2}
\end{eqnarray}
is twice continuously differentiable  (note that the identities \eqref{tvspectraldensity} and \eqref{spec} define the coefficients $a_{j,0}(u)$).
\item[3)] There exists a constant $C \in \mathbb{R}^{+}$, which is independent of $u$ and $\lambda$, such that for $l \not= 0$ the conditions
\be
\sup_{u \in (0,1)}|\psi_l'(u)| &\leq & C\log |l||l|^{D-1},
\sup_{u \in (0,1)}|\psi_l''(u)| \leq   C\log^2 |l||l|^{D-1},
\label{2.1c} \\
\sup_{u \in (0,1)} \big|\frac{\partial}{\partial u} f(u,\lambda)\big| &\leq&  C |\log(\lambda)||\lambda|^{-2D},
\sup_{u \in (0,1)} \big|\frac{\partial^2}{\partial u^2} f(u,\lambda)\big| \leq  C \log^2(\lambda)|\lambda|^{-2D}
\notag
 \ee
are satisfied for all $\lambda \in [-\pi,\pi]$.
\end{itemize}
\end{assumption}

Similar locally stationary long-range dependent models have been investigated by \cite{beran2009}, \cite{palole2010} and \cite{vonsachs2011}
 and \cite{wuzhou2014}. It is also worthwhile to mention that in general \eqref{apprpsi} does not imply \eqref{spec}  and \eqref{gfkt}  and vice versa
 conditions \eqref{spec}  and \eqref{gfkt} do not imply \eqref{apprpsi}.  {Therefore, none of the conditions \eqref{apprpsi}, \eqref{spec} or \eqref{gfkt} can be omitted in Assumption \ref{ass1}.}
Note  also  hat in contrast to the standard framework of local stationarity introduced by \cite{dahlhaus1997} and extended to the long-memory case in \cite{palole2010}, condition \eqref{apprbed} is much weaker. For example, in contrast to these references the assumptions made here include tvFARIMA($p,d,q$)-models as well [see Theorem 2.2 in \cite{prevet2012}]. Moreover, we mention again that the assumption of Gaussianity is only imposed to simplify the technical arguments in the proofs  of our main results -  see  Remark \ref{prozess nicht-gauss} for more details.   The very specific form of the function $g$ in (\ref{gfkt}) implies that the process $\{X_{t,T}\}_{t=1, \ldots, T}$ can be locally  approximated by a  FARIMA($\infty,d,0$) process in the sense of \eqref{apprbed}. More precisely, we obtain with
\be \label{pot}
b_k(u)= \binom{k+d(u)-1}{k}\quad \text{and} \quad (\sum_{k=0}^\infty a_{k,0}(u)z^k)^{-1}= \sum_{k=0}^\infty a_{k,0}^{(-1)}(u) z^k
\ee
($ a_{0,0}=1$) the relation
\bea
\psi_l(u) = \sum_{k=0}^l a_{k,0}^{(-1)}(u) b_{l-k}(u)
\eea
between the approximating functions $\psi_l(u)$ and the time-varying AR-parameters [see the proof of Lemma 3.2 in \cite{koktaqqu} for more details]. The relation  \eqref{pot} can be used to calculate the coefficients $a^{-1}_{k,0}(u)$ from the functions $a_{k,0}(u)$ , i.e.
$$a^{(-1)}_{0,0}(u)= {1 \over a_{0,0}(u)} ,~
a^{(-1)}_{1,0}(u)= -{a_{1,0}(u) \over  a_{0,0}^2(u)} ,  ~~\ldots $$

In order to further visualize  some properties of these kinds of locally stationary long-memory models we introduce for every fixed $u \in [0,1]$ the stationary process
\begin{eqnarray*}
X_t(u):=\mu(u)+\sum_{l=0}^\infty \psi_l(u) Z_{t-l}.
\end{eqnarray*}
One can show that condition (\ref{apprpsi}) implies the existence of bounded functions $y_i:[0,1] \rightarrow \boldsymbol{\mathbb{R}^{+}}$ $(i=1,2)$ such that the approximations
\be \label{cond1}
|\Cov(X_t(u),X_{t+k}(u))| \sim y_1(u)k^{2d_0(u)-1} \quad \text{as } k \rightarrow \infty
\ee
and
\begin{eqnarray} \label{cond2}
f(u,\lambda) \sim y_2(u)|\lambda|^{-2d_0(u)} \quad \text{ as } \lambda \rightarrow 0
\end{eqnarray}
hold [see \cite{palole2010} for details].
Consequently, the autocovariances  $\gamma_k(u,k)=\Cov(X_0(u),X_k(u))$ are  not absolutely summable if the function $a(u)$ in \eqref{apprpsi} is not vanishing, and in this case the time varying spectral density $f(u,\lambda)$ has a pole at $\lambda=0$ for any $u \in [0,1]$ for which $d_0(u)$ is positive.
Note that in general
the  statements \eqref{cond1} and \eqref{cond2}  are not equivalent [see \cite{yong1974} for a discussion of this problem in the stationary case]. \\
In the framework of these long-range dependent locally stationary processes we now investigate the null hypothesis that the time-varying ``long memory'' parameter $d_0(u)$ vanishes for all $u \in [0,1]$, i.e. there is no long-range dependence in the locally stationary process $X_{t,T}$. {The alternative is defined by the property that the function $d_0$ is nonnegative on the interval $[0,1]$ and positive on a subset of positive Lebesgue measure.} Since the function $d_0$ is continuous and non-negative we obtain that the hypotheses
\begin{eqnarray} \label{H_0}
\mbox{H}_0: d_0(u)=0  \hspace{.3cm}\forall u \in [0,1] \  \text{ vs.} \ \mbox{H}_1:
 d_0(u) \geq 0 \ \forall   u \in [0,1] \ \mbox {and} \
 d_0(u)>0 \hspace{.3cm} \text{for some }u \in [0,1]
\end{eqnarray}
are equivalent to
\begin{eqnarray}
\mbox{H}_0:  F=0 \quad \text{ vs. } \quad \mbox{H}_1: F >0,
\label{H_0 und H_1}
\end{eqnarray}
where the quantity $F$ is defined by
\begin{eqnarray}
F:=\int_{0}^{1}d_0(u) du.
\label{Ffkt}
\end{eqnarray}

In Section \ref{sec3} we will develop a nonparametric estimator of the function $d_0$ and
 the integral $F$.
   {   Roughly speaking,  the sample size $T$ is decomposed into $M$ blocks with length $N$ (i.e. $T=NM$), where
   $M$ is some positive integer. We define the corresponding midpoints in both the time and rescaled time domain by
   $t_j=N(j-1)+N/2$, $u_j=t_j/T$, respectively,  and calculate  an estimator $\hat d_N(u_j)$ of the long range dependence
   parameter  at the point $u_j$ on each of the $M$ blocks (for  the exact definition of the estimator see Section \ref{sec3}). The test statistic is then obtained as
\be \label{fhat1}
\hat F_T=\frac{1}{M}\sum_{j=1}^{M} \hat{d}_N(u_{j})
\ee
and could be considered as a Riemann sum of the integral $\int_0 \hat d_N(u)du$, which approximates the integral in \eqref{Ffkt}. }

\begin{rem} \label{remrev1} {\rm ({\it some boundary issues}) 
Note that  for each $u \in [0,1]$ the local long range dependence parameter $d_0(u)$ is a boundary point of the parameter
space $[0, 1/2)$ defined by  the two hypotheses in \eqref{H_0}. However, we will not use this property for the construction of the  estimates
  $\hat{d}_N(u_{j}) $ of the quantities $d_0(u_{j})$, which are aggregated in the statistic \eqref{fhat1}.
  For this purpose we consider a sequence of approximating
  tvFARIMA$(k,d,0)$  models, where  the parameter $k=k(T)$ converges to infinity  as the sample size increases
and   the corresponding long range
dependence parameters are allowed to vary in intervals of the form $[-\gamma_k, 1/2-\delta_k]$, where $(\gamma_k)_{k \in \en} $.
and  $(\delta_k)_{k \in \en} $ are positive sequences converging to $0$.
 We will prove in Theorem \ref{uniform} below
that this provides a uniformly consistent estimate of the function $d_0$ and that
an average of  these statistics provides a consistent  and asymptotically normal distributed estimate of the integral $F$ (see Theorem \ref{asymp teststat}  and
Theorem \ref{asymp teststatalt} below).
  As a consequence we obtain a consistent level-$\alpha$ test for the presence of long-range dependence in non-stationary time series  by rejecting the null hypothesis for large values of the estimator of $\hat F_T$. \\
  On a first glance these
properties are surprising because we use unconstrained  (i.e. potentially negative)  estimators of the long range dependence parameters in the approximating tvFARIMA models to estimate the non-negative function $d_0$, but the statements become clear
from the rather technical arguments given in the proofs of Section \ref{appendix}.
The situation is similar to the problem of testing the hypothesis $H_0: \mu=0$ versus $H_1: \mu>0$ for the mean of a
 sample of  i.i.d.  random variables $X_1,\ldots,X_n$. A test which rejects $H_0$, whenever $\sqrt{n} \overline{X}_n > \hat \sigma_n u_{1- \alpha}$ (here $\hat \sigma_n$ is an estimator of the variance and $u_{1- \alpha}$ the $(1- \alpha)$-quantile of the standard normal distribution) has asymptotic level $\alpha$ and is consistent. Moreover,
 in Section \ref{sec2A}  we consider  an  example of testing for a positive signal in a nonparametric regression model and
 demonstrate  that  the aggregation of local  statistics of the  type $ \overline{X}_n$  might have substantial  advantages compared
to the aggregation of  local estimators of the from  $\max \{  \overline{X}_n, 0\}$, which reflect the constraint $\mu \geq 0$ in its definition.
}
\end{rem}

\section{Testing for a positive nonparametric signal}  \label{sec2A}
\def\theequation{3.\arabic{equation}}
\setcounter{equation}{0}

In this section we provide some heuristic explanation for the phenomenon described in the previous paragraph,
which is also available to readers which are less familiar with the technical machinery used
for the asymptotic analysis of non-stationary long range dependent processes. We will also demonstrate that there are situations where more powerful tests can be
obtained by ignoring particular constraints in the estimation procedure. This situation occurs in particular if different estimators are aggregated as described in \eqref{fhat1}.

For this purpose  we consider
 the problem of testing the hypothesis of a   vanishing regression function against the alternative that the regression function is positive on the interval $[0,1]$ in the common nonparametric regression model
    $$
 Y_{ji} = \mu (t_j) + \varepsilon_{ji}; \qquad j=1,\dots,M;  \qquad i=1,\dots,N .
 $$
Here  $\varepsilon_{11},\dots, \varepsilon_{MN}$ are i.i.d. standard normal distributed (centered) random variables
 (this assumption is in fact not necessary but
 makes some of the following arguments much simpler), $t_j = t_{j,M} = j/M$ and $\mu$ is a smooth non-negative   Lipschitz continuous
 function  on the interval $[0,1]$. We are interested in testing the hypothesis
\be  \label{h0}
 H_0 : \mu (t) \equiv 0 \qquad  ~\mbox{ versus } \qquad H_1 : \mu (t) > 0 \qquad \qquad  \mbox {for all  } t \in [0,1]~.
\ee
Note that the alternative could also be considered on the set of all non-negative functions which are positive on  a subset of positive Lebesgue measure, say ${\cal U} \subset [0,1]$.
{As this generalization does not change any of the subsequent arguments (only integrals of the form $\int^1_0 \mu(t)dt$ and sums of the form $\frac {1}{M} \sum^M_{j=1} \mu (t_j)$ have to be replaced by $\int_{\mathcal{U}}\mu (t) dt$ and $\frac {1}{M \lambda (\mathcal{U})} \sum^M_{j=1} 1_{\mathcal{U}} (\frac {j}{M})\mu (t_j))$ (here $\lambda(\mathcal{U})$ denotes the Lebesgue measure of the set $\mathcal{U}$) we restrict ourselves to the case $\mathcal{U} = [0,1]$ for the sake of transparency.}

\medskip

{\bf 3.1 Tests based on unconstrained estimators:} The idea used in Section \ref{sec3} below
 for testing hypotheses of this type translates in the nonparametric regression model
  to the following procedure. We first define (unconstrained)
estimators for the values $\mu (t_j)$, that is
$
\hat \mu_j = \frac {1}{N} \sum^N_{i=1} Y_{ji},$
 ($ j = 1,\dots,M$),
and then consider the average
$$
T_M = \frac {1}{M} \sum^M_{j=1} \hat \mu_j ~=~ \frac {1}{MN} \sum^M_{j=1}  \sum^N_{i=1} Y_{ji} .
$$
Note that $
\mathcal{S}_M = \sqrt{N} ( T_M- \frac {1}{M} \sum^M_{j=1} \mu (t_j) ) = \frac {\sqrt{N}}{M} \sum^M_{j=1} (\hat \mu_j - \mu(t_j))
$
is a sum of independent identically distributed random variables with variance
$
\mbox{Var} (\mathcal{S}_M) = {1}/{M}.
$
Consequently, using a  central limit theorem for triangular arrays, shows that
$
\sqrt{M} \, \mathcal{S}_M \stackrel{\mathcal{D}}{\rightarrow} \mathcal{N} (0,1)
$ as $M \to \infty, N \to \infty$.
Moreover, since  $\mu$ is Lipschitz  continuous, this implies
$$
\sqrt{MN}\, \Bigl( T_M - \int^1_0 \mu(t) dt \Bigr) \stackrel{\mathcal{D}}{\rightarrow} \mathcal{N} (0,1),
$$
whenever $N=o(M) $.
Thus a consistent and  asymptotic level $\alpha$ test
for the hypothesis \eqref{h0} is obtained by rejecting the null hypothesis $H_0$, whenever
\be  \label{test1}
\sqrt{MN} T_M  > u_{1-\alpha},
\ee
 where $u_{1-\alpha}$ is the $(1-\alpha)$-quantile of the standard normal distribution.

\medskip

{\bf 3.2 Tests based on constrained estimators:}  Alternatively, and - on a first glance - more reasonable strategy is to use a constrained  estimator
  which addresses the boundary condition $\mu(t) \geq 0$.
This  gives
$$
\tilde \mu_j = \max \{ 0, \hat \mu_j \}
$$
as estimator for  $\mu(t_j)$, and  we obtain with the notation
$\tilde T_M = \frac {1}{M} \sum^M_{j=1} \tilde \mu_j$ the representation
\begin{eqnarray} \label{a1}
\mathcal{\tilde  S}_M &=& \sqrt{N} \ \Bigl(\tilde T_M - \frac {1}{M} \sum^M_{j=1} \mu (t_j) \Bigr)  
= \frac {\sqrt{N}}{M} \Bigl( \sum^M_{i=j} Z_j + \sum^M_{j=1} \delta_j \Bigr),
\end{eqnarray}
where $
Z_j =\max (0, \hat \mu_j) -  \mathbb{E} [ \max (0, \hat \mu_j) ]$, $
\delta_j = \mathbb{E} [ \max (0, \hat \mu_j) ] - \mu (t_j).$
Note that $\hat \mu_j  \sim {\cal N} (\mu(t_j) , 1/N)$, which yields
\begin{eqnarray}
\delta_j
&=& \frac {1}{\sqrt{2 \pi N}} \exp \Bigl(- \frac {N \mu^2 (t_j)}{2} \Bigr) - {\mu (t_j) \over\sqrt{\pi} }
\int^\infty_{\mu (t_j)\sqrt{N/2}} \exp(-z^2)dz\label{a1b} .
\end{eqnarray}
This  term is of order $o(1)$ (exponentially in $N$
and independent of $M$, provided  that $\mu (t) \ge c > 0$ on $[0,1]$). Note also that
\begin{eqnarray*}
\mathbb{E} [(\max (0, \hat \mu_j))^2] &=&
\mu^2 (t_j) + {1 \over N } +\frac { \mu (t_j) }{\sqrt{2 \pi N}} \exp \Bigl(- \frac {N \mu^2 (t_j)}{2} \Bigr)
- { 1+N \mu^2 (t_j)  \over N \sqrt{\pi}} \int^\infty_{\mu (t_j)\sqrt{N/2}} \exp(-z^2)dz~.
\end{eqnarray*}
 This gives for the variance of the random variable $Z_j$
$$
 \mathbb{E}  [Z_j^2]= \mbox{Var}  (\max (0, \hat \mu_j)) =
 \left\{ \begin{array}{cc}
 {1\over N}\big({1\over 2} - {1\over 2\pi} \big) &  \mbox{ if }  \mu(t_j)=0 \\
 {1\over N}(1+o(1)) &  \mbox{ if }  \mu(t_j) > 0
 \end{array}
 \right.
$$
 Ljapunoff's  central limit theorem now shows that
$
\sqrt{MN/\sigma_N^2} ( \tilde T_M - \frac {1}{M} \sum^M_{j=1} \mu (t_j)  - B_{M,N} )  \stackrel{\mathcal{D}}{\rightarrow} \mathcal{N} (0,1),
$
where $\sigma^2_N = N \mathbb{E} [Z_i^2]$  and
$$
B_{M,N}= {1\over M}  \sum_{i=j}^M \delta_j =
\left\{ \begin{array}{cc}
{1\over M}  \sum_{j=1}^M \delta_j  & \mbox{ if } \mu (t) > 0   \mbox{ for all } t \\
{1\over \sqrt{2\pi N} } & \mbox{ if } \mu (t) = 0   \mbox{ for all } t
\end{array}~.
\right.
$$
This  implies (observing the Lipschitz continuity of the regression function and $N=o(M)$)
$$
\sqrt{NM / \sigma^2_N } \Bigl( \tilde T_M -  \int^1_0 \mu (t) dt  - B_{M,N}  \Bigr) \stackrel{\mathcal{D}}{\rightarrow} \mathcal{N} (0,1).
$$
 {Note that the statistic is asymptotically normal distributed, although we average constrained estimators.}
Under the null hypothesis things are simplifying. In particular we obtain  $\sigma^2_{H_0}= \sigma^2_{N,H_0}=   N\mathbb{E}_{H_0}  [Z_i^2]= {1\over 2} - {1\over 2\pi} $
and  a test based on $ \tilde T_M $    rejects the null hypothesis  $H_0$, whenever
\be  \label{test2}
\tilde T_M  > {1 \over \sqrt{2 \pi N} } + \sigma_{H_0} {u_{1-\alpha}  \over \sqrt{MN} }=  {1 \over \sqrt{2 \pi N} } +
{  \sqrt{ {1\over 2}   - {1\over  2\pi} }}  {u_{1-\alpha}  \over \sqrt{MN} }.
\ee
This test has asymptotic level $\alpha$ and is consistent.
We conclude this section mentioning once again that the assumption of i.i.d. standard normal distributed errors was made to minimize the technical arguments. All statements remain true for arbitrary centered errors which have moments of order $4$. This observation is a simple consequence of the central limit theorem, and in the following finite sample comparison we actually use non-normal error distributions.

\medskip

{\bf 3.3 A comparison of the two tests:}
The use of different estimators for the quantities $\mu (t_i)$   yields to the different tests \eqref{test1}  and  \eqref{test2}  the hypotheses  in \eqref{h0}. Both test statistics have an asymptotic
normal distribution under the null hypothesis and the alternative.  A finite sample comparison
is given in  Table \ref{tab1a} where  we report  simulation results for the functions
\be
\label{m1} \mu_1 (t) &\equiv & 0, \\
\label{m2}  \mu_2 (t) &=& 0.1, \\
\label{m3}  \mu_3 (t) &=& 0.1 + 0.1t.
\ee
The  sample sizes are  $M=N=20$ and $M=N=50$ and we use $10 000$ simulation runs to estimate
the rejection probabilities of the tests \eqref{test1} and \eqref{test2}.
For the distribution of the errors distribution we use a $(\mathcal{X}^2_5-5)/\sqrt{10}$ distribution, in order to demonstrate that the previous findings do not depend on the assumption of normal distributed errors.
We observe that the test  \eqref{test1}  based on the statistic $T_M$  (which uses the unconstrained estimators of the regression function)
outperforms the method  \eqref{test2}  which uses the constrained estimators.

\begin{table}[h]
{
\begin{center}
   \begin{tabular}{|c|c|c|c|c|c|c||c|c|c|c|c|c|}
   \hline
     & \multicolumn{6}{|c||}{$M=N=20$}
       & \multicolumn{6}{|c|}{$M=N=50$} \\
       \hline
     model & \multicolumn{2}{|c|}{\eqref{m1}} & \multicolumn{2}{|c|}{\eqref{m2}} & \multicolumn{2}{|c||}{\eqref{m3}}& \multicolumn{2}{|c|}{\eqref{m1}} & \multicolumn{2}{|c|}{\eqref{m2}} & \multicolumn{2}{|c|}{\eqref{m3}}  \\
     \hline
     level & $5 \%$ & $10 \% $ & $5 \%$ & $10 \% $ & $5 \%$ & $10 \% $ & $5 \%$  & $10 \% $ & $5 \%$ & $10 \% $ & $5 \%$ & $10 \% $ \\
     \hline
     \hline
     { test \eqref{test1}} & 0.052 & 0.103 & 0.634 & 0.764 & 0.918 & 0.966 & 0.056 & 0.109 & 0.999 & 1.000 & 1.000 & 1.000 \\
     \hline
     { test  \eqref{test2}} & 0.070 & 0.118 & 0.576 & 0.684 & 0.873 & 0.926 & 0.065 & 0.112 & 0.997 & 0.998 & 1.000 & 1.000 \\
     \hline
\end{tabular}
\caption{\label{tab1a} {\it Simulated power of the tests  \eqref{test1}  and  \eqref{test2}  for the hypothesis \eqref{h0}
in model  \eqref{m1}  -  \eqref{m3}.}}
\end{center}
}
\end{table}

 We  can also give  a ``theoretical'' argument for the advantages of the unconstrained approach. Note that for a positive function $\mu$ the power of test \eqref{test1} is approximately given by
 \begin{equation} \label{b1}
 \mathbb{P}_{H_1} \Bigl( T_M > \frac {u_{1 -\alpha}}{\sqrt{MN}} \Bigr) \approx \Phi \Bigl(  \sqrt{MN} \int^1_0 \mu (t)
 dt - u_{1 -\alpha} \Bigr)~,
 \end{equation}
 where $\Phi$ denotes the distribution function of the standard normal distribution.
 This formula is remarkably precise. For example, if $N=M=20, \ \mu(t)=0.1$ we obtain for the power of the test \eqref{test1} $0.638$, while the result of the simulation is $0.643$.
Similarly, the power of the test \eqref{test2} is approximately given by
 \be \label{powertest2}
 \mathbb{P}_{H_1} \Bigl(\tilde T_M > \frac {1}{\sqrt{2 \pi N}} + \sigma_{H_0} \frac{u_{1- \alpha}}{\sqrt{MN}} \Bigr) \approx \Phi \Bigl({\sqrt{MN} \over \sigma_N}  \int^1_0 \mu (t)  dt + r_{N,M}  \Bigr),
\ee
 where the term $r_{N,M}$ is defined by
 $$
 r_{N,M} = \frac {\sqrt{NM}}{\sigma_N} B_{M,N} - \sqrt{\frac {M}{2 \pi \sigma_N^2}} - \frac {\sigma_{H_0}}{\sigma_N} u_{1 - \alpha}.
 $$
 Now, note that $\sigma_N^2=1+o(1)$ and that $B_{M,N}=o(1)$ of exponential order (uniformly) if $\mu(t) \geq c > 0$   for all $t \in [0,1]$
 as $M,N \to \infty$ Consequently, the term $ r_{N,M}$ will be negative for reasonable large $M,N$. It actually
diverges to $-\infty$, but at a lower rate as the dominating term ${ \sqrt{MN} \over \sigma}  \int^1_0 \mu (t)  dt $ in \eqref{powertest2}, which converges to $\infty$.
This  means that the test \eqref{test1}  based on unconstrained estimation is more powerful than the test \eqref{test2}, which uses constrained estimation. \\
 A similar argument for the superiority of the test \eqref{test1} based on the unconstrained estimators of the regression function
 can be given for local alternatives of the form $\mu_{M,N} (t) = c(t) / \sqrt{MN}$, where  $c:[0,1]\to \mathbb{R}$ is  a Lipschitz continuous function. 
 More precisely,  the asymptotic power of the tests (\ref{test1}) and (\ref{test2}) is given by
 $$
 \Phi \Big ( \int^1_0 c(t) dt - u_{1 - \alpha} \Big )
 $$
 and
 $$
 \Phi \Big ( \int^1_0 c(t) dt  \big / \sqrt{2 - 2/\pi} - u_{1 - \alpha} \Big ),
 $$
 respectively.
As $\sqrt{2- 2 \pi} \approx 1.1676 >1$, it follows that  the unconstrained  test \eqref{test1} also outperforms the test \eqref{test2} 
 under local alternatives.
 Exemplarily, we display in
  Table \ref{table2}  the power of the two tests under the local alternatives
$ \mu_{M,N} (t) =  ({1+t})/{\sqrt{MN}}.
$

\begin{table}[h]
\begin{center}
   \begin{tabular}{|c|c|c|c|c|c|c|c|c|c|c|}
   \hline
     & \multicolumn{2}{|c|}{$M=N=20$}
       & \multicolumn{2}{|c|}{$M=N=50$}
       &   \multicolumn{2}{|c|}{$M=N=100$}
       &   \multicolumn{2}{|c|}{$M=N=200$}
       &   \multicolumn{2}{|c|}{$M=N=500$} \\
       \hline
       level & $5 \% $ & $10 \% $ & $5 \% $ & $10 \% $ & $5 \% $ & $10 \% $ & $5 \% $ & $10 \% $ & $5 \% $ & $10 \% $ \\
       \hline
 test       \eqref{test1} & 0.458 &  0.600 & 0.444 & 0.591 & 0.442 & 0.588 & 0.437 & 0.587 & 0.451 & 0.588 \\
     test   \eqref{test2} & 0.420 & 0.536 & 0.393 & 0.516 & 0.382 & 0.508 & 0.378 & 0.506 & 0.379 & 0.509 \\
       \hline
\end{tabular}
\caption { \label{table2} \it  Simulated power  of the tests  \eqref{test1}  and  \eqref{test2}  under local alternatives in model \eqref{m3}.}
\end{center}
\end{table}

In the following section we will use a similar approach based on averages of unconstrained estimates of the function $d_0(\cdot)$ in  sequence
of approximating tvFARIMA models. The proofs in  Section \ref{appendix} show that this approach provides a consistent and
asymptotic level $\alpha$ test for the hypotheses \eqref{H_0 und H_1}.

\section{Testing short- versus long-memory}  \label{sec3}
\def\theequation{4.\arabic{equation}}
\setcounter{equation}{0}
In order to estimate the integral $F$ we use a sequence  of semi-parametric models approximating the processes $\{X_t(u) \}_{t \in \mathbb{Z}}$ with time varying spectral density (\ref{spec}) and proceed in several steps. First we  choose an  increasing sequence $k=k(T) \in \mathbb{N}$, which diverges 'slowly' to infinity as the sample size $T$ grows, and fit a  tvFARIMA($k$,$d$,0) model to the data. To be precise, we consider  a locally stationary long-memory model with time varying spectral density $f:[0,1]\times [-\pi,\pi]\rightarrow \mathbb{R}_{0}^{+}$ defined by
\be \label{lmmodel}
f_{\theta_k(u)}(\lambda)=|1-\exp(i \lambda)|^{-2d(u)}   g_k(u,\lambda) , 
\ee
where
\begin{eqnarray*}
 g_k(u,\lambda)= \frac{1}{2\pi} |1+\sum_{j=1}^k a_j(u) \exp(-i\lambda j)|^{-2}
\end{eqnarray*}
and  $\theta_k=(d,a_{1},\ldots ,a_{k}):[0,1] \rightarrow \mathbb{R}^{k+1}$ is a vector valued function.
 We emphasize again that  $k=k(T)$ depends on the sample size and refer to Assumption \ref{ass2} for the precise condtions regarding its
growth rate.
We then estimate the function $\theta_k(u)$ by a localized Whittle-estimator,  that is
\begin{eqnarray}
\hat{\theta}_{N,k}(u) &=& \arg \min_{\theta_k \in \Theta_{\lfloor uT\rfloor/T,k}} \mathcal{L}_{N,k}^{\hat \mu}(\theta_k, u),
\label{Whittleestimator}
\end{eqnarray}
where
\begin{eqnarray}
\mathcal{L}_{N,k}^{\hat \mu}(\theta_k, u):= \frac{1}{4\pi} \int_{-\pi}^{\pi} \Big( \log(f_{\theta_k}(\lambda))+ \frac{I^{\hat \mu}_N(u,\lambda)}{f_{\theta_k}(\lambda)}  \, \Big) d\lambda
\label{whittle}
\end{eqnarray}
denotes the (local) Whittle likelihood [see \cite{optimallength} or \cite{dahlpolo2009}] and for each $u \in [0,1]$
the parameter space $\Theta_{u,k} \subset \mathbb{R}^{k+1}$ is a compact set which will be specified in Assumption \ref{ass2}.
In  \eqref{Whittleestimator} and (\ref{whittle}) the quantity
\begin{eqnarray}
I_N^{\hat \mu}(u,\lambda)&:=&\Big| \frac{1}{\sqrt{2\pi N}}\sum_{p=0}^{N-1}\Big[X_{\lfloor uT\rfloor-N/2+1+p,T}-\hat \mu(\lfloor uT\rfloor-N/2+1+p,T)\Big] e^{-i p\lambda} \Big|^2,
\label{per}
\end{eqnarray}
denotes  the mean-corrected local  periodogram, $N$ is an even window-length which is 'small' compared to $T$
 and  $\hat{\mu}$ is an asymptotically unbiased estimator of the mean function $\mu:[0,1] \rightarrow \mathbb{R}$,
see \cite{dahlhaus1997}. Here and throughout this paper we use the convention
$X_{j,T}=0$ for $j \not\in \{1,...,T\}$. We finally obtain an estimator $\hat d_N(u)$ for the time-varying long-memory parameter by taking the first component of the $(k+1)$ dimensional vector $\hat \theta_{N,k} (u)$ defined in (\ref{Whittleestimator}). We emphasize that the tvFARIMA models are only used to define the estimator  $\hat d_N(u)$
as the solution of the optimization problem \eqref{Whittleestimator}.  \\
It will be demonstrated in Theorem \ref{uniform}
below that -  provided that the ``true'' underlying process can be approximated by  tvFARIMA models -
this approach results in a uniformly consistent estimator of the time-varying long-memory parameter.
For this purpose we define $\theta_{0,k}(u):=(d_0(u),a_{1,0}(u),...,a_{k,0}(u))$ as the $(k+1)$ dimensional vector containing
the long memory parameter $d_0(u)$ and the first $k$ AR-parameter functions $a_{1,0}(u),...,a_{k,0}(u)$  of the approximating process $\{ X_t(u)\}_{t \in \mathbb{Z}}$ defined by the representation (\ref{spec}) and \eqref{gfkt}. Here and
throughout this paper, $A_{11}$ denotes the element in the position (1,1)  and $\|A \|_{sp}$ the spectral norm of the matrix $A=(a_{ij})_{i,j=1}^{k}$, respectively.
We state the following technical assumptions.

\begin{assumption} \label{ass2} {\rm
 Let $k=k(T)$ be a sequence converging to infinity for increasing sample size $T$ and let
$(\gamma_\ell)_{\ell \in \en}$ and  $(\delta_\ell)_{\ell \in \en}$ denote positive sequences in the interval $(0, \min\{ 1/4,1/2-D\} )$
such that
\begin{eqnarray*}
&&
\liminf_{T\to \infty} \gamma_{k(T)} \log T  > 0~, ~~\liminf_{k\to \infty} \delta_{k(T)} \log T > 0, \\
&& \lim_{T\to \infty} \gamma_{k(T)}   =  0~, ~~   \lim_{k\to \infty} \delta_{k(T)} = 0.
\end{eqnarray*}
For each $u \in [0,1]$ and $k \in \{k(T), T \in \mathbb{N}\}$ define  $\Theta_{u,k}=[-\gamma_k,1/2-\delta_k] \times \Theta_{u,k,1} \times  \ldots \times \Theta_{u,k,k}$, where
the constant $D$ is the same as in Assumption \ref{ass1}. For each $i=1, \ldots, k$
 $\Theta_{u,k,i}$ is a compact set with a finite number (independent of $u,k,i$)  of connected components with
 positive  Lebesgue measure. Let $\Theta_k$ denote the space of all four times continuously differentiable functions $\theta_k:[0,1] \rightarrow \mathbb{R}^{k+1}$ with $\theta_k(u)
  \in \Theta_{u,k}$  for all $u \in [0,1]$. If $\theta_k(u)$ and $\theta_k^{'}(u)$ are distinct elements of  $\Theta_{u,k}$, we assume that the set $\{\lambda: f_{\theta_k(u)}(\lambda) \not = f_{\theta^{'}_k(u) }(\lambda) \}$ has positive Lebesgue measure.   We assume that the following conditions hold for each $k \in \{k(T), T \in \mathbb{N}\}$:
  \begin{itemize}
\item[(i)] The functions $g_k$ in \eqref{lmmodel} are bounded from below by a positive constant (which is independent of $k$) and are four times continuously differentiable with respect to $\lambda$ and $u$, where all partial derivates of $g_k$  up to the order four are bounded with a constant independent of $k$.
\item[(ii)] For each $u \in [0,1]$ the parameter $\tilde \theta_{0,k}(u)=\arg \min_{\theta_k \in \Theta_{u,k}} \mathcal{L}_k(\theta_k,u)$ exists
and is uniquely determined, where
\begin{eqnarray*}
\mathcal{L}_k(\theta_k,u):= \frac{1}{4\pi} \int_{-\pi}^{\pi} \Big( \log(f_{\theta_k}(\lambda))+ \frac{f(u,\lambda)}{f_{\theta_k}(\lambda)} \, \Big) d\lambda.
\end{eqnarray*}
Moreover, for each $u \in [0,1]$ the vectors $\tilde \theta_{0,k}(u)$ and $\theta_{0,k}(u)$ are  interior points of $\Theta_{u,k}$.
\item[(iii)] Define
\begin{eqnarray}
\Gamma_k(\theta_k) &=& \frac{1}{4\pi}\int_{-\pi}^{\pi}f_{\theta_k}^2(\lambda) \nabla f_{\theta_k}^{-1}(\lambda)\nabla f_{\theta_k}^{-1}(\lambda)^T \,d\lambda,
\label{gamk}\\
V_k(\theta_k,u) &=& \frac{1}{4\pi}\int_{-\pi}^{\pi} f^2(u,\lambda) \nabla f_{\theta_k}^{-1}(\lambda)\nabla f_{\theta_k}^{-1}(\lambda)^T d \lambda,
\notag
\end{eqnarray}
[here $\nabla$ denotes the derivative with respect to the parameter-vector $\theta_k$], then the matrix $\Gamma_k(\theta_{0,k})$ is
non-singular  for every $u \in [0,1]$, $k \in \{k(T), T \in \mathbb{N}\}$, and
\begin{eqnarray}
\lim_{T \rightarrow \infty}\int_0^1[\Gamma_k^{-1}(\theta_{0,k}(u))]_{1,1}du \Big / \int_0^1[ \Gamma_k^{-1}(\theta_{0,k}(u))V_k(\theta_{0,k}(u),u)\Gamma_k^{-1}(\theta_{0,k}(u)) ]_{1,1}du  = 1
\label{lim1}
\end{eqnarray}
as $T \rightarrow \infty$. Furthermore, condition  \eqref{lim1} is also satisfied if the function $\theta_{0,k}(u)$ is replaced by any  sequence $\tilde \theta_{T} (u)$ such that $\sup_{u \in [0,1]}|\tilde \theta_{T}(u)-\theta_{0,k}(u)| \rightarrow 0$. For such a sequence we additionally assume that the condition
\bea
\lim_{T \rightarrow \infty}\int_0^1[\Gamma_k^{-1}(\theta_{0,k}(u))]_{1,1}du/\int_0^1[\Gamma_k^{-1}(\tilde \theta_{T}(u))]_{1,1}du = 1
\eea
is satisfied as $T \rightarrow \infty$.
\item[(iv)] Let $\Theta_{R,k}=\bigcup_{u \in [0,1]} \Theta_{u,k}$ be compact and
$$ \sup_{\theta_k \in \Theta_{R,k}}\|\Gamma_k^{-1}(\theta_k)\|_{sp}=O(k)~,~~\liminf_{T \rightarrow \infty}\int_0^1[ \Gamma_k^{-1}(\theta_{0,k}(u))]_{1,1}du \geq c >0.$$
\end{itemize}
}
\end{assumption}

In order do illustrate the construction of the sets $\Theta_{u,k,i}$ in Assumption \ref{ass2}, consider exemplarily the case where for some $\delta > 0$
 the polynomial $z \rightarrow 1+\sum_{j=1}^\infty a_{j,0}(u)z^j$ with the coefficients from \eqref{gfkt} is bounded away from zero inside the disc
 $D_\delta := \{ z: |z| \leq 1 + \delta \}$  (uniformly with respect to $u$). In this case the sets  $\Theta_{u,k,1} \times ... \times \Theta_{u,k,k}$
 can   be chosen as the intersection of the set $\{(\theta_{u,k,1},\ldots ,\theta_{u,k,k}) \in \er^k |
|1+\sum_{j=1}^k \theta_{u,k,j}z^j| > C_1 > 0 \phantom{.} \forall z \in D_\delta \}$ with
the set  $$\{ (a_1,...,a_k) \in \er^k: \text{ there exists  a sequence } (a_i)_{i > k} \text{ such that } (a_i)_{i\in \en} \in A_0 \}.
$$
Here the set $A_0$ is defined by
\bea
A_0 := && \Big \{ (a_i)_{i \in \en} ~  \Big | ~
\text{ the polynomial } p(z):= 1+\sum_{j=1}^\infty a_j z^j \text{ satisfies } |p(z)|> C_2 > 0  \\
  && \phantom{..} \text{ and } |p^{(l)}(z)| \leq C_3 \text{ for all } z \in D_\delta  \text{ and  } 0 \leq l \leq 4    \Big\},
\eea
 the constants $C_2, C_3$ are chosen such that $C_1< C_2$ and such that the sequence $(a_{j,0})_{j \in \en}$ is an inner point of the set $A_0$.  \\
Assumption (i) and (ii)  are rather standard in a semi-parametric locally stationary time series model [see for example \cite{optimallength} or \cite{dahlpolo2009}
among others]. Note that the number of parameters  $k$  grows  with increasing sample size in order to obtain a consistent estimate of the function $u \rightarrow d(u)$ in model  (\ref{tvspectraldensity}).
 The restriction  on the spectral norm in part (iv)  was verified for a large number of long-range dependent models
 by \cite{kokozkarate} [see equation (4.4) in this reference].   Note that these assumptions  solely depend on the "true" underlying model. \\
 On the other hand, an important step of our approach is  the approximation
of  the spectral density $f(u,\lambda)$ in \eqref{spec}  by the  truncated analogue
\begin{eqnarray*}
|1-e^{i \lambda}|^{-2d_0(u)} |1+\sum_{j=1}^k a_{j,0}(u) e^{-i\lambda j} |^{-2},
\end{eqnarray*}
and
the following  assumption guarantees that the corresponding  approximation error
converges to $0$ with  reasonable rate.
As a consequence it provides a link between  the growth rate of $k=k(T)$ and $N$
as the sample size $T$ increases.

\begin{assumption} \label{ass3}
Suppose that $N \rightarrow \infty$, $N \log(N) = o(T)$ and
\be
\sup_{u \in [0,1]}\sum_{j=k+1}^\infty | a_{j,0}(u)| &=&O(N^{-1+\epsilon})
\label{null 4}
\ee
for some $0 < \epsilon < 1/6$ as $T \rightarrow \infty$.
\end{assumption}
Note that
\be \label{neukr}
 f(u,\lambda)-f_{\theta_{0,k}(u)}(\lambda) = |1-e^{i \lambda}|^{-2d_0(u)}\Big ( \big |1+\sum_{j=1}^\infty a_{j,0}(u) e^{-i\lambda j} \big|^{-2}
 - \big |1+\sum_{j=1}^k a_{j,0}(u) e^{-i\lambda j} \big |^{-2} \Big)~,
\ee
and an application of  Lemma 2.4 in \cite{kreisspappol2011} to the second factor (corresponding to the
 "short memory" part)  shows
 that  condition \eqref{null 4} with $0 < \gamma_k < 1/2-D$ implies
\bea \label{null 4 implication}
\sup_{u \in [0,1]}\int_{-\pi}^\pi |1-e^{i\lambda}|^{-2\gamma_k} \big | f(u,\lambda)-f_{\theta_{0,k}(u)}(\lambda) \big | d\lambda &=&O(N^{-1+\epsilon}).
\eea
As a consequence  Assumption \ref{ass2} (iii) is rather intuitive, because the parametric model
 \eqref{lmmodel} can be considered as an approximation of the ``true'' model   defined in terms of the time varying spectral density \eqref{tvspectraldensity}.
We finally note that  condition \eqref{null 4} is satisfied for a large number of tvFARIMA($p,d,q$) models, because it can be shown by similar arguments as
in the proof of Theorem 2.2 in \cite{prevet2012} that
 the coefficients $a_{j,0}(u)$ are geometrically decaying. This yields $\sum_{j=k+1}^\infty \sup_u | a_{j,0}(u)| = O(q^k)$ for some $q \in (0,1)$ resulting in a
  logarithmic growth rate for $k$, which is in line with the findings of \cite{kokozkarate}. Similarly, one can include processes whose AR coefficients decay
   such that $\sum_{j=0}^\infty \sup_u |a_{j,0}(u)| j^r < \infty$ is satisfied for some $r \in \mathbb{N}_0$. In this case $k$ needs to grow at some
  specific polynomial rate.

Our first main result states a uniform convergence rate for the difference between $\hat{\theta}_{N,k}(u)$ and its true counterpart $\theta_{0,k}(u)$. As a consequence it implies that the estimator $\hat d_N$
obtained in the approximating models  is uniformly consistent for the (time
varying) long-range dependence parameter of the locally stationary process.

\begin{satz}
\label{uniform}
Let Assumption \ref{ass1}, \ref{ass2} and \ref{ass3} be satisfied and suppose that the estimator of the mean function $\mu$ satisfies
\begin{eqnarray}
N^{\epsilon} k^{3} \max_{t=1, \ldots, T}\big|\mu(t/T)-\hat{\mu}(t/T)\big| &=& o_p(1)
\label{convergencemu}
\end{eqnarray}
for some $0 < \epsilon < \min\{1/4-D/2,1/6 \}$. If $N^{5/2}/T^2 \rightarrow 0$ and $k^4\log^2(T)N^{-\epsilon/2} \rightarrow 0$, then
\begin{eqnarray} \label{uniform0}
\sup_{u\in [0,1] }  \big\|  \hat{\theta}_{N,k}(u)-  \theta_{0,k}(u)\big\|_2 = O_P\big (k^{3/2} N^{-1/2+\epsilon}+N^{\epsilon} k^{3/2} \max_{t=1, \ldots, T}\big|\mu(t/T)-\hat{\mu}(t/T)\big| \big ).
\end{eqnarray}
In particular
\begin{eqnarray*}
\sup_{u\in [0,1] } |  \hat{d}_{N}(u)-  d_{0}(u)| = O_P\big (k^{3/2} N^{-1/2+\epsilon}+N^{\epsilon} k^{3/2} \max_{t=1, \ldots, T}\big|\mu(t/T)-\hat{\mu}(t/T)\big| \big).
\end{eqnarray*}
\end{satz}

\begin{rem} \label{remcons} {\rm
It follows from the proof of Theorem \ref{mean} below that there exists an estimator $\hat \mu$ with
\bea
N^{1/2-D-\alpha}\max_{t=1, \ldots, T} \big|\mu(t/T)-\hat{\mu}(t/T)\big|=o_p(1)
\eea
for every  $\alpha > 0$. Under the addional assumption
\be
\sup_{u \in [0,1]}\sum_{j=k+1}^\infty | a_{j,0}(u)| &=&O(q^k)
\label{nullhol1}
\ee
for some $q \in (0,1)$  a logarithmic rate for the dimension $k$ of the tvFARIMA models can be used
such that assumption (\ref{convergencemu}) is satisfied
[for a  broad class of models, where the stronger condition \eqref{nullhol1} is in fact satisfied, we refer
to   the discussion following (\ref{neukr})].
}
\end{rem}

In order to obtain an estimator of the quantity  $F$ in \eqref{Ffkt} we assume without loss of generality that the sample size $T$ can be decomposed into $M$ blocks with length $N$ (i.e. $T=NM$), where  $M$ is some positive integer. We define the corresponding midpoints in both the time and rescaled time domain by  $t_j=N(j-1)+N/2$, $u_j=t_j/T$, respectively,  and calculate  $\hat d_N(u_j)$ on each of the $M$ blocks as described in the previous paragraph. The test statistic is then obtained as
\be \label{fhat}
\hat F_T=\frac{1}{M}\sum_{j=1}^{M} \hat{d}_N(u_{j}).
\ee


The following two results specify the asymptotic behaviour of the statistic $\hat F_T$ under the null hypothesis and alternative.

\begin{satz}
\label{asymp teststat}  Assume that the null hypothesis $H_0$ (of no long-range dependence) is true.
Let Assumptions \ref{ass1}, \ref{ass2} and \ref{ass3} be satisfied, define $W_T=[\int_0^1\Gamma_k^{-1}(\theta_{0,k}(u))du]_{1,1}$ and
suppose that the estimator $\hat \mu $ of  the mean function satisfies
\begin{eqnarray}
\max_{t=1, \ldots, T}\big|\mu(t/T)-\hat{\mu}(t/T)\big| &=&  O_p(N^{-1/2+\epsilon/2})
\label{null 2},\\
\max_{t=1, \ldots, T}\Big|\Big\{\mu\Big(\frac{t-1}{T}\Big)-\hat{\mu}\Big(\frac{t-1}{T}\Big)\Big\}-\Big\{{\mu}\Big(\frac{t}{T}\Big)-\hat{\mu}\Big(\frac{t}{T}\Big)\Big\}\Big| &=& O_p(N^{-1/2-2\epsilon}T^{-1/2})
\label{null 3},
\end{eqnarray}
where $\epsilon$ is the constant in  Assumption \ref{ass3} satisfying $0 < \epsilon<1/6$. Moreover, if the conditions
\begin{eqnarray}
k^6\sqrt{T}/N^{1-\epsilon}  \rightarrow 0 ,\hspace{.4cm}
k^4\log^2(T)/N^{\epsilon/2}  \rightarrow 0 ,\hspace{.4cm}
 k^2\log(T)/T^{1/6-\epsilon}  \rightarrow 0,&&\hspace{-.2cm}
k^2N^2/T^{\frac{3}{2}}\rightarrow 0
\notag
\end{eqnarray}
hold as $T \rightarrow \infty$, then we have
\begin{eqnarray}
&&\sqrt{T} \hat F_T /\sqrt{W_T}
\stackrel{D}{\rightarrow} \mathcal{N}(0, 1).
\label{lim}
\end{eqnarray}
\end{satz}

Note that $\hat F_T$ is an average of the estimates of the long-range dependence parameter in the approximating
 tvFARIMA model. By Assumption \ref{ass2} the point $0$ is an interior point of the canonical projection of the parameter space $\Theta_{u,k}$ onto the first component, which motivates the asymptotic normality obtained in Theorem \ref{asymp teststat}. More precisely, we show in Section \ref{appendix} that the leading term in the stochastic expansion of $\hat F_T$ is given by
\begin{eqnarray*}
-\frac{1}{M}\sum_{j=1}^{M}\frac{1}{4\pi} \int_{-\pi}^{\pi} \big(I^{\mu}_N(u_j,\lambda)-f_{\theta_{0,k}(u_j)}(\lambda) \big) [\Gamma_k^{-1}(\theta_{0,k}(u_j))\nabla f^{-1}_{\theta_{0,k}(u_j)}(\lambda)]_{1} \,d\lambda,
\end{eqnarray*}
where $[a]_1$ denotes the first element of the $(k+1)$ dimensional vector $a$.
Asymptotic normality follows because the individual terms in this sum are asymptotically independent (see Section \ref{appendix} for details and Section \ref{sec2A} for a similar result in a simplified model).
\begin{satz}
\label{asymp teststatalt}
Assume that the  alternative $H_1$  of long-range dependence is true.
Let Assumptions \ref{ass1}, \ref{ass2} and \ref{ass3} be satisfied and
suppose that the estimator $\hat \mu $ of  the mean function satisfies
\begin{eqnarray}
N^{\epsilon} k^{3} \max_{t=1, \ldots, T}\big|\mu(t/T)-\hat{\mu}(t/T)\big| &=& o_p(1)
\label{alt 2},
\end{eqnarray}
where $\epsilon$ is the constant in Assumption \ref{ass3} satisfying $0 <  \epsilon < \min\{1/4-D/2,1/6\}$. Moreover, if  the conditions
\begin{eqnarray}
k^{6}/N^{1-2\epsilon} \rightarrow 0 , \hspace{.3cm}
k^4\log^2(T)/N^{\epsilon/2} \rightarrow 0 , \hspace{.3cm}
k^{4}/N^{1-2D-2\epsilon} \rightarrow 0 , \hspace{.3cm}
 k^2 N^{5/2}/T^2 \hspace{-.1cm}&\rightarrow& 0
\notag
\end{eqnarray}
are satisfied as $T \rightarrow \infty$,  then we have
\begin{eqnarray*}
&& \hat F_T  \stackrel{P}{\rightarrow} F >0.
\end{eqnarray*}
\end{satz}

\begin{rem} \label{remneu} {\rm ({\it more transparent conditions})
If assumption \eqref{nullhol1} is satisfied, more transparent conditions for  Theorem \ref{uniform}, \ref{asymp teststat} and \ref{asymp teststatalt} can be given. To be precise assume that  \eqref{nullhol1} holds
for some $q \in (0,1)$  and choose
$$k =\lfloor-a {\log T \over \log q}\rfloor
$$
 for some $a \in (1/2,1)$. If $D < 1/6$, then it follows by straightforward but tedious calculations that Theorem    \ref{uniform}, \ref{asymp teststat} and \ref{asymp teststatalt}
hold for $N=T^\beta$ with any $\beta$ satisfying $a < \beta < \min \{ {6 \over 5}a, {3\over 4} \} $ (note that this conditions provides a further restriction for the choice of the constant  $a$).  Similarly, if $1/6 \le D < 1/2$
the results hold, whenever  $a < \beta < \min \{ {4a \over 3 + 2D}, {3\over 4} \} $.
}
\end{rem}

\begin{rem}
\label{prozess nicht-gauss}
{\rm  ({\it the non-Gaussian case) }
It is worthwhile to mention that in most of articles cited in this paper the assumption of Gaussianity for the innovation process in (\ref{ass1}) is required. In the present case this assumption is not necessary  and is only imposed here to simplify technical arguments in the proof of Theorem \ref{zgws}. This observation is a consequence of method of proof used in Section \ref{appendix}. In fact, asymptotic normality is established by the method of moments showing that all cumulants of the statistic under consideration converge to those of a normal distribution. In the definition of all cumulants one needs the existence of all moments of $Z_i$ (which is obviously true in the Gaussian case).
 The main simplification under the assumption of Gaussianity consists in the fact that one does not have to work with partitions including cumulants of any possible order. The extension to non Gaussian innovations does  not change the main argument in  the proofs, but the calculations become
substantially more complicated, and the details are omitted for the sake of brevity. \\
As a consequence all results of this section remain true as long as the innovations are independent with all moments existing, mean zero and $\E(Z_{t}^2)=\sigma^2(t/T)$ for some twice continuously differentiable function $\sigma: [0,1] \rightarrow \mathbbm{R}$.
To be more precise, in order to address for non-Gaussian innovations
the variance $V_T$ in Theorem \ref{zgws} (which is one of the main ingredients for the proofs in Section \ref{appendix}) has to be replaced by
\bea
V_{T,general}=V_T+\frac{1}{TM}\sum_{j=1}^M \kappa_4(u_j)/\sigma^4(u_j) \Big( \int_{-\pi}^\pi f(u_j,\lambda) \phi_T(u_j, \lambda) d \lambda \Big)^2,
\eea
where  $V_T$ is defined in \eqref{var} and $\kappa_4(u)$ denotes the fourth cumulant of the innovations, i.e. $\kappa(t/T)=\E(Z_{t}^4)-3\sigma^4(t/T)$ for all $t=1,\ldots, T$. In the proof of Theorem \ref{asymp teststat} we apply this result with  $\phi_T(u_j, \lambda)= (4\pi)^{-1}[\Gamma_k^{-1}(\theta_{0,k}(u_j))\nabla f^{-1}_{\theta_{0,k}(u_j)}(\lambda)]_{1}$ . Consequently, we obtain that  in the non-Gaussian case the asymptotic normality in Theorem \ref{asymp teststat} holds, where the matrix
$W_T=[\int_0^1\Gamma_k^{-1}(\theta_{0,k}(u))du]_{1,1}$ has to be replaced by
\be \label{asymgen}
W_{T,general}=W_T+\frac{1}{TM}\sum_{j=1}^M \kappa_4(u_j)/\sigma^4(u_j) \Big( \int_{-\pi}^\pi f(u_j,\lambda) \phi_T(u_j, \lambda) d \lambda \Big)^2.
\ee
Thus, under the null hypothesis it follows that
\begin{eqnarray}
&&{\sqrt{T} \hat F_T \over \sqrt{W_{T,general}}}
\stackrel{D}{\rightarrow} \mathcal{N}(0, 1).
\label{lim1}
\end{eqnarray}
 }
\end{rem}
\medskip

\begin{rem} \label{gentlest}  {\rm ({\it the final test}) 
Note that the first term $W_T$ in \eqref{asymgen}  can be   consistently estimated by
\bea
\hat W_T=\Bigl[ \frac{1}{M} \sum_{j=1}^{M} \Gamma_k^{-1}(\hat \theta_{N,k}(u_j)) \Bigr]_{11}.
\eea
This gives as an estimator for $V_{T,{\rm general}}$ the statistic
\bea
\hat W_{T,general}=\hat W_T+\frac{1}{ M}\sum_{j=1}^M \hat \kappa_4(u_j)/\hat \sigma^4(u_j) \Big( \int_{-\pi}^\pi f_{\hat \theta_{N,k}(u_j)}(\lambda) [\Gamma_k^{-1}(\hat \theta_{N,k}(u_j))\nabla f^{-1}_{\hat \theta_{N,k}(u_j)}(\lambda)]_{1} d \lambda \Big)^2,
\eea
where $\hat \sigma(u_j)$ and $\hat \kappa(u_j)$ are obtained by calculating the empirical second and fourth moment $\hat \mu_{2,Z}(u_j)$, $\hat \mu_{4,Z}(u_j)$ of the residuals
\bea
Z_{t,res}=X_{t,T}-\sum_{i=2}^k [\hat \theta_{N,k}(u_j)]_i X_{t-i+1,T}, \quad t=t_j-N/2+k+1,t_j-N/2+k+2,...,t_j+N/2,
\eea
and setting $\hat \sigma^2(u_j)=\hat \mu_{2,Z}(u_j)$, $\hat \kappa(u_j)=\hat \mu_{4,Z}(u_j)-3\hat \mu^2_{2,Z}(u_j)$.
Since  $\hat W_{T,general} / W_{T,general} \xrightarrow{P} 1$,  an asymptotic level $\alpha$-test is obtained from (4.18) by rejecting the null hypothesis (\ref{H_0 und H_1}),
whenever
\begin{eqnarray} \label{testrule}
&&\sqrt{T}\hat F_T /\sqrt{\hat W_{T,general}}  \geq u_{1-\alpha},
\end{eqnarray}
where $u_{1-\alpha}$ denotes the $(1-\alpha)$-quantile of the standard normal distribution (in the Gaussian case $\hat W_{T,general}$ can be replaced by
$\hat W_{T}$). It then follows from  Remark  \ref{prozess nicht-gauss} and Theorem \ref{asymp teststatalt} that for any estimator of the mean function $\mu$ satisfying \eqref{null 2}, \eqref{null 3} and
\eqref{alt 2}, the test, which rejects $H_0$ whenever \eqref{testrule} is satisfied, is a consistent level-$\alpha$ test for the null hypothesis stated in \eqref{H_0 und H_1}.
The finite sample properties of this resulting test are investigated in Section \ref{sec4}.
}
\end{rem}

\bigskip

A popular estimate of the mean function is given by
the  the local-window estimator
 \begin{eqnarray}
 \hat \mu_L (u)=\frac{1}{L}\sum_{p=0}^{L-1} X_{\lfloor uT\rfloor-L/2+1+p,T}~,
\label{muhat}
 \end{eqnarray}
 where $L$ is a window-length which does not necessarily coincide with the corresponding parameter in the calculation of the local periodogram. Note that also $I_N^{\hat \mu}(u,\lambda)$ is an asymptotically unbiased estimator for $f(u,\lambda)$ if $N \rightarrow \infty$ and $N/T \rightarrow 0$.
The final result of this section shows that this estimator
satisfies the assumptions of Theorem \ref{asymp teststat} and \ref{asymp teststatalt} if $L$ grows at a 'slightly' faster rate than $N$. This means, it can be used in the asymptotic level $\alpha$ test defined by \eqref{testrule}.
\begin{satz}  \label{mean} ~~\\
\vspace{-.8cm}
\begin{itemize}
\item[a)] Suppose that the assumptions of Theorem \ref{asymp teststat} hold and
additionally $N^{1+4\epsilon}/L^{1-\delta} \rightarrow 0$ and \linebreak $L^{5/2-\delta}/T^{3/2} \rightarrow 0$ are satisfied for some $\delta>0$, where $\epsilon >0$ denotes the constant in Theorem \ref{asymp teststat}. Then the local-window estimator $ \hat \mu_L$ defined in
\eqref{muhat}  satisfies \eqref{null 2} and \eqref{null 3}.
\item[b)] Suppose the assumptions of Theorem \ref{asymp teststatalt} hold. If additionally $N^{\epsilon} k^{5}/L^{1/2-D-\delta} \rightarrow 0$ and $L^{5/2-D} /T^2 \rightarrow 0$ for some $0< \delta < 1/2-D-\epsilon$ (with the constant  $\epsilon$ from Theorem \ref{asymp teststatalt}), then the local-window estimator  $ \hat \mu_L$ defined in  \eqref{muhat}  satisfies \eqref{alt 2}.
\end{itemize}
\end{satz}

\begin{rem} {\rm   ({\it parametric models})
Analogues of Theorem \ref{asymp teststat} and \ref{asymp teststatalt} can be obtained in a parametric framework. To be  precise, assume that the approximating processes
$\{ X_t(u)\}_{t \in  \mathbb{Z}}$   has  a time varying spectral density of the form \eqref{lmmodel}, where $k$ is fixed and known. In this case it is not necessary
that the dimension $k$ is increasing with the sample size  $T$ and Assumption \ref{ass2}(iii) and \ref{ass3} are  not required. All other stated assumptions are rather standard in this framework of a semi-parametric locally stationary time series model [see for example \cite{optimallength} or \cite{dahlpolo2009} among others].
With these modifications Theorem \ref{asymp teststat} and \ref{asymp teststatalt}  remain valid and as  a consequence we obtain an alternative test to the likelihood ratio test proposed in  \cite{Yau2012}, which operates in the spectral domain and
can be used  for  more general null hypotheses as
considered by these authors.
}
\end{rem}

\begin{rem} \label{lok alt}
{\rm ({\it local alternatives})
Theorem \ref{asymp teststat} remains valid under local alternatives converging to the null hypothesis at a rate $\sqrt{T/k}$. To be precise let $d_{0,T}(u)= a(u)\sqrt{W_{T, \rm {general}}/T}$ where $a:[0,1] \rightarrow [0,\infty)$ is a twice continuously differentiable function such that $\int_{0}^{1}a(u) \,du>0$. Then it follows by similar arguments as
given  in the proof of Theorem \ref{asymp teststat}, that
\begin{eqnarray*}
\sqrt{T} \Big(\frac{\hat F_T-\int_0^1 a(u) \,du}{\sqrt{W_{T, \rm {general}}}}\Big) \stackrel{D}{\rightarrow} \mathcal{N}(0,1)
\end{eqnarray*}
(note that $W_{T, \rm {general}}=O(k)$ due to Assumption \ref{ass2}(iv) and that $W_T$ does not depend on the long-memory parameter function $d_0$). This indicates that (asymptotically) the power of the test (\ref{testrule}) is increasing with $\int_{0}^{1}a(u) \,du$, which can also be observed in the simulation study presented in the following Section.
}
\end{rem}

\begin{rem} \label{weakendcondition}
{\rm ({\it some technical comments})
 The (uniform) smoothness conditions stated in Assumption \ref{ass1} are commonly made in the literature
 [see for example  \cite{palole2010}] and are also required in the present context
 to obtain the uniform consistency of the estimator for the function $d_0$.
 However, it is worthwhile to mention  that the asymptotic properties of the
 proposed test can also be derived  under weaker  assumptions. To be more precise, Theorem \ref{asymp teststat}  remains valid if the conditions on
 the function $\psi_l(u)$ and its derivatives stated in Assumption \ref{ass1} are replaced by
 \bea
 |\psi_l'(u)| &\leq & C(u) \log |l||l|^{D-1},
 |\psi_l''(u)| \leq   C(u)\log^2 |l||l|^{D-1},\\
\big|\frac{\partial}{\partial u} f(u,\lambda)\big| &\leq&  C(u) |\log(\lambda)||\lambda|^{-2D},
 \big|\frac{\partial^2}{\partial u^2} f(u,\lambda)\big| \leq  C (u) \log^2(\lambda)|\lambda|^{-2D}
 \eea
 for all $\lambda \in [-\pi,\pi]$ and $u\in (0,1)$.
Here   $C: (0,1) \rightarrow \er$ denotes a function such that    $\int_0^1 |C(u)|^p du< \infty$ for all $p \in \en$. The proof of this statement can be performed by
similar arguments as given in  the proof of Theorem \ref{asymp teststat} with additional technical arguments for the more delicate
estimates of the  error terms. \\
Moreover, we conjecture that, the conditions can be further weakened such that the function
$C$ is only  integrable up to a specific order. A detailed verification of such a statement, however, is an open problem and far beyond the scope of
the present paper.
}
\end{rem}

\section{Finite sample properties}  \label{sec4}
\def\theequation{5.\arabic{equation}}
\setcounter{equation}{0}

The application of the test \eqref{testrule} requires the choice of several parameters. Based on an extensive numerical investigation
we recommend the following rules. For the choice of the parameter
$L$ in the local window estimate $\hat \mu_L $  of the mean function [for a precise definition see  (\ref{muhat})] we use
$L = N^{1.05}$. Because the procedure is based on a  sequence  of approximating tvFARIMA$(k,d,0)$-processes the choice of the order
$k$ is essential, and we suggest  the AIC criterion  for this purpose, that is
\begin{eqnarray}
\hat k = \arg \min_k \frac{1}{T} \sum_{j=1}^{T/2} \Big( \log(h_{\hat \theta_{k,s}}(\lambda_j))+ \frac{I^{\hat \mu }(\lambda_j)}{h_{\hat \theta_{k,s}}(\lambda_j)}\Big)+ \frac{k+1}{T},
\label{AIC}
\end{eqnarray}
where $\lambda_j= 2\pi j /T \hspace{.1cm} (j=1, \ldots, T)$, and $h_{\hat \theta_{k,s}}(\lambda)$ is the estimated spectral density of a stationary FARIMA$(k,d,0)$ process and $I^{\hat \mu_L }(\lambda)$ is the mean-corrected periodogram given by
\begin{eqnarray*}
I^{\hat \mu_L }(\lambda)&:=&\Big| \frac{1}{\sqrt{2\pi N}}\sum_{t=1}^{T}\Big[X_{t,T}-\hat \mu_L(t/T)\Big] e^{-i t\lambda} \Big|^2.
\end{eqnarray*}

Note that we choose the same order $k$ for each of the $M$ blocks. An alternative choice is to use tvFARIMA models of different order for each block. In our numerical experiments we investigated both methods and we observed substantial advantages for the rule  $(\ref{AIC})$ (the results of this comparison are not displayed for the sake of brevity). Because this approach also has additional computational advantages we recommend to choose the same approximating tvFARIMA(k,d,0) model for all blocks.
Finally, the performance of the test depends on the choice of $N$, and this dependency will be carefully investigated in the following discussion.
\subsection{Simulation of level and power}
\label{Simulation results}

\begin{figure}
\begin{center}
\includegraphics[width=0.3\textwidth]{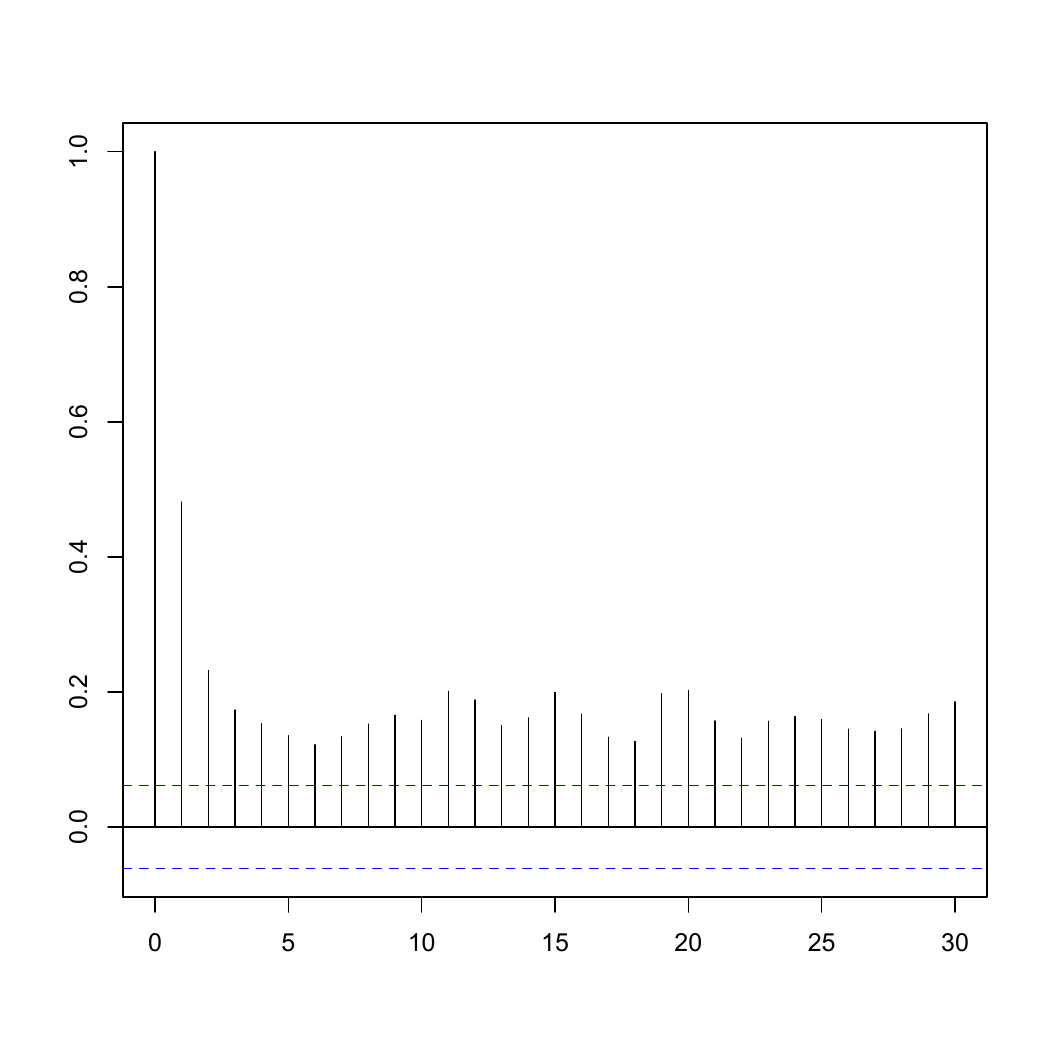}~~
\includegraphics[width=0.3\textwidth]{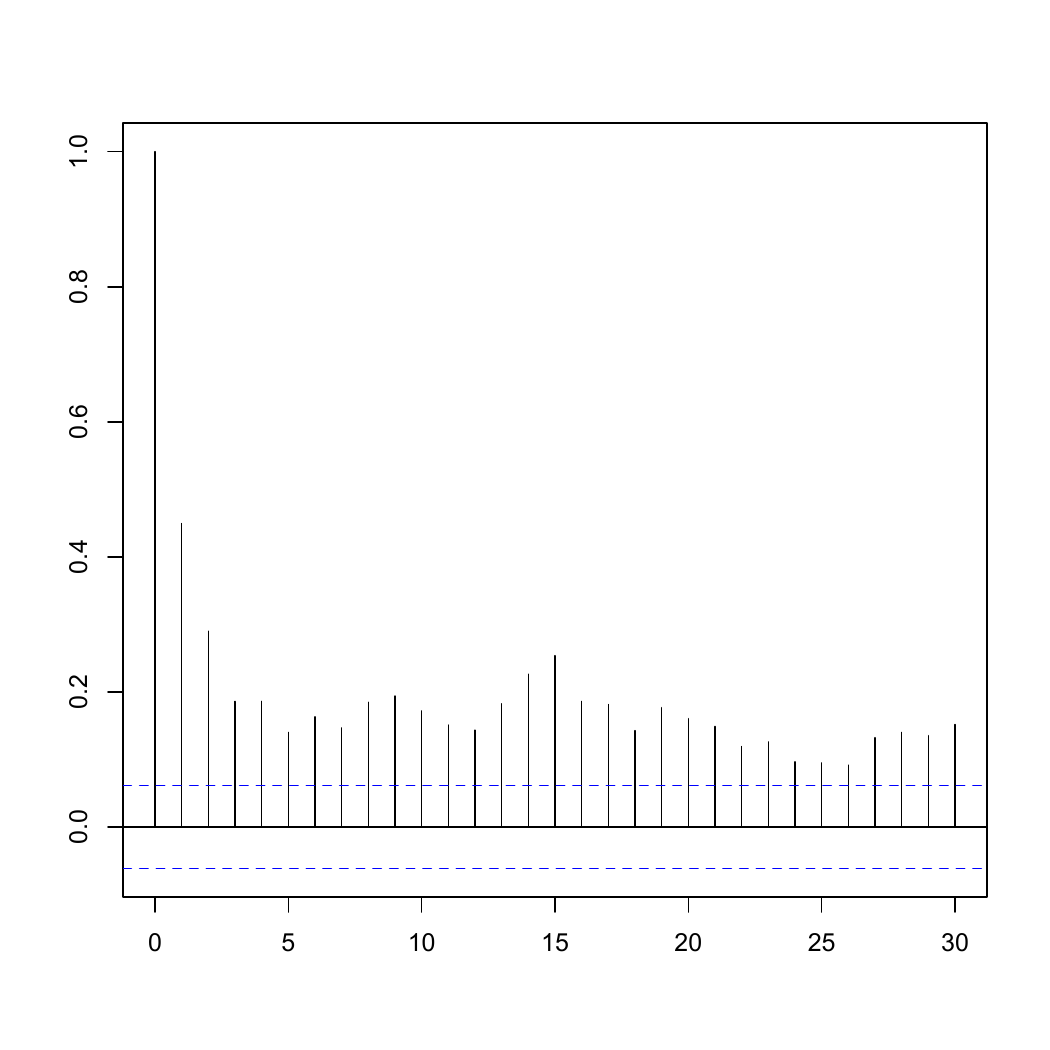}~~
\includegraphics[width=0.3\textwidth]{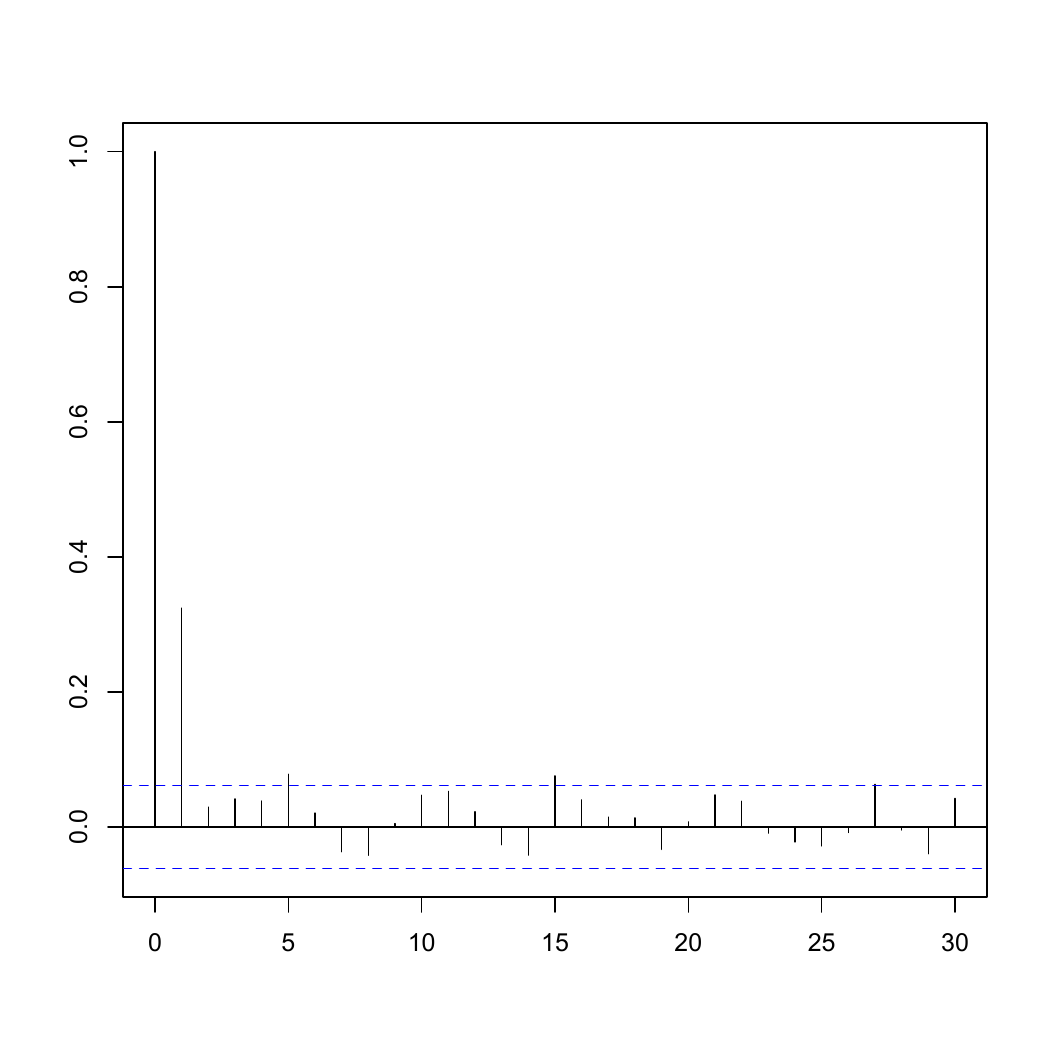}~~
  \caption{ \label{ACF_H_0} {
  \it Sample autocovariance functions of model (\ref{AR-Prozess}) with mean function (\ref{cont}) (left panel), (\ref{jump}) (middle panel) and of model (\ref{MA1}) (right panel).
  The sample size is T=1024.} }
\end{center}
\end{figure}

\begin{figure}
\begin{center}
\includegraphics[width=0.3\textwidth]{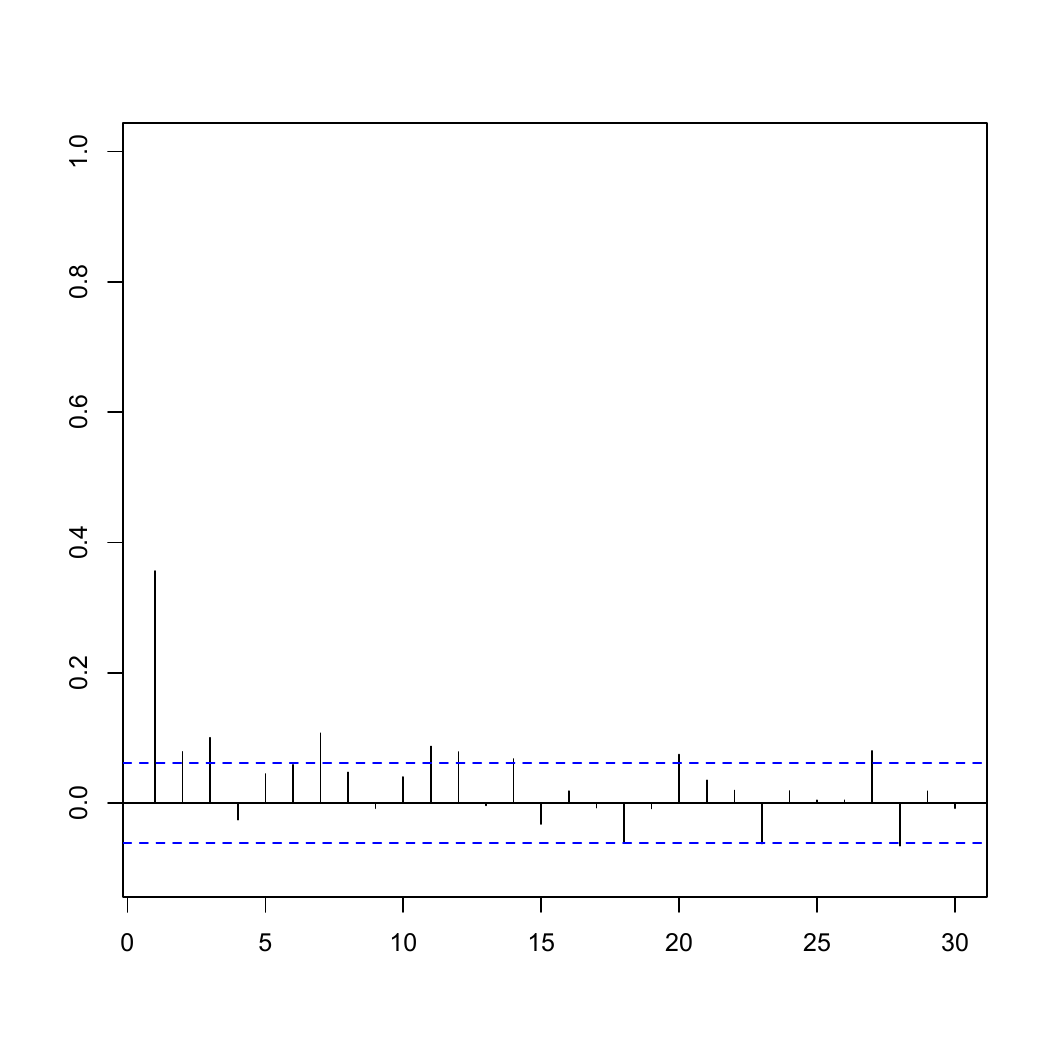}~~
\includegraphics[width=0.3\textwidth]{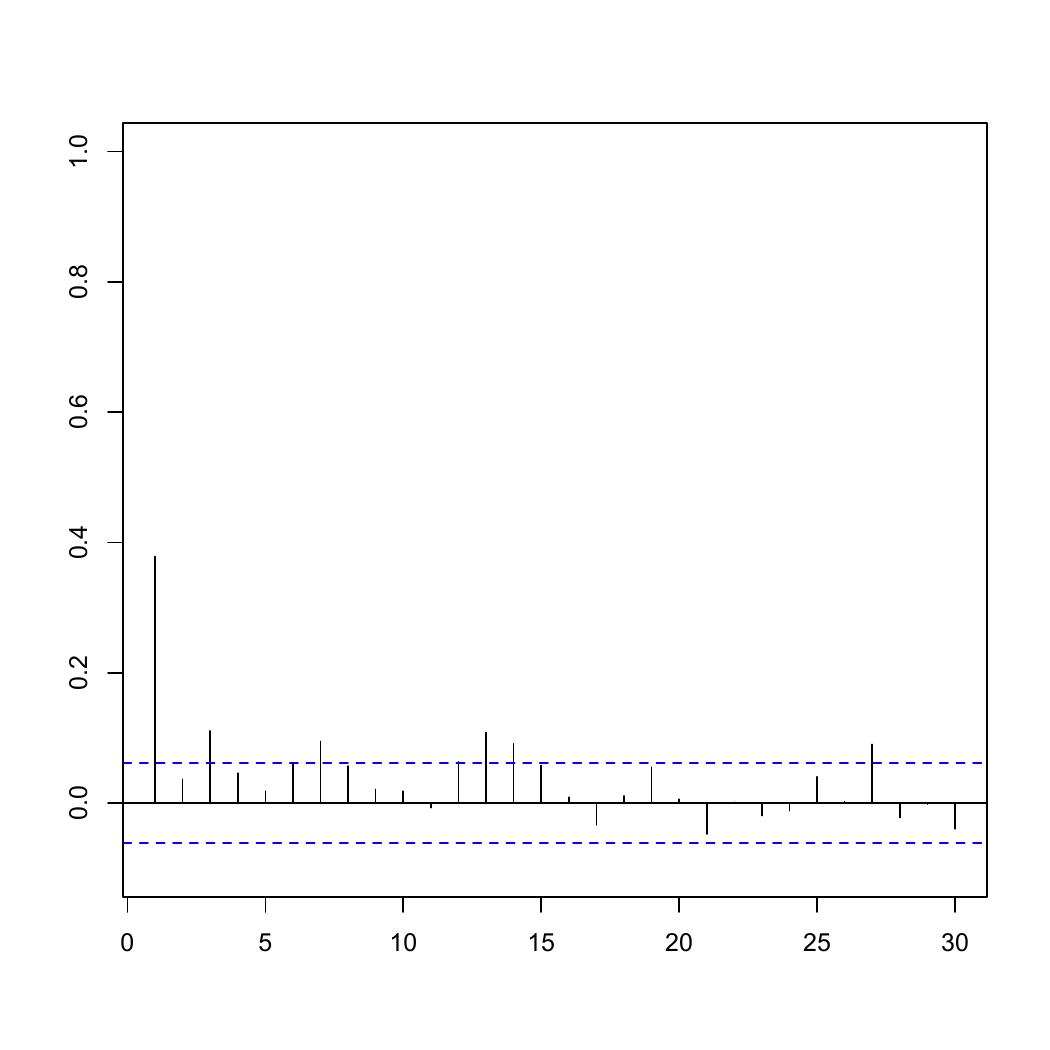}~~
\includegraphics[width=0.3\textwidth]{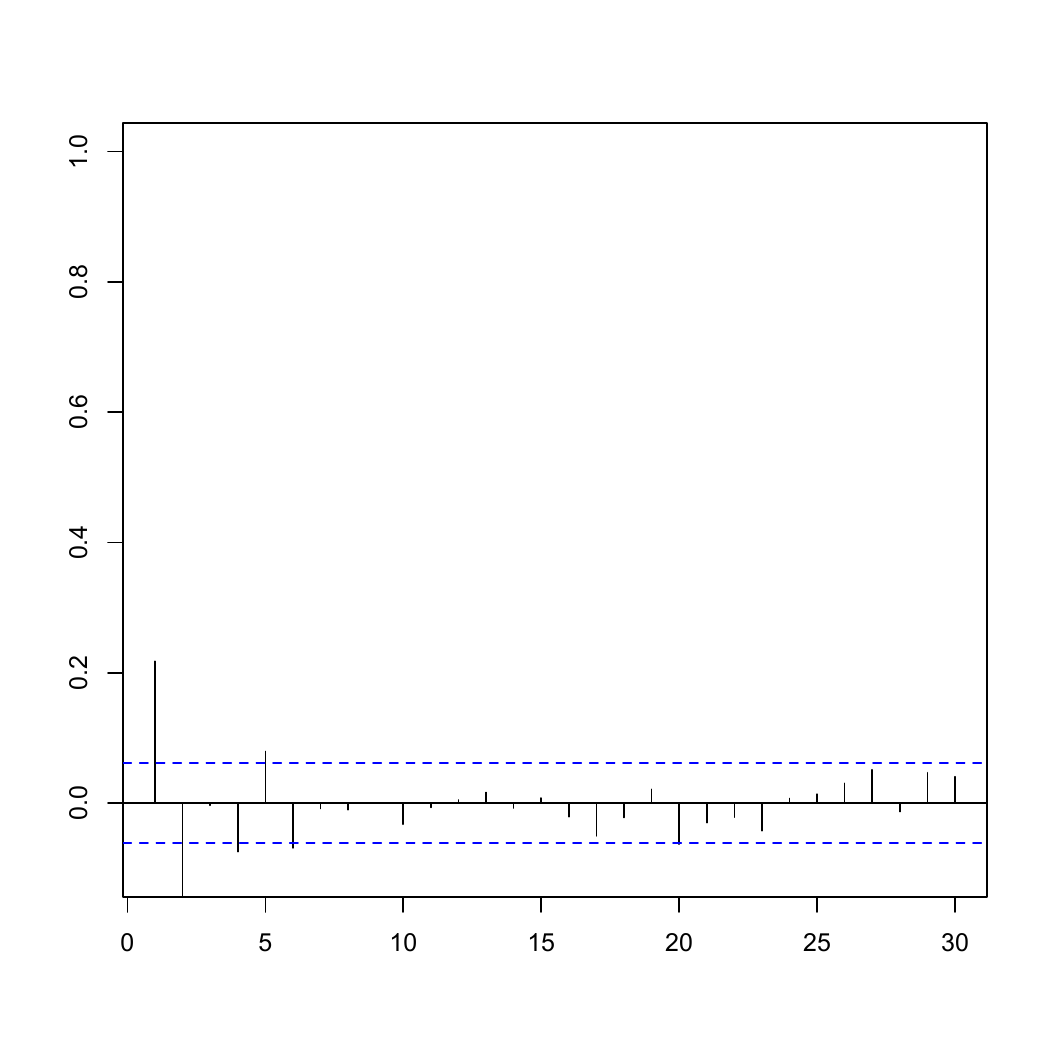}~~
  \caption{ \label{PACF_H_0} {
  \it Sample partial autocovariance functions of model (\ref{AR-Prozess}) with mean function (\ref{cont}) (left panel), (\ref{jump}) (middle panel) and of model (\ref{MA1}) (right panel).
  The sample size is T=1024.} }
\end{center}
\end{figure}

All results presented in this section are based on $1000$ simulation runs, and we begin with an investigation of the approximation of the nominal level  of the test \eqref{testrule} considering  three examples.
The first example is given by a location model with a tvAR(1)-process, that is
\begin{eqnarray}
X_{t,T}&=& \mu_i(t/T)+ Y_{t,T}, \quad t=1, \ldots, T,
\label{AR-Prozess}
\end{eqnarray}
where
\begin{eqnarray}
Y_{t,T}= 0.6\frac{t}{T}Y_{t-1,T}+Z_{t}, \quad t=1, \ldots, T.
\label{tvAR1}
\end{eqnarray}
The innovations $\{ Z_{t}\}_{t=1, \ldots, T}$ in (\ref{tvAR1}) are either
i.i.d. standard normal or
i.i.d. chi-square distributed normalized such that $E[Z_i]=0, Var(Z_i)=1$, i.e.    $Z_i \sim (\chi^2_5-5)/\sqrt{10}$.
Two cases are investigated for the mean function representing a smooth change and abrupt change in the mean effect, i.e.
\begin{eqnarray}
\mu_1(t/T) &=&1.2 \frac{t}{T}\label{cont}, \\
\mu_2(t/T) &=& \left \{
\begin{array}{cl}
0.65& \mbox{for } t=1, \ldots, T/2\\
1.3& \mbox{for } t=T/2+1, \ldots T.
\end{array}
\right.
\label{jump}
\end{eqnarray}
The mean function (\ref{jump}) is not smooth and used to investigate the impact of a violation of the assumptions in the procedure.
Our third example consists of a tvMA(1)-process given by
\begin{eqnarray}
X_{t,T}=Z_{t}+0.55\sin\Big(\pi \frac{t}{T}\Big)Z_{t-1}, \quad t=1, \ldots, T,
\label{MA1}
\end{eqnarray}
where   $\{Z_{t}\}_{t=1, \ldots, T}$ is again a sequence of i.i.d. normal or chi-square distributed random variables normalized to have mean 0 and variance 1.
Figure \ref{ACF_H_0} and \ref{PACF_H_0} show the sample autocovariance and the sample partial autocovariance functions of 1024 observations generated by the  models (\ref{cont}), (\ref{jump}) and (\ref{MA1}), respectively, from which it is clearly visible that the mean functions in (\ref{cont}) and  (\ref{jump}) are causing a long-memory type behaviour.
In Table \ref{tab1}, we show for these models the simulated level of the test \eqref{testrule}
for various choices of $N$. We observe in model (\ref{AR-Prozess}) and (\ref{MA1}) a  reasonable approximation of the nominal level whenever
$M=T/N \approx 4$ and the sample size $T$  is larger or equal than $512$. Here the results are similar for normal and chi-square distributed innovations. On the other hand in  model \eqref{AR-Prozess} with
mean function \eqref{jump} the assumptions of the asymptotic theory are violated and the situation is different. For moderate sample sizes the specification $M=T/N \approx 4$ yields to an overestimation of the nominal level.
Moreover, the approximation of the nominal level  becomes worse with increasing sample size. We conjecture that the performance of the test
could be improved by using estimators addressing the problem of jumps in the mean function.

In order to investigate the power of the test \eqref{testrule} and to compare it with the  competing procedures
proposed by \cite{horvarth2006}, \cite{baekpipi} and \cite{Yau2012}, we simulated data from  a tvFARIMA($1,d,0)$-process
\begin{eqnarray}
(1+0.2\frac{t}{T} B)(1-B)^{d(t/T)}X_{t,T}= Z_{t}, \quad t=1, \ldots, T,
\label{power1}
\end{eqnarray}
and a  tvFARIMA($0,d,1)$-process
\begin{eqnarray}
(1-B)^{d(t/T)}X_{t,T}=(1-0.35\frac{t}{T} B)Z_{t},\quad t=1, \ldots, T,
\label{power2}
\end{eqnarray}
where $B$ is the backshift operator that is $B^{j}X_{t,T}:=X_{t-j,T}$. In both cases the long-memory function is given by $d(t/T)=0.1+0.3t/T$.
Because all competing procedures are designed to detect stationary long-range dependent alternatives, we also  simulated data from a stationary FARIMA(1,$d$,1)-process
\begin{eqnarray}
(1+0.25B)(1-B)^{0.1}X_{T}= (1-0.3B)Z_{t}, \quad t=1, \ldots, T.
\label{power3}
\end{eqnarray}
The corresponding results for the new test \eqref{testrule} and its competitors are  presented  in Table \ref{tab2}-\ref{tab4}. In Table \ref{tab2} and \ref{tab2x}  we show the simulated power in model (\ref{power1}) for (standardized) normal and chi-square distributed innovations. We do not observe substantial differences in the power of the new test under different distributional assumptions and for this reason Table \ref{tab3} and \ref{tab4}  only contain results for normal distributed innovations. In the first column the rejection probabilities of the new test are displayed and  we observe a reasonable power in all models under consideration. Interestingly, the differences in power between the tvFARIMA($1,d,0)$ and the tvFARIMA($0,d,1)$-model are rather small (see second column in Table \ref{tab2} and \ref{tab3}). The results in Table \ref{tab4} show a loss in power, which corresponds to intuition because the ``average'' long-memory effect in model (\ref{power3}) is 0.1, while it is $\int_{0}^{1}(0.1+0.3u)\,du= 0.25$ in model (\ref{power1}) and (\ref{power2}) [see also Remark \ref{lok alt} and the discussion at the end of this section].
In order to compare the new test with existing approaches we next investigate the performance of the procedures  proposed by \cite{horvarth2006}, \cite{baekpipi} and \cite{Yau2012}, which are designed for a
test of the null hypothesis  ''the process has the  short memory property
with a structural break in the mean" against the alternative "the process is stationary and has the  long memory property".
The third columns of Table  \ref{tab2}-\ref{tab4} show the power of  the test in  \cite{baekpipi}, which also operates in the spectral domain.
These authors estimate
the change in the mean with a break point estimator and remove
this mean effect (which is responsible for the observed local stationarity) from the time series. Then they calculated the local Whittle estimator introduced by \cite{rob95} for the self similarity parameter and reject the null hypothesis for large value of this estimate.
Note that the calculation of the local Whittle estimator requires the specification of the number of ``low frequencies'' and we used $m=\sqrt{T}$ as  \cite{baekpipi} suggested in their simulation study.
We observe that the new test  \eqref{testrule}  yields larger power than the procedure of \cite{baekpipi}  in nearly all cases under consideration. This improvement becomes more substantial with increasing sample size.\\
Next we study the performance of the procedure proposed by   \cite{horvarth2006} in models \eqref{power1}-\eqref{power3}. These authors use a CUSUM statistic to construct an estimator, say $\hat k^*$, for a (possible) change point $k^*$ in a time series.
Then two CUSUM statistics are computed for the first $\hat k^*$ elements of the time series and the remaining ones, respectively.
The test statistic is given by the maximum of those two.
For the choice of the bandwidth function we use $q(n)=15\log_{10}(n)$ as  suggested by these authors in
Section 3 of their article. The results are depicted in the fourth columns of Table \ref{tab2}-\ref{tab4} and demonstrate that this test is not able to detect long-range dependence in both the stationary and locally stationary case. These findings coincide with the results of \cite{baekpipi} who also remarked that the test in \cite{horvarth2006} has very low power against long-range dependence alternatives.\\
The method proposed by \cite{Yau2012} consists of a parametric
 likelihood ratio test  assuming two (not necessarily equal) ARMA$(p,q)$ models before and after
the breakpoint of the mean function. Their method requires a specification of  the order of these two models and
we used ARMA$(1,1)$-models under the null hypothesis and a FARIMA$(1,d,1)$ model under the alternative hypothesis.
The corresponding results for this test are depicted in the fifth columns of Table \ref{tab2}-\ref{tab4} corresponding to non-stationary and stationary
long-range dependent alternatives, respectively. We observe that in these cases the new test   \eqref{testrule}  outperforms
the test proposed in \cite{Yau2012} if the sample size is larger than $512$ and that both tests have similar power for sample size $256$ (see the fifth column of Table \ref{tab2} and \ref{tab3}). On the other hand, in the case of the long-range dependent stationary alternative (\ref{power3}) the test of \cite{Yau2012} yields slightly better rejection probabilities than the new test (\ref{testrule}) for smaller sample sizes  while we observe advantages of the proposed test in this paper for sample sizes $512$ and $1024$. These results are remarkable, because the test of \cite{Yau2012} is especially designed to detect stationary alternatives of FARIMA($1,d,1$) type, but the new semi-parametric test still yields an improvement in many cases.

Finally, as it was pointed out by a reviewer, it is also of interest to systematically investigate the power of the test \eqref{testrule} as a function of the quantity $F=\int_{0}^{1}d(u)\,du$. The arguments in Remark \ref{lok alt} indicate that the power is increasing with $F$, and we will now investigate if these properties can also be observed in finite samples. For this purpose we simulated data from the tvFARIMA(0,d,1)-process in (\ref{power2}) with different choices for the long-memory function $d$:
\begin{eqnarray}
d_1(t/T)&=&1/8 \quad,\quad d_2(t/T)=t/4T,
\label{long1}\\
d_3(t/T)&=&\left \{
\begin{array}{cl}
0 & \mbox{for } 1 \leq t \leq T/3\\
3/8 & \mbox{for }T/3 <  t  \leq 2T/3 \\
0 & \mbox{for } 2T/3 < t \leq T,\\
\end{array}
\right.
\label{long2}\\
d_4(t/T)&=&0.3 \quad, \quad d_5(t/T)=1.8t/T(1-t/T).
\label{long3}
\end{eqnarray}
For the functions $d_1, d_2,$ and $d_3$ the quantity $F=\int_{0}^{1}d(u)\,du$ is given by $1/8$ while $F=3/10$ for $d_4$ and $d_5$. The corresponding results are shown in Table \ref{tab5}. We mainly discuss the case $M=4$ (because it yields to the best approximation of the nominal level) and mention that the interpretation of the results for other choice of $M$ is very similar. For a fixed $F=1/8$ we do not observe substantial differences between the functions $d_1$ and $d_2$ in the case $M=4$, while the function $d_3$ yields to a larger power. This observation can be explained by the fact that the integral in (\ref{Ffkt}) is approximated by a Riemann sum $\frac{1}{M}\sum_{j=1}^{M} d(u_{j})$ at points $u_j=\frac{j-1}{M}+\frac{1}{2M}$. Consider exemplarily the case $M=4$ (which is recommended, because it yields to a good approximation of the nominal level). While for the function $d_1(u)=1/8$ all estimates roughly yield the same contribution of size $1/8$, we observe that for the function $d_3$ two points (namely $u_2$ and $u_3$) yield a contribution of size $3/8$ and the other points $u_1,u_4$ yield the value $d_3(u_j)=0 \hspace{.1cm} (j=1,4)$. Nevertheless the total contribution in this case is $3/16$, while it is only $1/8$ for $d_1$. This explains the improvement in power observed for the function $d_3$. We expect that these advantages vanish asymptotically, because the approximation of $F=\int_{0}^{1}d(u)\,du$ by its Riemann sum becomes more accurate with increasing $M$. Finally, a comparison of columns 1-3 (corresponding to the case $F=1/8$ with columns 4-5 in Table \ref{tab5}
(corresponding to the case $d=3/10$) shows that the monotonicity of the power as a function of the integral $F=\int_{0}^{1}d(u)\,du$ can also be observed in samples of realistic size.
\begin{table}
\begin{center}
{\scriptsize
\begin{tabular}{|ccc|cc|cc|cc|cc|cc|cc|}
\hline
 & $$ & $$ &   \multicolumn{6}{c|}{$Z_t \sim \mathcal{N}(0,1)$}  &   \multicolumn{6}{c|}{$Z_t \sim (\chi^2_5-5)/\sqrt{10}$}     \\
\hline
 & $$ & $$ &   \multicolumn{2}{c|}{(\ref{AR-Prozess}), (\ref{cont})}  &   \multicolumn{2}{c|}{(\ref{AR-Prozess}),(\ref{jump})}     & \multicolumn{2}{c|}{(\ref{MA1})} &   \multicolumn{2}{c|}{(\ref{AR-Prozess}), (\ref{cont})}&\multicolumn{2}{c|}{(\ref{AR-Prozess}),(\ref{jump})}& \multicolumn{2}{c|}{(\ref{MA1})}    \\
\hline
$T$&$N$&$M$& $5 \%$ & $10 \%$   & $5 \%$ & $10 \%$   & $5 \%$ & $10 \%$& $5 \%$ & $10 \%$    & $5 \%$ & $10 \%$  & $5 \%$ & $10 \%$ \\
\hline
256&		64  	&	4  	&.090	&.128	&.094	&.145	&	.085	& .122	&.094 	& .142	&.100	 &.162	&084		&.118\\
256&		32  	&	8	&.151	&.228	&.165	&.255	&	.182	&.261	&.218	&.319	&.249	 &.335	&.187	&.258\\
512&		128	&	4	&.061	&.095	&.070	&.114	&	.069	&.099	&.066  	& .100	&.062	 &.098	&.068	&.090\\
512&		64  	&	8 	&.089	&.130	&.089	&.126	&	.081	& .107	&.086  	& .144	&.102	 &.140	&.074	&.118\\
1024&	256	&	4 	&.046	&.072	&.077	&.119	&	.069	& .106	&.042	& .076	&.080	&.126	 &.080	&.114\\
1024&	128	&	8	&.059	&.087	&.061	&.088	&	.064	&.093	&.058	& .082 	&.090	&.124	 &.066	&.106\\
2048&	512	&	4 	&.048	&.090	&.094	&.148	&	.074	&.122	&.048	& .078	&.116	&.154	 &.086	&.116\\
2048&	256	&	8 	&.026	&.034	&.026	&.058	&	.062	&.084	&.020	& .026	&.040	&.068	 &.046	&.074\\
4096& 	1024&	4	&.056	&.094 	&.164	&.248	&	.076	&.112	&.052	& .098	&.196	&.264	 &.085	&.127\\
4096&	512	&	8 	&.014	&.030	&.026	&.056	&	.060	&.080	&.026	& .044	&.046	&.062	 &.058	&.090\\
\hline
\end{tabular}}
\caption{\textit{\label{tab1} \small{Simulated level of the test \eqref{testrule} for different processes and choices of T,N and M. }}}
\end{center}
\end{table}

\begin{table}
\begin{center}
{\scriptsize
\begin{tabular}{|ccc|cc|cc|cc|cc|cc|}
\hline
 & $$ & $$ &   \multicolumn{2}{c|}{\eqref{testrule} }  &   \multicolumn{2}{c|}{Baek/Pipiras}     & \multicolumn{2}{c|}{Berkes et. al} & \multicolumn{2}{c|}{Yau/Davis}  \\
\hline
$T$&$N$&$M$& $5 \%$ & $10 \%$   & $5 \%$ & $10 \%$   & $5 \%$ & $10 \%$  & $5 \%$ & $10\%$  \\ \hline
256&		64  &4  	&0.288	&0.354	&	0.248	&0.330	&0.037	&0.080	&0.250	&0.306	\\
256&		32  &	8	&	0.290&0.436	&			&		&		&		&		&		\\
512&		128&	4	&	0.530&0.590	&	0.356	&0.468	&0.006	&0.041	&0.182	&0.226	\\
512&		64  &	8 	&	0.348&0.458	&			&		&		&		&		&		\\
1024&	256&	4 	&	0.746&0.770	&	0.562	&0.656	&0.026	&0.102	&0.204	&0.267	\\
1024&	128&	8	&	0.412&0.512	&			&		&		&		&		&	\\
2048& 	512&4	&	0.882&0.900	&0.724		&0.816	&0.152	&0.222	&0.376	&0.452	\\
2048&	256&8	&	0.625&0.683	&			&		&		&		&		&	\\
4096&	1024&4	&	0.974&0.978	&0.892		&0.928	&0.318	&0.460	&0.740	&0.782	\\
4096& 	512&8	&	0.892&0.910	&			&		&		&		&		&	\\
\hline
\end{tabular}}
\caption{\textit{\label{tab2} \small{Rejection frequencies of the test \eqref{testrule} and three competing procedures under the alternative (\ref{power1}) for different choices of T,N and M. The innovations are standard normal distributed.}}}
\end{center}
\end{table}

\begin{table}
\begin{center}
{\scriptsize
\begin{tabular}{|ccc|cc|cc|cc|cc|cc|}
\hline
 & $$ & $$ &   \multicolumn{2}{c|}{\eqref{testrule} }  &   \multicolumn{2}{c|}{Baek/Pipiras}     & \multicolumn{2}{c|}{Berkes et. al} & \multicolumn{2}{c|}{Yau/Davis}  \\
\hline
$T$&$N$&$M$& $5 \%$ & $10 \%$   & $5 \%$ & $10 \%$   & $5 \%$ & $10 \%$  & $5 \%$ & $10\%$  \\ \hline
256&		64  &4  	&0.340	&0.436	&0.244	&0.343	&0.034	&0.082		&0.262	&0.335	\\
256&		32  &	8	&0.373	&0.492	&		&		&		&		&		&	\\
512&		128&	4	&0.550	&0.600	&0.434	&0.510	&0.005	&0.021	&0.228	&0.276	\\
512&		64  &	8 	&0.362	&0.476	&		&		&		&		&		&	\\
1024&	256&	4 	&0.714	&0.756	&0.527	&0.641	&0.047	&0.130	&0.197	&0.240\\
1024&	128&	8	&0.446	&0.522	&		&		&		&		&		&	\\
2048& 	512&4	&0.910	&0.926	&0.721	&0.805	&0.143	&0.244	&0.263	&0.334	\\
2048&	256&8	&0.634	&0.708	&		&		&		&		&		&	\\
4096&	1024&4	&0.974	&0.976	&0.889	&0.938	&0.311	&0.408	&0.713	&0.741	\\
4096& 	512&8	&0.923	&0.938	&		&		&		&		&		&	\\
\hline
\end{tabular}}
\caption{\textit{\label{tab2x} \small{Rejection frequencies of the test \eqref{testrule} and three competing procedures under the alternative (\ref{power1}) for different choices of T,N and M. The innovations are $ (\chi^2_5-5)/\sqrt{10}$ distributed.   }}}
\end{center}
\end{table}

\begin{table}
\begin{center}
{\scriptsize
\begin{tabular}{|ccc|cc|cc|cc|cc|cc|}
\hline
 & $$ & $$ &   \multicolumn{2}{c|}{\eqref{testrule} }  &   \multicolumn{2}{c|}{Baek/Pipiras}     & \multicolumn{2}{c|}{Berkes et. al} & \multicolumn{2}{c|}{Yau/Davis}  \\
\hline
$T$&$N$&$M$& $5 \%$ & $10 \%$   & $5 \%$ & $10 \%$   & $5 \%$ & $10 \%$  & $5 \%$ & $10\%$  \\ \hline
256&		64  &4  	&0.260	&0.330	&0.230	&0.322	&0.039	&0.088	&0.296	&0.366		\\
256&		32  &	8	&0.276		&0.394		&		&		&		&		&		&		\\
512&		128&	4	&0.528	&0.590	&0.342	&0.456	&0.010	&0.036	&0.268	&0.322	\\
512&		64  &	8 	&0.314	&0.414	&		&		&		&		&		&	\\
1024&	256&	4 	&0.774	&0.796	&0.546	&0.656	&0.024	&0.086	&0.228	&0.292\\
1024&	128&	8	&0.414	&0.492	&		&		&		&		&		&	\\
2048& 	512&4	&0.900	&0.913	&0.758	&0.820	&0.168	&0.268	&0.320	&0.404\\
2048&	256&8	&0.608	&0.665	&		&		&		&		&		&	\\
4096&	1024&4	&0.994	&0.996	&0.900	&0.940	&0.332	&0.444	&0.649	&0.697	\\
4096& 	512&8	&0.982	&0.990	&		&		&		&		&		&	\\
\hline
\end{tabular}}
\caption{\textit{\label{tab3} \small{Rejection frequencies of the test \eqref{testrule} and three competing procedures under the alternative (\ref{power2}) for different choices of T,N and M. The innovations are standard normal distributed.}}}
\label{}
\end{center}
\begin{center}
{\scriptsize
\begin{tabular}{|ccc|cc|cc|cc|cc|cc|}
\hline
 & $$ & $$ &   \multicolumn{2}{c|}{\eqref{testrule} }  &   \multicolumn{2}{c|}{Baek/Pipiras}     & \multicolumn{2}{c|}{Berkes et. al} & \multicolumn{2}{c|}{Yau/Davis}  \\
\hline
$T$&$N$&$M$& $5 \%$ & $10 \%$   & $5 \%$ & $10 \%$   & $5 \%$ & $10 \%$  & $5 \%$ & $10\%$  \\ \hline
256&		64  &4  	&0.094	&0.136	&0.087	&0.149	&0.045	&0.093	&0.178	&0.210		\\
256&		32  &	8	&0.138	&0.216	&		&		&		&		&	&		\\
512&		128&	4	&0.146	&0.196	&0.119	&0.177	&0.022	&0.055	&0.140	&0.176	\\
512&		64  &	8 	&0.138	&0.214	&		&		&		&		&		&		\\
1024&	256&	4 	&0.328	&0.406	&0.127	&0.197	&0.018	&0.079	&0.152	&0.206		\\
1024&	128&	8	&0.152	&0.218	&		&		&		&		&	&	\\
2048& 	512&4	&0.646	&0.710	&0.174	&0.266	&0.052	&0.116	&0.374	&0.470	\\
2048&	256&8	&0.312	&0.388	&		&		&		&		&	&	\\
4096&	1024&4	&0.854	&0.888	&0.232	&0.466	&0.064	&0.162	&0.736&0.792	\\
4096& 	512&8	&0.716	&0.742	&		&		&		&		&	&	\\
\hline
\end{tabular}}
\caption{\textit{\label{tab4} \small{Rejection frequencies of the test \eqref{testrule} and three competing procedures under the alternative (\ref{power3}) for different choices of T,N and M. The innovations are standard normal distributed.}}}
\label{}
\end{center}
\end{table}
\begin{table}
\begin{center}
{\scriptsize
\begin{tabular}{|ccc|cc|cc|cc|cc|cc|ccl}
\hline
 & $$ & $$ &   \multicolumn{2}{c|}{$d_1$}  &   \multicolumn{2}{c|}{$d_2$}     & \multicolumn{2}{c|}{$d_3$} &   \multicolumn{2}{c|}{$d_4$}& \multicolumn{2}{c|}{$d_5$}  \\
\hline
$T$&$N$&$M$& $5 \%$ & $10 \%$   & $5 \%$ & $10 \%$   & $5 \%$ & $10 \%$& $5 \%$ & $10 \%$    & $5 \%$ & $10 \%$ \\ \hline
256&		64  &4  	&0.118	&0.174	&0.118	&0.168	&0.146	&0.222 	&0.374 	&0.450 	&0.356	&0.432	 \\
256&		32  &	8	&0.184	&0.270	&0.167	&0.241	&0.183	&0.261	&0.359	&0.452	&0.335	&0.464	 \\
512&		128&	4	&0.198	&0.272	&0.216	&0.296	&0.350	&0.412 	&0.592 	&0.638 	&0.622	&0.662	 \\
512&		64  &	8 	&0.092	&0.146	&0.134	&0.188	&0.104	&0.188 	&0.372	&0.490 	&0.400	&0.518	 \\
1024&	256&	4 	&0.430	&0.504	&0.402	&0.504	&0.648	&0.716 	&0.792 	&0.808 	&0.776	&0.808	 \\
1024&	128&	8	&0.160	&0.226	&0.136	&0.200	&0.230	&0.290 	&0.506 	&0.608 	&0.500	&0.586	 \\
2048&	512&4 	&0.746	&0.790	&0.804	&0.844	&0.868	&0.880	&0.876 	&0.892 	&0.912	&0.932	\\
2048&	256&8 	&0.434	&0.492	&0.454	&0.510	&0.534	&0.578	&0.670 	&0.754 	&0.678	&0.744	\\
4096& 	1024&4	&0.932	&0.940	&0.930	&0.940	&0.943	&0.953	&0.982	&0.985	&0.992	&0.992	\\
4096&	512&8 	&0.914	&0.922	&0.910	&0.918	&0.910	&0.925 	&0.967	&0.978 	&0.895	&0.923	\\
\hline
\end{tabular}}
\caption{\textit{\label{tab5} \small{Rejection frequencies of the test \eqref{testrule} under the alternative (\ref{power2}) for different choices of the long-memory function $d$ defined in (\ref{long1})-(\ref{long3}). The innovations are standard normal distributed.}}}
\end{center}
\end{table}
\begin{rem}
{\rm
It is well known that fitting tvFAR or tvFARIMA to tvAR or tvARMA models yields to confounded estimates of the AR/MA coefficients and the long-memory parameter. As a consequence the approximation of the nominal level becomes less accurate if the AR polynomial $|1+\sum_{j=1}^k a_j(u) e^{-i\lambda j}|^{2}$ has roots close to the unit disc. For example, motivated by a comment of a reviewer, we have conducted a further simulation study investigating a tvAR(1) model. These results are not depicted for the sake of brevity but they clearly show that the approximation of the nominal level of the new test is not accurate if the AR coefficients vary in the interval $(0.85,1)$. In this case the level is overestimated, and the test \eqref{testrule} decides too often for a long-memory process.
}
\end{rem}
\subsection{Simulation of prediction error} \label{predictionerrorsection}
In this subsection we investigate the question what one loses by fitting a short-range dependent non-stationary model to data that is truly non-stationary and long-range dependent. For this purpose we simulate data from the tvFARIMA($0,d,1$)-process in (\ref{power2}) with long-memory functions $d_1$ and $d_4$ in (\ref{long1}) and (\ref{long3}), respectively. We separately fit a tvARMA(1,1) model and a tvFARIMA(1,d,1) model to the data and use the state space framework in \cite{palolefer2013} in order to predict future values and then compare the prediction errors of these two fitted models. To be more precise we consider the sample size $T=1024$  with block length $N=256$ (resulting in $M=4$  blocks) and use the local Whittle estimator from Section \ref{sec3} to  estimate on each block the locally varying AR and MA coefficients for the tvARMA(1,1) model and the AR, MA and long-memory parameters for the tvFARIMA($1,d,1$) model. With these time-varying coefficients we use the Kalman filter equations in \cite{palolefer2013} and calculate $5, 10$ and $25$-step predictors with each of these two models. The prediction error is calculated by sum of squared residuals
$$
\sum_{\ell=1}^k \big ( X_{t,+\ell,T} - \widehat X_{t+\ell,T} \bigr)^2 ~,~\ell = 5, 10,25.
$$
In Table \ref{tab_pred} we display the median  and  median absolute  deviation of the prediction errors obtained
in $ 1000$ simulation runs. We observe that the predictions, which take the long memory property into account are substantially more accurate.
\begin{table}
\begin{center}
{\scriptsize
\begin{tabular}{|c|cc|cc|cc|cc|cc|}
\hline
  $$ &     \multicolumn{4}{c|}{ model  (\ref{power2}) with $d_1$}   &   \multicolumn{4}{c|}{model (\ref{power2}) with $d_4$}       \\
\hline
&    \multicolumn{2}{c|}{tvARMA(1,1) }   &  \multicolumn{2}{c|}{tvFARIMA(1,d,1) }   &    \multicolumn{2}{c|}{tvARMA(1,1)
}   &  \multicolumn{2}{c|}{tvFARIMA(1,d,1)}  \\ \hline
$h$-step prediction& med & dev   & med  & dev   & med &  dev  & med & dev  \\ \hline
5  	&19.1&52.4	&	4.8	&3.6	&12.3	&42.3	&4.5		&3.3	\\
10	&25.2&56.5	&	10.7	&5.0		&18.1	&44.0	&10.1	&4.7		\\
25	&43.2&54.6	&	25.4	&7.8		&36.5	&40.1	&26.3	&10.8	\\
\hline
\end{tabular}}
\caption{\textit{\label{tab_pred} \small{ Prediction error by a fit of  tvARMA(1,1) and  tvFARIMA(1,d,1)  models (median and median
absolute deviation obtained by $1000$ simulation runs).}}}
\end{center}
\end{table}

\subsection{Data examples}

\begin{figure}
\begin{center}
\includegraphics[width=0.5\textwidth]{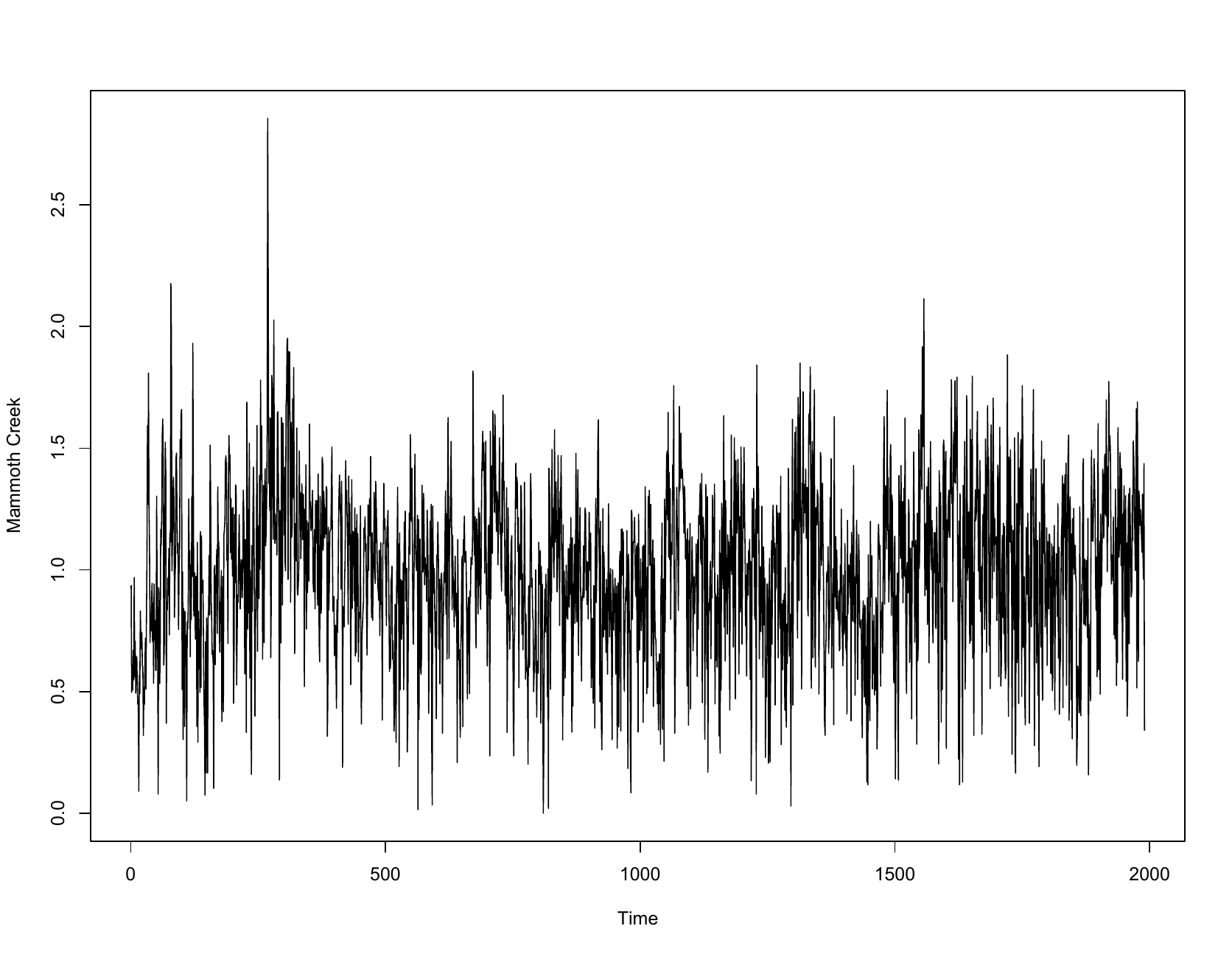}~~
\includegraphics[width=0.5\textwidth]{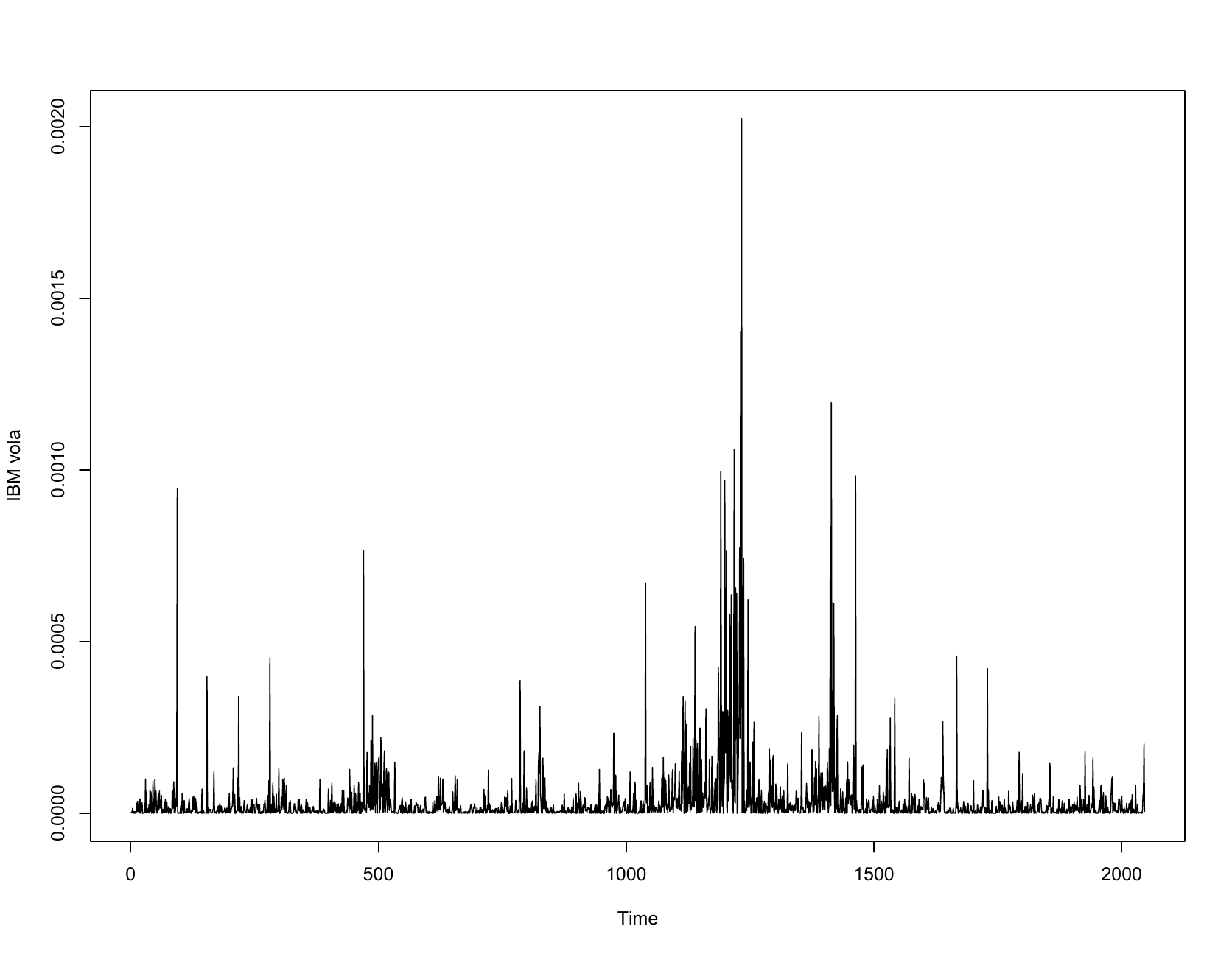}~~
 \caption{ \label{twodatasets} {\it Left panel: plot of the 1990 annual pinus longaeva tree ring width measurements at Mammoth Creek, Utah from 0 A.D. to 1989 A.D.; Right panel: plot of the squared log-returns of the IBM stock between July $15$th $2005$ and  August $30$th $2013$.} }
\end{center}
\end{figure}

{\bf Testing:} As an illustration we apply the new test to two different datasets, where in both examples the mean function has been estimated as described in Section \ref{sec3}. As pointed out in the previous section the quality of predictions can be improved, if long range dependence
is  present in non stationary data and  considered in the predictions. For this reason the test proposed in this paper can be useful to obtain more accurate forcasts. \\
The first data set contains annual pinus longaeva tree ring width measurements at Mammoth Creek, Utah, from 0 A.D. to 1989 A.D.  while the second data set contains 2048 squared log-returns of the IBM stock between July $15$th $2005$ and  August $30$th $2013$ which was already discussed in the introduction. Both time series are depicted in Figure \ref{twodatasets}, and in the case of the tree ring data our test statistic $\sqrt{T} \hat F_T/\sqrt{\hat W_{T}}$ equals 17.8 for $M=4$ and yields a p-value  $\approx 0$. This implies that the null hypothesis of a non-stationary short-memory model has to be rejected for this dataset, which coincides with the results of the tests in \cite{baekpipi}  and \cite{Yau2012}. Their test statistics have the values 3.49 and 9.37 and p-values of $0.00024$ and $0$ corresponding to the local Whittle and likelihood ratio approach, respectively. The CUSUM procedure of \cite{horvarth2006} yields a value of $0.906$ for the test statistic and does not reject the null hypothesis at even $10\%$ nominal level. This result
 is possibly due to the low power of this test as remarked in Section \ref{Simulation results}.

In the situation of the squared log-returns of the IBM stock, the assumption of Gaussianity is too restrictive and we therefore apply the more general test described in Remark \ref{prozess nicht-gauss}. The values of the test statistic $\sqrt{T} \hat F_T/\sqrt{\hat W_{T,general}}$ are $5.67$ and $9.48$ for  $M=4$ and $M=8$, respectively, yielding that the p-value is smaller than $2.87 \cdot 10^{-7}$ for both choices of the segmentation. This means that the assumption of no long-range dependence is clearly rejected. If we apply the likelihood ratio test of \cite{Yau2012} to this dataset, we obtain a value for the statistic of $15.77$ which is then compared with the quantiles of the standard normal distribution. This yields also to a rejection of the null hypothesis. On the other hand, the CUSUM procedure of \cite{horvarth2006} only rejects the null hypothesis of no long-range dependence at a $10\%$ but not at a $5\%$ level. This observation is, however, not surprising given the low power of this test in the finite sample situations presented in the previous section. The test of \cite{baekpipi} rejects the null hypothesis with a p-value $8.65 \cdot 10^{-12}$, yielding the same result as our approach and the one of \cite{Yau2012}.
\medskip

\textbf{Prediction:}
The result of the test \eqref{testrule} has important consequences for the subsequent data analysis as it advices the statistician to use short memory or long memory (non-stationary) models.
In the final part of this section  we demonstrate how the information of the test can be employed to obtain superior forecasting results in
the two datasets analyzed in the previous paragraph. For this purpose, we divide both datasets into two parts. One  part contains the first $3/4 \times T$ observations of the corresponding dataset while the second part contains the remaining $T/4$ data points. The  new testing procedure  \eqref{testrule} proposed in this paper  is applied to the first part
of the data, and -
depending on the result of the test - forecasts are performed using either a tvFARIMA(1,d,1) or a ARMA(1,1)
with the window of size $N=256$ in the localized Whittle estimator [see also Section \ref{predictionerrorsection}]. In order to compare the forecasting performance of the short- with the long-memory model,
we define the prediction error on the second part of each dataset by
$$
PE(h)=\sum_{t=3/4T+1}^T \sum_{\ell=1}^h \big ( X_{t+\ell,T} - \widehat X_{t+\ell,T} \bigr)^2 ~,~h = 5, 10,25,
$$
and denote with $PE_{short}(h)$ and $PE_{long}(h)$ the prediction error for the short- and long-memory approach respectively. The expression
$$R(h)={PE_{long}(h)\over PE_{short}(h)}
$$
 then serves as a  measure for the comparison. It   is smaller than one if the long-range dependence approach yields superior predictions, while it is larger than one in the other case. As in the previous paragraph (where we applied the test to the total sample), an application of
 the  test  \eqref{testrule} to the first $3/4 \times T$ observations of the Mammoth Creek and the IBM dataset yields $p$-values much smaller than one percent
 in both cases. Consequently one  would perform data analysis on the basis of a non-stationary long range dependent model.
The advantages  of this approach are clearly visible in Table \ref{R_pred} where we depict the ratio of the prediction error from a short and long range dependent model.
We observe that the long-range dependence approach, in fact, yields substantially smaller prediction errors. In all cases the  prediction error from
 the long-range dependent
model  is less than one third of the corresponding error from the short-memory model  (for both datasets and all considered values of $h$).
This demonstrates that the difference in forecasting performance is huge and highlights the importance of powerful  tests to discriminate between long- and short-range dependence.

\begin{table}
\begin{center}
{\scriptsize
\begin{tabular}{|c|ccc|}
\hline
  $$ &     \multicolumn{3}{c|}{ $R(h)$} \\
\hline
dataset &    $h=5$   &  $h=10$   &    $h=25$ \\ \hline
Mammoth Creek Data  	&$0.21$& $0.24$	&	$0.26$	\\
IBM Data  	&$0.09$& $0.15$	&	$0.33$	\\
\hline
\end{tabular}}
\caption{\textit{\label{R_pred} \small{Comparison  of prediction errors
for the Mammoth Creek and IBM dataset with different values of $h$. A value
smaller  than one indicates a better performance of the long range dependent model.
}}}
\end{center}
\end{table}

\section{Conclusions} \label{seccon}
\def\theequation{6.\arabic{equation}}
\setcounter{equation}{0}

In this paper we have  developed  a test for  weak against  strong (long-range) dependence for non-stationarity
time series. Our approach is based on an average of unconstrained Whittle-likelihood estimates
   of the (nonnegative)  local long-range dependence
parameter from a sequence of approximating time varying FARIMA models [see equation  \eqref{fhat} for its definition]. It is demonstrated that
a standardized version of this average is asymptotically normal distributed, which
provides a very simple asymptotic level $\alpha$ and consistent
test for discriminating between short and long range dependence  of   a non-stationary time series. \\
As an alternative to the statistic $\hat F_{T} $ in \eqref{fhat} one could form an average of constrained Whittle-likelihood estimates, say $\hat d_{N,c} (u_i)$. Constrained
parameter estimation has found  considerable attention
in the literature [see for example \cite{chernoff1954} or \cite{andrews1999} among  many others], but - to our best knowledge - it has not been considered so far
in locally stationary processes. The ``classical''   results indicate  that for a fixed
value $u$ the asymptotic distribution  of
 $\hat d_{N,c} (u)$ is given by a function of a multivariate
normal distribution (in the simplest case a half normal type distribution). However, we expect that - due to averaging - the (standardized) statistic
 $\hat F_{T,c} =   {1 \over M} \sum_{i=1}^M  \hat d_{N,c} (u_i)$  is still asymptotically normal distributed.  An interesting direction for future
 research  is the development of an asymptotic theory for  constrained estimators in  locally stationary (long memory) processes
and to use it for  a  rigorous investigation of  the asymptotic properties of the statistic   $\hat F_{T,c} $.  Moreover,  the results of Section
\ref{sec2A} for the nonparametric regression model with independent errors  indicate  some  advantages of unconstrained
 over unconstrained averages,  and it will be of interest  to investigate if the superiority of $\hat F_{T} $ over  $\hat F_{T,c} $
  can also be observed for the testing problem considered in this paper. \\
  It is also notable that this paper has its focus on discriminating between short and standard long-range dependence,
  which corresponds to a pole of the local spectral density at frequency $0$. However, it was pointed out by several authors
  [see for example \cite{arteche2000,hidalgo2004,reisen2006} among others]  that - due to strong  cyclic components -
  strong dependency can also  occur as a pole in the spectral density at any other frequency (reflecting strong seasonal long range dependence).
  In this case  the analogue of the model \eqref{spec} is given by
\be \label{spec1}
f(u,\lambda)=  |1-e^{i(\lambda-\lambda_0) }|^{-d_0(u)}|1-e^{i(\lambda+ \lambda_0) }|^{-d_0(u)} g(u,\lambda),
\ee
  where $\lambda_0$  denotes the unknown pole [see \cite{hidalgo2004}], and  a  further interesting direction of future research
  is the construction of tests for the hypothesis \eqref{H_0} in the more general model \eqref{spec1}. \\
{We finally note that several authors have analyzed financial data under
  linearity  assumptions as made in equation  \eqref{proc} [see \cite{mikosch2004}, \cite{motivation1} or \cite{motivation2}
  among others]. On the other hand it is also argued in the literature  that this assumption
might not be reasonable in some cases. Long range dependent processes have
mainly been investigated in models with linear representations. A  nonlinear (nonparametric) extension does not seem to be obvious
as indicated by the results of \cite{grusur2014}, who proposed a linear representation with random coefficients.
Therefore, an interesting problem for future research is to investigate if the methodology
suggested in this paper is also valid for processes with nonlinear representations. }

\noindent
\\
{\bf Acknowledgements.}
This work has been supported in part by the
Collaborative Research Center ``Statistical modeling of nonlinear
dynamic processes'' (SFB 823, Teilprojekt A1, C1) of the German Research Foundation (DFG). 
We would also like to thank Tobias Kley and Kevin Kokot  for computational assistance, Lutz D\"umbgen for  some helpful discussions 
about Section 
\ref{sec2A}  and
Wilfredo Palma for some constructive discussion on the subject and for making the R-Code in \cite{palolefer2013} available to us.

\bibliographystyle{apalike}

 \bibliography{presendet}

\section{Appendix: Proofs}
\label{appendix}
\subsection{Preliminary results}
\def\theequation{7.\arabic{equation}}
\setcounter{equation}{0}
We begin stating two results, which will be the main tools in the asymptotic analysis of the proposed estimators and the test statistic. For this purpose, we let
$\phi_T: [0,1] \times [-\pi,\pi] \rightarrow \mathbb{R}$ denote a function which (might) depend on the the sample size $T$ and
 define
\begin{eqnarray*}
G_T(\phi_T)&=&  \frac{1}{M} \sum_{j=1}^{M} \int_{-\pi}^{\pi} f(u_j,\lambda) \phi_T(u_j,\lambda) \,d\lambda, \\
\hat G_T(\phi_T)&=&\frac{1}{M}\sum_{j=1}^{M} \int_{-\pi}^{\pi} I_N^{ \mu}(u_j,\lambda)\phi_T(u_j,\lambda) \,d\lambda,
\end{eqnarray*}
where $I_N^{\mu}$ is the analogue of the local periodogram (\ref{per}) where the estimator $\hat \mu$ has been replaced by the ``true'' mean function $\mu$.
\begin{satz} \label{zgws}~~\\
\vspace{-.8cm}

 \begin{itemize}
\item[a)] Let Assumption \ref{ass1} be fulfilled and  assume that $\phi_T(u,\lambda): [0,1] \times [-\pi,\pi] \rightarrow \mathbb{R}$ is symmetric in $\lambda$, twice continuously  differentiable with uniformly bounded partial derivatives such that for all $u \in [0,1]$, $\lambda \in [-\pi,\pi]$, $T \in \mathbb{N}$
\begin{eqnarray}
\phi_T(u,\lambda) & \leq & Cg(k) |\lambda|^{2d_0(u)-\epsilon} ,
\label{Ordnung theta} \\
\frac{\partial }{\partial \lambda}\phi_T(u,\lambda)& \leq & C g(k) |\lambda|^{2d_0(u)-1-\epsilon}  ,
\label{Ordnung theta 1.Ableitung} \\
\frac{\partial^2 }{\partial \lambda^2}\phi_T(u,\lambda) & \leq & Cg(k) |\lambda|^{2d_0(u)-2-\epsilon} ,
\label{Ordnung theta 2.Ableitung}
\end{eqnarray}
where $C>0, 0<\epsilon<1/2-D$ are constants and $g: \mathbb{N} \rightarrow (0,\infty)$ is a given function. Then we have
\be
\E [\hat G_T(\phi_T)] &=& G_T(\phi_T)+ O\Big(\frac{g(k)}{N^{1-\epsilon}}\Big) + O \Big(\frac{g(k)N^2}{T^2} \Big) \label{ew}, \\
\Var[\hat G_T(\phi_T)]&=&  V_T + O\Big(\frac{1}{T}\frac{ g^2(k)}{N^{1-2D-2\epsilon}}\Big) +O \Big(\frac{ g^2(k)N^2}{T^3} \Big) \label{var}
\ee
where
\begin{eqnarray*}
V_T=\frac{1}{T}\frac{4\pi}{M} \sum_{j=1}^{M} \int_{-\pi}^{\pi}f^2(u_j,\lambda) \phi_T^2(u_j, \lambda)\,d\lambda .
\end{eqnarray*}
\item[b)] Suppose the assumptions of part a) hold with $D=0$, $\epsilon< 1/6$ and additionally $\liminf_{T \rightarrow \infty }T \cdot  V_T \geq c,$
\begin{eqnarray*}
&&N \rightarrow \infty, \hspace{.3cm}
g(k)\sqrt{T}/N^{1-\epsilon}  \rightarrow 0 , \hspace{.3cm}
 g(k)\log(T)/T^{1/6-\epsilon}  \rightarrow 0,\hspace{.3cm}
 \mbox{ and } \hspace{.3cm}
g(k)N^2/T^{\frac{3}{2}}\rightarrow 0.
 \end{eqnarray*}
Then we have $\sqrt{T}\big(\hat G_T(\phi_T)-G_T(\phi_T)\big) / \sqrt{ V_T} \stackrel{\mathcal{D}}{\rightarrow} \mathcal{N}(0,1).$
\end{itemize}
\end{satz}

\textbf{Proof:} In order to prove part a)  Theorem \ref{zgws}
we define $\tilde{t}_j:= t_j-N/2+1$, $\tilde{\psi}_l(u_{j,p}):=\psi_l(\frac{\tilde{t}_j+p}{T})$, $Z_{a,b}:=Z_{a-N/2+1+b}$ and obtain
\begin{eqnarray*}
\E [\hat G_T(\phi_T)]&=&\frac{1}{2 \pi N}\frac{1}{M} \sum_{j=1}^M  \sum_{p,q=0}^{N-1} \sum_{l,m=0}^{\infty}
\psi_{\tilde{t}_j+p,T,l}\psi_{\tilde{t}_j+q,T,m}\int_{-\pi}^{\pi}\phi_T(u_j,\lambda)e^{-i(p-q)\lambda}\,d\lambda \E(Z_{t_j,p-l}Z_{t_j,q-m}) \\
&=& E_{N,T}+A_{N,T}+B_{N,T}
\end{eqnarray*}
where
\begin{eqnarray*}
E_{N,T}&:=&\frac{1}{2 \pi N}\frac{1}{M} \sum_{j=1}^M  \sum_{p,q=0}^{N-1} \sum_{l,m=0}^{\infty}
\psi_l(u_j)\psi_m(u_j)\int_{-\pi}^{\pi}\phi_T(u_j,\lambda)e^{-i(p-q)\lambda}\,d\lambda \E(Z_{t_j,p-l}Z_{t_j,q-m}) ,\\
A_{N,T}&:=&\frac{1}{2 \pi N}\frac{1}{M} \sum_{j=1}^M  \sum_{p,q=0}^{N-1} \sum_{l,m=0}^{\infty}
\int_{-\pi}^{\pi}\phi_T(u_j,\lambda)e^{-i(p-q)\lambda}\,d\lambda \E(Z_{t_j,p-l}Z_{t_j,q-m}) \\
&&\big\{\psi_l(u_j)\big(\tilde{\psi}_m(u_{j,q})-\psi_m(u_{j}) \big) + \big(\tilde{\psi}_l(u_{j,p})-\psi_l(u_j)\big)\tilde{\psi}_m(u_{j,q})\big\} ,\\
B_{N,T}&:=&\frac{1}{2 \pi N}\frac{1}{M} \sum_{j=1}^M  \sum_{p,q=0}^{N-1} \sum_{l,m=0}^{\infty}
\int_{-\pi}^{\pi}\phi_T(u_j,\lambda)e^{-i(p-q)\lambda}\,d\lambda \E(Z_{t_j,p-l}Z_{t_j,q-m}) \\
&&\big\{\big( \psi_{\tilde{t}_j+p,T,l}-\tilde{\psi}_l(u_{j,p})\big)\psi_{\tilde{t}_j+q,T,m} +\tilde{\psi}_l(u_{j,p})\big(\psi_{\tilde{t}_j+q,T,m} -\tilde{\psi}_m(u_{j,q})\big) \big\}.
\end{eqnarray*}
Note that $B_{N,T}$ and $A_{N,T}$ compromise the error arising in the approximation of $\psi_{\tilde{t}_j+p,T,l}$ by $\psi_l(\frac{\tilde{t}_j+p}{T})$ and $\tilde{\psi}_m(u_{j,q})$ by $\psi_m(u_j)$, respectively. In order to establish the claim (\ref{ew}), we prove the following statements:
\begin{eqnarray}
E_{N,T}&=&\frac{1}{M} \sum_{j=1}^M\int_{-\pi}^{\pi}f(u_j,\lambda) \phi_T(u_j,\lambda)  \,d\lambda+O\Big(\frac{g(k)}{N^{1-\epsilon}}\Big)
\label{E_{N,T}}\\
A_{N,T}&=&
O \Big( \frac{g(k) \log(N)}{N^{1-\epsilon} M } \Big)+O \Big (\frac{g(k)N^2}{T^2} \Big)
\label{A_{N,T}}\\
B_{N,T}&=& O\Big( \frac{g(k)\log(T)}{T}\Big).
\label{B_{N,T}}
\end{eqnarray}
\textit{Proof of (\ref{E_{N,T}}):} Due to the independence of the random variables $Z_t$, we only need to consider terms fulfilling $p=q+l-m$ (this means $0 \leq p=q+l-m \leq N-1$ because of $p \in \{0,1,2 \ldots , N-1 \}$) which in turn implies $|l-m| \leq N-1$. Therefore
\begin{eqnarray*}
E_{N,T}&=&\frac{1}{2 \pi N}\frac{1}{M} \sum_{j=1}^M \sum_{\substack{l,m=0 \\ |l-m| \leq N-1}}^{\infty}  \sum_{\substack{q=0 \\ 0 \leq q+l-m \leq N-1 }}^{N-1}  \psi_l(u_j) \psi_m(u_j)\int_{-\pi}^{\pi}\phi_T(u_j,\lambda)e^{-i(l-m)\lambda}\,d\lambda \\
&=&\frac{1}{2 \pi N}\frac{1}{M} \sum_{j=1}^M
\sum_{\substack{l,m=0\\ |l-m| \leq N-1}}^{\infty}
\psi_l(u_j)\psi_m(u_j)\int_{-\pi}^{\pi}\phi_T(u_j,\lambda)e^{-i(l-m)\lambda}\,d\lambda (N-|l-m|)\\
&=&\frac{1}{M} \sum_{j=1}^M\int_{-\pi}^{\pi} \phi_T(u_j,\lambda) f(u_j,\lambda) \,d\lambda  +E^1_{N,T}+E^2_{N,T},
\end{eqnarray*}
where
\begin{eqnarray*}
E^1_{N,T}&=& -\frac{1}{2 \pi }\frac{1}{M} \sum_{j=1}^M   \sum_{\substack{l,m=0 \\ N \leq  |l-m|}}^{\infty} \psi_l(u_j) \psi_m(u_j)\int_{-\pi}^{\pi}\phi_T(u_j,\lambda)e^{-i(l-m)\lambda}\,d\lambda , \\
E^2_{N,T}
&=&-\frac{1}{2 \pi N}\frac{1}{M} \sum_{j=1}^M  \sum_{\substack{l,m=0 \\ |l-m| \leq N-1}}^{\infty} \psi_l(u_j) \psi_m(u_j)\int_{-\pi}^{\pi}\phi_T(u_j,\lambda)e^{-i(l-m)\lambda}\,d\lambda  |l-m| .
\end{eqnarray*}
Using (\ref{apprpsi}), (\ref{Ordnung theta}) and Lemma \ref{fourier} in the online supplement, we obtain
\begin{eqnarray*}
|E^1_{N,T} |& \leq& C\frac{g(k)}{M} \sum_{j=1}^M   \sum_{\substack{l,m=1 \\ N \leq  |l-m|}}^{\infty} \frac{1}{l^{1-d_0(u_j)}} \frac{1}{m^{1-d_0(u_j)}}\frac{1}{|l-m|^{1+2d_0(u_j)-\epsilon}}(1+o(1)),
\end{eqnarray*}
where we used the fact that terms corresponding to $l=0$ or $m=0$ are of  smaller or  the same order (we will use this property frequently from now on without further mentioning it). We set  $h:=l-m$ and obtain from Lemma \ref{Lemma 7.1}a) in the online supplement that
\begin{eqnarray*}
& & \frac{g(k)}{M} \sum_{j=1}^M \sum_{\substack{h \in \mathbb{Z} \\ N \leq |h|}}^{}  \sum_{\substack{m=1 \\ h+m \geq 1}}^{\infty} \frac{1}{(h+m)^{1-d_0(u_j)}} \frac{1}{m^{1-d_0(u_j)}}\frac{1}{|h|^{1+2d_0(u_j)-\epsilon}} \leq C g(k)  \sum_{\substack{h \in \mathbb{Z} \\ N \leq |h|}}^{} \frac{1}{|h|^{2-\epsilon}}  =O\Big( \frac{g(k)}{N^{1-\epsilon}}\Big).
\end{eqnarray*}
By proceeding analogously we obtain that $E_{N,T}^2=O(g(k) N^{-1+\epsilon})$ which proves the assertion in (\ref{E_{N,T}}). \\

\textit{Proof of (\ref{A_{N,T}}):} Without loss of generality we only consider the first summand
$$
A_{N,T}(1) =\frac{1}{2 \pi N}\frac{1}{M} \sum_{j=1}^M  \sum_{p,q=0}^{N-1} \sum_{l,m=0}^{\infty}\psi_l(u_j)\big(\tilde{\psi}_m(u_{j,q})-\psi_m(u_{j}) \big)
\int_{-\pi}^{\pi}\phi_T(u_j,\lambda)e^{-i(p-q)\lambda}\,d\lambda \E(Z_{t_j,p-l}Z_{t_j,q-m}) \notag
$$
in $A_{N,T}$ (the second term is treated exactly in the same way). A Taylor expansion
and similar arguments as in the proof of (\ref{E_{N,T}}) yield
\begin{eqnarray*}
A_{N,T}(1)  &=&A^1_{N,T}+ A^2_{N,T}
\end{eqnarray*}
where
\begin{eqnarray*}
A^1_{N,T}&=&\frac{1}{2 \pi N}\frac{1}{M} \sum_{j=1}^M   \sum_{\substack{l,m=0\\ |l-m| \leq N-1 }}^{\infty}  \sum_{\substack{q=0 \\ 0 \leq q+l-m \leq N-1 }}^{N-1}
  \psi_l(u_j)\psi_{m}^{'}(u_{j}) \Big(\frac{-N/2+1+q}{T} \Big) \int_{-\pi}^{\pi}\phi_T(u_j,\lambda)e^{-i(l-m)\lambda}\,d\lambda, \\
A^2_{N,T} &= &\frac{1}{2 \pi N}\frac{1}{M} \sum_{j=1}^M   \sum_{\substack{l,m=0\\ |l-m| \leq N-1 }}^{\infty}  \sum_{\substack{q=0 \\ 0 \leq q+l-m \leq N-1 }}^{N-1}
  \psi_l(u_j)\psi^{''}_{m}(\eta_{m,j,q})\Big(\frac{-N/2+1+q}{T} \Big)^2\int_{-\pi}^{\pi}\phi_T(u_j,\lambda)e^{-i(l-m)\lambda}\,d\lambda
\end{eqnarray*}
and  $\eta_{m,j,q} \in (u_j-N/(2T),u_j+N/(2T))$.
Using (\ref{apprpsi}), (\ref{2.1c}), (\ref{Ordnung theta}), Lemma \ref{fourier} it follows
\begin{eqnarray*}
|A^1_{N,T}|&\leq &C\frac{g(k)}{N}\frac{1}{M} \sum_{j=1}^M  \sum_{\substack{l,m=1 \\ 1 \leq |l-m| \leq N-1}}^{\infty}
\frac{1}{l^{1-d_0(u_j)}}\frac{\log(m)}{m^{1-d_0(u_j)}}\frac{1}{|l-m|^{1+2d_0(u_j)-\epsilon}}\Big| \sum_{\substack{q=0 \\ 0 \leq q+l-m \leq N-1}}^{N-1}\Big(\frac{-N/2+1+q}{T} \Big)\Big| \\
&\leq &C\frac{g(k)}{ T}\frac{1}{M} \sum_{j=1}^M  \sum_{\substack{l,m=1 \\ 1 \leq |l-m| \leq N-1}}^{\infty}
\frac{1}{l^{1-d_0(u_j)}}  \frac{\log(m)}{m^{1-d_0(u_j)}}\frac{1}{|l-m|^{2d_0(u_j)-\epsilon}}\\
&= &C\frac{g(k)}{ T}\frac{1}{M} \sum_{j=1}^M  \sum_{\substack{s \in \mathbb{Z} \\ 1 \leq |s| \leq N-1}}^{}  \sum_{\substack{l=1 \\ 1 \leq l-s}}^{\infty}
\frac{1}{l^{1-d_0(u_j)}} \frac{\log(l-s)}{(l-s)^{1-d_0(u_j)}} \frac{1}{|s|^{2d_0(u_j)-\epsilon}}\\
&\leq& C\frac{g(k) \log(N)}{ T}\frac{1}{M} \sum_{j=1}^M \sum_{\substack{s \in \mathbb{Z} \\ 1 \leq |s| \leq N-1}}^{} \frac{1}{|s|^{1-\epsilon}}  = O\Big(\frac{ g(k)\log(N)}{N^{1-\epsilon} M }\Big)
\end{eqnarray*}
where we used Lemma \ref{Lemma 7.1}(c) in the online supplement for the last step.
Finally, (\ref{apprpsi}), (\ref{2.1c}), (\ref{Ordnung theta}), Lemma \ref{fourier} in the online supplement and the same arguments as above, show that the term $A^2_{N,T}$ is of order $O(g(k)N^2T^{-2})$. \\

\textit{Proof of (\ref{B_{N,T}}):} By employing (\ref{apprbed}) and the same arguments as above it can be shown that $B_{N,T}$ is of order $O(\frac{g(k)\log(T)}{T})$. \\

In the next step we prove the asymptotic representation
for the variance in  (\ref{var}). We obtain
\begin{eqnarray*}
\Var(\hat G_T(\phi_T))&=& \frac{1}{(2 \pi N)^2}\frac{1}{M^2} \sum_{j_1,j_2=1}^M \sum_{p,q,r,s=0}^{N-1}\sum_{l,m,n,o=0}^{\infty} \psi_l(u_{j_1})\psi_m(u_{j_1})\psi_n(u_{j_2})\psi_o(u_{j_2}) \\
&&\cum(Z_{t_{j_1},p-l}Z_{t_{j_1},q-m} ,Z_{t_{j_2},r-n}Z_{t_{j_2},s-o})\int_{-\pi}^{\pi} \phi_T(u_{j_1},\lambda_1)e^{-i(p-q)\lambda_1} \,d\lambda_1
\int_{-\pi}^{\pi} \phi_T(u_{j_2},\lambda_2)e^{-i(r-s)\lambda_2} \,d\lambda_2\\
&&+ O\Big(\frac{ g^2(k) \log(N)}{TN^{1-\epsilon}M}\Big)+O\Big(\frac{ g^2(k)N^2}{T^3} \Big),
\end{eqnarray*}
where we used assumption (\ref{apprbed}) and similar arguments as given in the proof of \eqref{ew}. Because of the Gaussianity of the innovations we obtain
\bea
\cum(Z_{t_{j_1},p-l}Z_{t_{j_1},q-m} ,Z_{t_{j_2},r-n}Z_{t_{j_2},s-o}) &=& \E(Z_{t_{j_1},p-l}Z_{t_{j_2},r-n}) \E(Z_{t_{j_1},q-m}Z_{t_{j_2},s-o} ) \\
&&+\E(Z_{t_{j_1},p-l}Z_{t_{j_2},s-o}) \E(Z_{t_{j_1},q-m}Z_{t_{j_2},r-n}).
\eea
This implies that the calculation of the (dominating part of the) variance splits into two sums, say $V_{N,T}^1$ and $V_{N,T}^2$. In the following discussion we will show that both terms converge to the same limit, that is
\begin{eqnarray*}
V^i_{N,T}&=& \frac{1}{T}\frac{2\pi}{M}\sum_{j=1}^M  \int_{-\pi}^{\pi}f^2(u_{j},\lambda) \phi_T^2(u_{j},\lambda) \,d\lambda + O\Big(\frac{1}{T}\frac{ g^2(k)}{N^{1-2D-2\epsilon}}\Big); \quad i=1,2
\end{eqnarray*}
For the sake of brevity we restrict ourselves to the case $i=1$. Because of the independence of the innovations $Z_t$, we obtain that the conditions $p=r+l-n+(j_2-j_1)N$ and $s=q+o-m+(j_1-j_2)N$ must hold, which, because of $p,s \in \{0,...,N-1\}$, directly implies  $|l-n+(j_2-j_1)N| \leq N-1$ and  $|o-m+(j_1-j_2)N| \leq N-1$. Thus, the term $V^1_{N,T}$ can be written as
\begin{eqnarray*}
&&  \frac{1}{(2 \pi N)^2}\frac{1}{M^2} \sum_{j_1=1}^M \sum_{\substack{q,r=0  }}^{N-1}\sum_{\substack{l,m,n,o=0  }}^{\infty} \sum_{\substack{j_2=1 \\ 0 \leq r+l-n+(j_2-j_1)N \leq N-1 \\ 0 \leq q+o-m+(j_1-j_2)N \leq N-1 \\ |l-n+(j_2-j_1)N| \leq N-1 \\ |o-m+(j_1-j_2)N| \leq N-1 }}^M \psi_l(u_{j_1}) \psi_m(u_{j_1})\psi_n(u_{j_2})  \psi_o(u_{j_2})\\
&&\times\int_{-\pi}^{\pi} \phi_T(u_{j_1},\lambda_1)e^{-i(r-q+l-n+(j_2-j_1)N)\lambda_1} \,d\lambda_1
\int_{-\pi}^{\pi} \phi_T(u_{j_2},\lambda_2)e^{-i(r-q+m-o+(j_2-j_1)N)\lambda_2} \,d\lambda_2  .
\end{eqnarray*}
Since $q \in \{ 0,1,2 \ldots, N-1\}$, we get from the condition $0 \leq q+o-m+(j_1-j_2)N \leq N-1$ that, if $q,o,m,j_1$ are fixed, there are at most two possible values for $j_2$ such that the corresponding term does not vanish. It follows from Lemma \ref{Fehler Varianz}  (i)--(iii) in the online supplement that there appears an error of order $O(\frac{1}{T}\frac{ g^2(k)}{N^{1-2D-2\epsilon}})$ if we drop the condition $0 \leq r+l-n+(j_2-j_1)N \leq N-1$ and assume that the variable $r$ runs from $-(N-1)$ to $-1$. Therefore, up to an error of order $O(\frac{1}{T}\frac{ g^2(k)}{N^{1-2D-2\epsilon}})$, the term $V^1_{N,T}$ is equal to
$$
D_{1,T}+D_{2,T},
$$
where
\begin{eqnarray*}
D_{1,T}&=& \frac{1}{(2 \pi N)^2}\frac{1}{M^2} \sum_{j_1=1}^M \sum_{\substack{q=0  }}^{N-1} \sum_{\substack{r=-(N-1)  }}^{N-1} \sum_{\substack{l,m,n,o=0 }}^{\infty} \sum_{\substack{j_2=1  \\ 0 \leq q+o-m+(j_1-j_2)N \leq N-1 \\ |l-n+(j_2-j_1)N| \leq N-1 \\ |o-m+(j_1-j_2)N| \leq N-1  }}^M \psi_l(u_{j_1}) \psi_m(u_{j_1})\psi_n(u_{j_2})  \psi_o(u_{j_2})\\
&&\times\int_{-\pi}^{\pi} \phi_T(u_{j_1},\lambda_1)\phi_T(u_{j_2},\lambda_1)e^{-i(r-q+l-n+(j_2-j_1)N)\lambda_1} \,d\lambda_1
\int_{-\pi}^{\pi} e^{-i(r-q+m-o+(j_2-j_1)N)\lambda_2} \,d\lambda_2  \\
D_{2,T}&=& \frac{1}{(2 \pi N)^2}\frac{1}{M^2}\sum_{j_1=1}^M \sum_{\substack{q=0  }}^{N-1} \sum_{\substack{r=-(N-1)  }}^{N-1} \sum_{\substack{l,m,n,o=0 }}^{\infty} \sum_{\substack{j_2=1  \\ 0 \leq q+o-m+(j_1-j_2)N \leq N-1 \\ |l-n+(j_2-j_1)N| \leq N-1 \\ |o-m+(j_1-j_2)N| \leq N-1  }}^M \psi_l(u_{j_1}) \psi_m(u_{j_1})\psi_n(u_{j_2})  \psi_o(u_{j_2})\\
&&\times\int_{-\pi}^{\pi} \phi_T(u_{j_1},\lambda_1)e^{-i(r-q+l-n+(j_2-j_1)N)\lambda_1}
\int_{-\pi}^{\pi}\Big[ \phi_T(u_{j_2},\lambda_2)-\phi_T(u_{j_2},\lambda_1) \Big]e^{-i(r-q+m-o+(j_2-j_1)N)\lambda_2} \,d\lambda_2  \,d\lambda_1 .
\end{eqnarray*}
We show
\begin{eqnarray}
D_{1,T} &=& \frac{2\pi}{N}\frac{1}{M^2} \sum_{j_1=1}^M  \int_{-\pi}^{\pi}f^2(u_{j_1},\lambda_1) \phi_T^2(u_{j_1},\lambda_1) \,d\lambda_1 + O\Big(\frac{1}{T}\frac{ g^2(k)}{N^{1-2D-2\epsilon}}\Big)
\label{d1T}\\
D_{2,T} &=&O\Big(\frac{1}{T}\frac{ g^2(k)}{N^{1-2D-2\epsilon}}\Big),
\notag
\end{eqnarray}
which then concludes the proof of (\ref{var}). For this purpose we begin with an investigation of the term $D_{1,T}$ for which the terms in the sum vanish if $r-q+m-o+(j_2-j_1)N\not=0$. Moreover, the following facts are correct:
\begin{itemize}
\item[I.] The variable $r$ runs from $0$ to $N-1$ since $r-q+m-o+(j_2-j_1)N=0$ and
$0 \leq q+o-m+(j_1-j_2)N \leq N-1$.
\item[II.] We can drop the condition $|l-n+(j_2-j_1)N| \leq N-1 $ by making an error of order $O( g^2(k)T^{-1}N^{-1+2D+2\epsilon})$ [this follows from Lemma  \ref{Fehler Varianz}(iv) in the online supplement].
\item[III.] There appears an error of order $O( g^2(k)T^{-1}N^{-1+2D+2\epsilon})$  if we omit the sum with $j_1 \not = j_2$ [we prove this in Lemma \ref{Fehler Varianz}(v) in the online supplement].
\item[IV.] We can afterwards omit the condition $0 \leq q+o-m \leq N-1$ since it is $0 \leq r \leq N-1$ and $r-q+m-o=0$ [note that, because of III., we assume $j_1=j_2$ from now on].
\item[V.] We can then drop the condition $ |o-m| \leq N-1$ since $r-q+m-o=0$ and $|r-q| \leq N-1$.
\end{itemize}
Thus, using the representation of $f(u_{j_1}, \lambda)$ in (\ref{tvspectraldensity}), the term $D_{1,T}$ can be written as (up to an error of order $O( g^2(k)T^{-1}N^{-1+2D+2\epsilon})$)
\begin{eqnarray*}
&& \frac{1}{N^2}\frac{1}{M^2} \sum_{j_1=1}^M\sum_{\substack{q,r=0  }}^{N-1} \int_{-\pi}^{\pi}f(u_{j_1},\lambda_1)  \phi^2_T(u_{j_1},\lambda_1)e^{-i(r-q)\lambda_1} \,d\lambda_1
\int_{-\pi}^{\pi}f(u_{j_1}, \lambda_2) e^{-i(r-q)\lambda_2} \,d\lambda_2   \\
&=&\frac{1}{N^2}\frac{1}{M^2} \sum_{j_1=1}^M\sum_{s=-(N-1)}^{N-1}  \int_{-\pi}^{\pi}f(u_{j_1},\lambda_1)  \phi^2_T(u_{j_1},\lambda_1)e^{-is\lambda_1} \,d\lambda_1
\int_{-\pi}^{\pi}f(u_{j_1}, \lambda_2) e^{-is\lambda_2} \,d\lambda_2  (N-|s|) \\
&=& D^{(1)}_{1,T}+D^{(2)}_{1,T}+D^{(3)}_{1,T} ,
\end{eqnarray*}
where
\begin{eqnarray*}
D^{(1)}_{1,T}&=&\frac{1}{N}\frac{1}{M^2} \sum_{j_1=1}^M\sum_{s=-\infty}^{\infty}  \int_{-\pi}^{\pi}f(u_{j_1},\lambda_1)  \phi_T^2(u_{j_1},\lambda_1)e^{-is\lambda_1} \,d\lambda_1
\int_{-\pi}^{\pi}f(u_{j_1},\lambda_2) e^{-is\lambda_2} \,d\lambda_2  \\
D^{(2)}_{1,T}&=&-\frac{1}{N}\frac{1}{M^2} \sum_{j_1=1}^M \sum_{\substack{s \in \mathbb{Z} \\   |s| \geq N}}^{} \int_{-\pi}^{\pi}f(u_{j_1},\lambda_1)  \phi_T^2(u_{j_1},\lambda_1)e^{-is\lambda_1} \,d\lambda_1
\int_{-\pi}^{\pi}f(u_{j_1},\lambda_2) e^{-is\lambda_2} \,d\lambda_2  \\
D^{(3)}_{1,T}  &=&-\frac{1}{N^2}\frac{1}{M^2} \sum_{j_1=1}^M\sum_{s=-(N-1)}^{N-1}   |s|\int_{-\pi}^{\pi}f(u_{j_1},\lambda_1)  \phi_T^2(u_{j_1},\lambda_1)e^{-is\lambda_1} \,d\lambda_1
\int_{-\pi}^{\pi}f(u_{j_1},\lambda_2) e^{-is\lambda_2} \,d\lambda_2  \\
\end{eqnarray*}
With Parseval's identity, we get
\begin{eqnarray*}
D^{(1)}_{1,T}&=&\frac{2\pi}{N}\frac{1}{M^2} \sum_{j_1=1}^M  \int_{-\pi}^{\pi}f^2(u_{j_1},\lambda_2) \phi_T^2(u_{j_1},\lambda_2) \,d\lambda_2  ,
\end{eqnarray*}
while Lemma \ref{fourier} in the online supplement yields (up to a constant) the inequalities
\begin{eqnarray*}
D^{(2)}_{1,T}& \leq &\frac{ g^2(k)}{N}\frac{1}{M^2} \sum_{j_1=1}^M\sum_{\substack{s \in \mathbb{Z} \\   |s| \geq N}}^{}  \frac{1}{|s|^{2-2\epsilon}} \leq \frac{ g^2(k)}{N^{2-2\epsilon}}\frac{1}{M}, \\
D^{(3)}_{1,T}& \leq &\frac{ g^2(k)}{N^2}\frac{1}{M^2} \sum_{j_1=1}^M \sum_{\substack{s \in \mathbb{Z} \\ 1 \leq |s| \leq N-1}}^{N-1} \frac{1}{|s|^{1-2\epsilon}} \leq \frac{ g^2(k)}{N^{2-2\epsilon}M},
\end{eqnarray*}
which proves (\ref{d1T}). We now consider the term 
$$
D_{2,T} =D^{(1)}_{2,T} + D^{(2)}_{2,T},
$$
where
\begin{eqnarray*}
D^{(1)}_{2,T} &=&   \frac{1}{(2 \pi N)^2}\frac{1}{M^2} \sum_{j_1=1}^M \sum_{\substack{q=0  }}^{N-1} \sum_{\substack{r=-\infty  }}^{\infty}\sum_{\substack{l,m,n,o=0 }}^{\infty}
\sum_{\substack{j_2=1 \\ 0 \leq q+o-m+(j_1-j_2)N \leq N-1 \\ |l-n+(j_2-j_1)N| \leq N-1 \\ |o-m+(j_1-j_2)N| \leq N-1}}^{\infty}
 \psi_l(u_{j_1}) \psi_m(u_{j_1})\psi_n(u_{j_2})  \psi_o(u_{j_2})\\
&&\int_{-\pi}^{\pi} \phi_T(u_{j_1},\lambda_1)e^{-i(r-q+l-n+(j_2-j_1)N)\lambda_1}
\int_{-\pi}^{\pi}\big[ \phi_T(u_{j_2},\lambda_2)-\phi_T(u_{j_2},\lambda_1) \big]e^{-i(r-q+m-o+(j_2-j_1)N)\lambda_2} \,d\lambda_2 \,d\lambda_1 \\
D^{(2)}_{2,T}   &= &- \frac{1}{(2 \pi N)^2}\frac{1}{M^2}\sum_{j_1=1}^M \sum_{\substack{q=0  }}^{N-1} \sum_{\substack{r \in \mathbb{Z}\\|r| \geq N }}^{}\sum_{\substack{l,m,n,o=0  }}^{\infty}
\sum_{\substack{j_2=1  \\ 0 \leq q+o-m+(j_1-j_2)N \leq N-1 \\ |l-n+(j_2-j_1)N| \leq N-1 \\ |o-m+(j_1-j_2)N| \leq N-1}}^{\infty}
\psi_l(u_{j_1}) \psi_m(u_{j_1})\psi_n(u_{j_2})  \psi_o(u_{j_2})\\
&&\int_{-\pi}^{\pi} \phi_T(u_{j_1},\lambda_1)e^{-i(r-q+l-n+(j_2-j_1)N)\lambda_1} \int_{-\pi}^{\pi}\big[ \phi_T(u_{j_2},\lambda_2)-\phi_T(u_{j_2},\lambda_1) \big]e^{-i(r-q+m-o+(j_2-j_1)N)\lambda_2} \,d\lambda_2  \,d\lambda_1 .
\end{eqnarray*}
Here $D^{(1)}_{2,T}$ corresponds to the sum over all $r$ and  vanishes by Parseval's identity. $D^{(2)}_{2,T}$ stands for the resulting error term which is of order   $O(T^{-1}  g^2(k) N^{-1+2D+2\epsilon})$ because of  Lemma \ref{Fehler Varianz} (vi) in the online supplement.
\\

Part b) follows with par a) if we show
\be
\cum_l[\sqrt{T}\hat G_T(\phi)]&=&O \big( g(k)^l T^{l(\epsilon-1/2+2D)+(1-4D)} \log(T)^l   \big) \quad \text{ for } l \geq 3 \mbox{ and }D < 1/4. \label{cum}
\ee
For a proof of this statement where we proceed (with a slight modification) analogously to the proof of Theorem 6.1 c) in \cite{prevet2012}. Note that these authors
work with functions $\phi_T$ such that
\begin{eqnarray}
\frac{1}{N} \sum_{k=1}^{N/2}\phi_T(u,\lambda_k)e^{ih\lambda_k}  &=& O\Big(\frac{1}{|h \text{ modulo } N/2|}\Big)
\label{fourier prevet}
\end{eqnarray}
while $\int _{-\pi}^{\pi}\phi_T(u,\lambda)e^{ih\lambda} d \lambda=O(h^{-1})$ for the integrated case. The authors then derive the exact same order as in \eqref{cum} with the only difference that $\epsilon=0$ and $g(k) \equiv 1$. In our situation, assumption (\ref{Ordnung theta}) and Lemma \ref{fourier} in the online supplement imply
\begin{eqnarray}
\int_{-\pi}^{\pi}\phi_T(u,\lambda)e^{ih\lambda}\,d\lambda   &=& O\Big(\frac{g(k)}{|h|^{1+2d_0(u)-\epsilon}}\Big) = O\Big(T^{\epsilon}\frac{g(k)}{|h|}\Big)
\label{fourier presendet}
\end{eqnarray}
and we can therefore proceed completely analogously to the proof of Theorem 6.1 c) in \cite{prevet2012} but using (\ref{fourier presendet}) instead of (\ref{fourier prevet}). The details are omitted for the sake of brevity. $\hfill \Box$

\bigskip

For the formulation of the next result we define the set
\begin{eqnarray*}
{\cal G}_T (s,\ell)   &=& \{ \tilde \phi_T: [-\pi,\pi] \rightarrow \mathbbm{R} ~| ~ \tilde \phi_T \mbox{ is symmetric, there exists a polynomial  }
P_\ell  \mbox{ of degree }  ~\ell  ~ \mbox{ and a } \\
&&  ~~~~~~~~~~~~~~~~~~~~~~~~~ \mbox{  constant } d \in [-\gamma_k,1/2) \mbox{ such that }  \tilde \phi_T(\lambda) = \log^s(|1-e^{i\lambda}|)|1-e^{i\lambda}|^{2d} |P_\ell (e^{i\lambda}) |^2 \}
\end{eqnarray*}
and state the following result.

\begin{satz} \label{uniformtheorem}
Suppose Assumption \ref{ass1} and \ref{ass3} are fulfilled, $N^{5/2}/T^{2}\rightarrow 0$  and  $0  < \epsilon <1/4-D/2$ is the constant of Assumption \ref{ass3}. Let $\Phi_T$ denote a class of functions
$\phi_T: [0,1] \times [-\pi,\pi] \rightarrow \mathbb{R}$
consisting of elements, which are  twice continuously  differentiable with uniformly bounded partial derivates with respect to $u,\lambda,T$ and satisfy \eqref{Ordnung theta}--\eqref{Ordnung theta 2.Ableitung} with $g(k) \equiv 1$, where the constant $C$ does not depend on $\Phi_T$, $T$. Furthermore, we assume that  for all $ u \in [0,1]$  the condition $ \phi_T(u,\cdot )  \in {\cal G}_T (s,qk)   $ holds,
 where  $q,s \in \mathbbm{N}$  are fixed and $k=k(T)$ denotes a sequence satisfying $k^4\log^2(T)N^{-\epsilon/2} \rightarrow 0$. Then
\bea
\sup_{u \in [0,1]} \sup_{\phi_T \in \Phi_T} \Big | \int_{-\pi}^\pi (I_N^{\mu}(u,\lambda)-f(\lfloor uT \rfloor /T,\lambda))\phi_T(\lfloor uT \rfloor /T,\lambda) d\lambda\Big |=o_P(N^{-1/2+\epsilon/2}).
\eea
\end{satz}
\textbf{Proof:}
We define $\Phi_{T}^*$ as the set of functions which we obtain by multiplying all elements $\phi_T \in \Phi_T$ with $1_{\{u=t/T\}}(u,\lambda)$,
that is $ \phi_T^*(u,\lambda) =  \phi_T (t/T,\lambda)$  for some $t=1,...,T$ and  $\phi_T \in \Phi_T$, and consider
\bea
\hat D_{T,1}(\phi_{T}^*):= \sum_{t_1=1}^T \int_{-\pi}^\pi I_N^{\mu}(t_1/T,\lambda)\phi_{T}^*(t_1/T,\lambda) d\lambda, \quad  \phi_{T}^* \in \Phi_{T}^*.
\eea
It follows from Theorem 2.1 in \cite{unifeconometrica} that the assertion
of Theorem \ref{uniformtheorem} is a consequence of the statements:
\begin{itemize}
\item[(i)] For every $\phi_{T}^*
\in \Phi_{T}^*$ we have
\begin{eqnarray}
\hat G_{T,1}(\phi_{T}^*):=N^{1/2-\epsilon/2} \Big(\hat D_{T,1}(\phi^*_{T})-\int_{-\pi}^{\pi} f(t/T,\lambda)\phi_T(t/T,\lambda) d\lambda \Big)=o_p(1)
\label{equi1}
\end{eqnarray}
\item[(ii)] For every $\eta>0$ we have
\begin{eqnarray}
\lim_{T \rightarrow \infty} P\bigl(\sup_{\phi^*_{T,1}, \phi^*_{T,2}  \in \Phi_{T}^*} |\hat G_{T,1}(\phi^*_{T,1}) - \hat G_{T,1}(\phi^*_{T,2})|> \eta \bigr) ~= 0.
\label{equi2}
\end{eqnarray}
\end{itemize}
In order to prove part (i) we use the same arguments as given in the proof of \eqref{ew} and \eqref{var} and obtain
\bea
\E [\hat D_{T,1}(\phi^*_{T}))] &=& \int_{-\pi}^{\pi} f(t/T,\lambda)\phi_T(t/T,\lambda) d\lambda+ O\Big(\frac{1}{N^{1-\epsilon-2\gamma_K}}\Big) + O \Big(\frac{N^2}{T^2} \Big),  \\
\Var[N^{1/2}\hat D_{T,1}(\phi^*_{T})]&=&   \int_{-\pi}^{\pi} f^2(t/T,\lambda)\phi_T^2(t/T,\lambda) d\lambda + O\Big(\frac{ 1}{N^{1-2D-2\epsilon-4\gamma_k}}\Big) +O \Big(\frac{ N^2}{T^2} \Big) ,
\eea
which yields (\ref{equi1}) observing the growth conditions on  $N$ and $T$. For the proof of part (ii) we note that it follows by similar arguments as given in the proof of Theorem 6.1 d) of \cite{prevet2012} that there exists a positive constant $C$  such that the inequlality
\bea
\E(|\hat G_{T,1}(\phi^*_{T,1})-\hat G_{T,1}(\phi^*_{T,2})|^l) \leq (2l)! C^l \Delta_{T,\epsilon}^l (\phi^*_{T,1},\phi^*_{T,2})
\eea
holds for all even $l\in \mathbb{N}$ and all $\phi^{*}_{T,1},\phi^{*}_{T,2}  \in \Phi_{T}^*$,
where
\bea
\Delta_{T,\epsilon}(\phi^*_{T,1},\phi^*_{T,2})= 1_{\{ t_1=t_2\}} N^{- \epsilon/2}  \sqrt{\int_{-\pi}^\pi (\phi_{T,1,1}(t_1/T,\lambda)-\phi_{T,1,2}(t_1/T,\lambda))^2 d\lambda } + A 1_{\{ t_1\not=t_2\}} N^{- \epsilon/2}
\eea
for a  constant $A$ which is sufficiently large such that
\begin{eqnarray*}
\sup_{\phi_{T,1,i} \in \Phi_{T}^*} \sqrt{\int_{-\pi}^\pi(\phi_{T,1,1}(t_1/T,\lambda)-\phi_{T,1,2}(t_1/T,\lambda))^2 d\lambda } \leq A.
\end{eqnarray*}
 By an application of Markov's inequality and a straightforward but cumbersome calculation [see the proof of Lemma 2.3 in \cite{dahlhaus1988} for more details] this yields
\bea
P(|\hat G_{T,1}(\phi^*_{T,1}) - \hat G_{T,1}(\phi^*_{T,2})|> \eta) \leq 96 \exp(-\sqrt{\eta \Delta_{T,\epsilon}^{-1}(\phi^*_{T,1},\phi^*_{T,2}) C^{-1}})
\eea

for all $\phi^*_{T,1},\phi^*_{T,2}  \in \Phi_{T}^*$. The statement (\ref{equi2}) then follows with the extension of the classical chaining argument as described in \cite{dahlhaus1988} if we show that the corresponding covering integral of $\Phi_{T}^*$ with respect to the semi-metric $\Delta_{T,\epsilon}$ is finite. More precisely, the covering number $N_T(u)$ of $\Phi_{T}^*$ with respect to $\Delta_{T,\epsilon}$ is equal to one for $u \geq A N^{- \epsilon/2}$ and bounded by $TC^{(qk)^2} u^{-qk}N^{-qk\epsilon/2}$ for some constant $C$ for $u<A N^{- \epsilon/2}$ [see Chapter VII.2. of \cite{pollard} for a definition of covering numbers]. This implies that the covering integral $J_T(\delta)=\int_0^\delta [\log(48 N_T(u)^2 u^{-1}]^2 du$ is up to a constant bounded by $k^4\log^2(T)N^{-\epsilon/2}$. The assertion
follows by  the assumptions on $k$ and $N$. $\hfill \Box$


\subsection{Proof of Theorem \ref{uniform}}

Introducing the notation
\begin{eqnarray*}
\mathcal{L}_{N,k}^{\mu}(\theta_k, u):= \frac{1}{4\pi} \int_{-\pi}^{\pi} \Big( \log(f_{\theta_k}(\lambda))+ \frac{I^{\mu}_N(u,\lambda)}{f_{\theta_k}(\lambda)} \Big) \,d\lambda, \quad u \in [0,1]
\end{eqnarray*}
we obtain with the same arguments as given in the proof of Theorem 3.6 in \cite{dahlhaus1997}
\bea
&& \max_{t=1, \ldots, T} \big | \mathcal{L}_{N,k}^{\hat \mu}(\theta_k, t/T) - \mathcal{L}_{N,k}^{\mu}(\theta_k, t/T) \big |  \\
&\leq & C \max_{t=1, \ldots, T} \max_{q=0, \ldots, N} \big\{ \big|\mu(t/T)-\hat{\mu}(t/T)\big|  \big | \int_{-\pi}^{\pi}
{d}^{X-\mu}_N(t/T,\lambda)   f^{-1}_{\theta_k}(\lambda) e^{iq\lambda}   \,d\lambda \big |  \big\} + CN^{\epsilon}\max_{t=1, \ldots, T}\big|\mu(t/T)-\hat{\mu}(t/T)\big|^2
\eea
for some constant $C \in \mathbb{R}$ and $d_N^{X-\mu}$ is defined by $\big|d^{X- \mu}_N(u,\lambda)\big|^2:=I_N^{\mu}(u,\lambda)$. By proceeding as in the proof of Theorem \ref{uniformtheorem} one verifies
\bea
\max_{t=1, \ldots, T} \max_{q=0,\ldots, N} \sup_{\theta_k \in \Theta_{R,k}}  \big | \int_{-\pi}^{\pi}
{d}^{X-\mu}_N(t/T,\lambda)   f^{-1}_{\theta_k}(\lambda) e^{iq\lambda}   \,d\lambda \big |=O(N^\epsilon),
\eea
and (\ref{convergencemu}) yields
\be
 \max_{t=1, \ldots, T} \sup_{\theta_k \in \Theta_{R,k}} \big | \mathcal{L}_{N,k}^{\hat \mu}(\theta_k, t/T) - \mathcal{L}_{N,k}^{\mu}(\theta_k, t/T) \big |
=  \max_{t=1, \ldots, T}\big|\mu(t/T)-\hat{\mu}(t/T)\big|O_p(N^{\epsilon})= o_p(k^{-5/2}) ,  \label{ratebeweisalt1}
\ee
and analogously we get
\be
\max_{t=1, \ldots, T} \sup_{\theta_k \in \Theta_{R,k}} \big \| \nabla \mathcal{L}_{N,k}^{\hat \mu}(\theta_k, t/T) - \nabla \mathcal{L}_{N,k}^{\mu}(\theta_k, t/T) \big \|_2 
= \max_{t=1, \ldots, T}\big|\mu(t/T)-\hat{\mu}(t/T)\big|O_p(k^{1/2}N^{\epsilon})= o_p(k^{-5/2}). ~~~~\label{ratebeweisalt2}
\ee
For each $u \in [0,1]$ let $\hat \theta_{N,k}(u)$ denote the Whittle-estimator defined in \eqref{Whittleestimator}. Then Theorem \ref{uniformtheorem} and similar arguments as in the proof of Theorem 3.2 in \cite{dahlhaus1997} yield
\be \label{konvergenzestim}
\sup_{u \in [0,1]} \big \|  \hat{\theta}_{N,k}(u)-  \theta_{0,k}(u)\big \|_2  &=& o_p(1).
\ee
We will now derive a refinement of this statement.
By an application of the mean value theorem, there exist  vectors $\zeta_{u}^{(k)}=(\zeta_{u,1}^{(k)}, \zeta_{u,2 }^{(k)}, \ldots, \zeta_{u,k+1}^{(k)}) \in \mathbb{R}^{k+1}$, $u \in \{1/T, 2/T, \ldots, 1 \}$, satisfying  $\| \zeta_{u}^{(k)}  - \theta_{0,k}(u)\|_2 \leq \| \hat{\theta}_{N,k}(u)-  \theta_{0,k}(u) \|_2$  such that
\begin{eqnarray*}
\nabla \mathcal{L}_{N,k}^{\hat \mu}(\hat{\theta}_{N,k}(u), u )-\nabla \mathcal{L}_{N,k}^{\hat \mu}( \theta_{0,k}(u), u)
&=& \nabla^2 \mathcal{L}_{N,k}^{\hat{\mu}}(\zeta_{u}^{(k)}, u) \big(  \hat{\theta}_{N,k}(u)-  \theta_{0,k}(u)\big),
\end{eqnarray*}
and the first term on the left-hand side vanishes due to \eqref{konvergenzestim}. This yields
\begin{eqnarray*}
E_T-\nabla \mathcal{L}_{N,k}^{\mu}( \theta_{0,k}(u), u)
&=& \nabla^2 \mathcal{L}_{N,k}^{\hat \mu}(\zeta_{u}^{(k)}, u) \big(  \hat{\theta}_{N,k}(u)-  \theta_{0,k}(u)\big) ,
\end{eqnarray*}
where $E_T$ denotes the difference between $\nabla \mathcal{L}_{N,k}^{\mu}( \theta_{0,k}(u), u)$ and $\nabla \mathcal{L}_{N,k}^{\hat \mu}( \theta_{0,k}(u), u)$, which is of order \linebreak $\max_{t=1, \ldots, T}\big|\mu(t/T)-\hat{\mu}(t/T)\big|O_p(k^{1/2}N^{\epsilon})$ by \eqref{ratebeweisalt2}. It follows from
\begin{eqnarray*}
\nabla \mathcal{L}_{N,k}^{\mu}( \theta_k,u) &=&\frac{1}{4\pi} \int_{-\pi}^{\pi} \big[I^{\mu}_N(u,\lambda)-f_{\theta_k}(\lambda) \big] \nabla f_{\theta_k}^{-1}(\lambda) \,d\lambda
\label{Log-Likeli 1. Ableitung}
\end{eqnarray*}
and Theorem \ref{uniformtheorem} that  $\max_{u \in \{1/T, \dots 1 \}} \|\nabla \mathcal{L}_{N,k}^{\mu}( \theta_{0,k}(u), u)\|_2=O_p(\sqrt{k} N^{-1/2+\epsilon/2})$ so it remains to show
\bea
P(\nabla^2 \mathcal{L}_{N,k}^{\hat \mu}(\zeta_{u}^{(k)}, u)^{-1} \text{ exists and } \|\nabla^2 \mathcal{L}_{N,k}^{\hat \mu}(\zeta_{u}^{(k)}, u)^{-1} \|_{sp} \leq C k \text{ for all } u \in \{1/T, \ldots, 1 \}) \rightarrow 1
\eea

for some positive constant $C$. This, however, follows with a Taylor expansion, \eqref{konvergenzestim}, Theorem \ref{uniformtheorem} and Assumption \ref{ass2} (iv) for the corresponding expression with $\hat \mu$  replaced by $\mu$. The more general case is then implied by the convergence-assumptions on $\hat \mu$. $\hfill \Box$

\subsection{Proof of Theorem \ref{asymp teststat} and Theorem \ref{asymp teststatalt}}
We will show in Section \ref{sec33} that under the null hypothesis  $H_0$ the estimate
\be
\max_{j=1,\ldots,M}\big\|  \hat{\theta}_{N,k}(u_j)-  \theta_{0,k}(u_j)\big\|_2  &=& O_p(k^{3/2}N^{-1/2+\epsilon/2})
\label{null 1}
\ee
is valid, while Theorem \ref{uniform} and \eqref{alt 2} imply
\be
k^{3/2}\max_{j=1,\ldots,M}\big\|  \hat{\theta}_{N,k}(u_j)-  \theta_{0,k}(u_j)\big\|_2  &=& o_p(1) \label{alt 1}
\ee
under the alternative $H_1$. As in the proof of Theorem \ref{uniform} there exist vectors $\zeta_{j}^{(k)}=(\zeta_{j,1}^{(k)}, \zeta_{j,2 }^{(k)}, \ldots, \zeta_{j,k+1}^{(k)}) \in \mathbb{R}^{k+1}$, $j=1, \ldots, M$, satisfying  $\| \zeta_{j}^{(k)}  - \theta_{0,k}(u_j)\|_2 \leq \| \hat{\theta}_{N,k}(u_j)-  \theta_{0,k}(u_j) \|_2$  such that
\begin{eqnarray*}
-\nabla \mathcal{L}_{N,k}^{\hat \mu}( \theta_{0,k}(u_j), u_j)
&=& \nabla^2 \mathcal{L}_{N,k}^{\hat{\mu}}(\zeta_{j}^{(k)}, u_j) \big(  \hat{\theta}_{N,k}(u_j)-  \theta_{0,k}(u_j)\big)
\end{eqnarray*}
holds because of Assumption \ref{ass2} (ii) and \eqref{null 1} (under $H_0$) or \eqref{alt 1} (under $H_1$). By rearranging and summing over every block, it follows that
\begin{eqnarray} \label{summ1}
\frac{1}{M}\sum_{j=1}^{M}\big(  \hat{\theta}_{N,k}(u_j)-  \theta_{0,k}(u_j)\big)&=&
R_{0,T}-R_{1,T}-R_{2,T}-R_{3,T}-R_{4,T}
\end{eqnarray}
where
\begin{eqnarray*}
R_{0,T}&:=&-\frac{1}{M}\sum_{j=1}^{M}\Gamma_k^{-1}(\theta_{0,k}(u_j))\nabla \mathcal{L}_{N,k}^{\mu}( \theta_{0,k}(u_j), u_j) ,
\end{eqnarray*}
$\Gamma_k^{-1}$ is defined in (\ref{gamk}) and the terms $R_{i,T}   (i=1\ldots,4)$ are given by
\begin{eqnarray*}
R_{1,T}&:=& \frac{1}{M}\sum_{j=1}^{M}\Gamma_k^{-1}(\theta_{0,k}(u_j))\big(\nabla \mathcal{L}_{N,k}^{\hat{\mu}}( \theta_{0,k}(u_j), u_j)-\nabla \mathcal{L}_{N,k}^\mu( \theta_{0,k}(u_j), u_j) \big),\\
R_{2,T}&:=&\frac{1}{M}\sum_{j=1}^{M}\Gamma_k^{-1}(\theta_{0,k}(u_j))\big(\nabla^2 \mathcal{L}_{N,k}^{\hat{\mu}}(\zeta_{j}^{(k)},  u_j)-\nabla^2 \mathcal{L}_{N,k}^{\mu}(\zeta_{j}^{(k)}, u_j)\big) \big(  \hat{\theta}_{N,k}(u_j)-  \theta_{0,k}(u_j)\big) , \\
R_{3,T}&:=& \frac{1}{M}\sum_{j=1}^{M}\Gamma_k^{-1}(\theta_{0,k}(u_j))\big(\nabla^2 \mathcal{L}_{N,k}^{ \mu}(\zeta_{j}^{(k)}, u_j)-\nabla^2 \mathcal{L}_{N,k}^{\mu}(\theta_{0,k}(u_j),  u_j)\big)  \big(  \hat{\theta}_{N,k}(u_j)-  \theta_{0,k}(u_j)\big),\\
R_{4,T}&:=&\frac{1}{M}\sum_{j=1}^{M}\Gamma_k^{-1}(\theta_{0,k}(u_j))\big(\nabla^2 \mathcal{L}_{N,k}^{\mu}(\theta_0(u_j), u_j)-\Gamma_k(\theta_{0,k}(u_j))\big) \big(  \hat{\theta}_{N,k}(u_j)-  \theta_{0,k}(u_j)\big).
\end{eqnarray*}
We obtain for the first summand in \eqref{summ1}
\begin{eqnarray*}
R_{0,T}=
-\frac{1}{M}\sum_{j=1}^{M}\frac{1}{4\pi} \int_{-\pi}^{\pi} \big[I^{\mu}_N(u_j,\lambda)-f_{\theta_{0,k}(u_j)}(\lambda) \big]  \Gamma_k^{-1}(\theta_{0,k}(u_j))\nabla f^{-1}_{\theta_{0,k}(u_j)}(\lambda) \,d\lambda
\end{eqnarray*}
and with the notation $\phi_T(u_j, \lambda)=1/(4\pi) [\Gamma_k^{-1}(\theta_{0,k}(u_j))\nabla f^{-1}_{\theta_{0,k}(u_j)}(\lambda)]_{1}$,  it is easy to see that Assumption \ref{ass2} (i)--(iv) imply the conditions of Theorem \ref{zgws} b) with $g(k)=k^2$.
Moreover, observing the definition of $V_T$ and $W_T$ in Theorem \ref{zgws} and \ref{asymp teststat}, (\ref{lim1}) yields $V_T /W_T \rightarrow 1$. Consequently, under the assumptions of Theorem \ref{asymp teststat} it follows (observing \eqref{null 4 implication} and the growth conditions on $N$, $T$)
\begin{eqnarray*}
\frac{\sqrt{T}}{M}\sum_{j=1}^{M}\big[\Gamma_k^{-1}(\theta_{0,k}(u_j))\nabla \mathcal{L}_N^\mu( \theta_{0,k}(u_j), u_j)\big]_{1}/\sqrt{W_T}\stackrel{\mathcal{D}}{\rightarrow} \mathcal{N}(0, 1).
\end{eqnarray*}
Since $d_{0}(u)$ is the first element of the vector $\theta_{0,k}(u)$, Theorem \ref{asymp teststat} is a consequence of the fact
$
\frac{1}{M} \sum_{j=1}^M d_0(u_j)=F + O(M^{-2})
$
[this can be proved by a second order Taylor expansion] if we are able to show that
\begin{eqnarray*}
R_{i,T}=o_p(T^{-1/2}); \quad i=1, \ldots, 4.
\end{eqnarray*}
Analogously, Theorem \ref{asymp teststatalt} follows from \eqref{ew} and \eqref{var} if the estimates
\begin{eqnarray*}
R_{i,T}=o_p(1) \quad i=1, \ldots, 4.
\end{eqnarray*}
can be established.
It can be shown analogously to the proof of Theorem 3.6 in \cite{dahlhaus1997}, that, under assumptions \eqref{null 2} -- \eqref{null 3},  both terms $R_{1,T}$ and $R_{2,T}$ are of order $O_p(k^2N^{-\epsilon}T^{-1/2}+k^2N^{\epsilon-1})$, while, under assumption  \eqref{alt 2}, the order is $o_p(1)$ [see the proof of \eqref{ratebeweisnull1} and \eqref{ratebeweisalt1}, respectively, for more details].
Therefore it only remains to consider the quantities $R_{3,T}$ and $R_{4,T}$. For this purpose note that
\begin{eqnarray}
\nabla^2 \mathcal{L}_{N,k}^{\mu}( \theta_k(u_j),u_j)&=&\frac{1}{4\pi} \int_{-\pi}^{\pi} \big[I^{\mu}_N(u_j,\lambda)-f_{\theta_k(u_j)}(\lambda)\big]  \nabla^2 f_{\theta_k(u_j)}^{-1}(\lambda)\,d\lambda  +\Gamma_k(\theta_k(u_j))
\label{Log-Likeli 2. Ableitung}
\\
\nabla^3 \mathcal{L}_{N,k}^\mu( \theta_k(u_j),u_j)&=&\frac{1}{4\pi} \int_{-\pi}^{\pi} \big[ I^{\mu}_N(u_j,\lambda) - f_{\theta_k(u_j)}(\lambda)\big] \bigg[\frac{\partial^3 f_{\theta_k(u_j)}^{-1}(\lambda)}{\partial \theta_{j,t} \partial \theta_{j,s} \partial \theta_{j,r}} \bigg]_{r,s,t=1,\ldots, k+1}\,d\lambda
\notag\\
&&- \frac{1}{4\pi} \int_{-\pi}^{\pi} \bigg[ \frac{\partial f_{\theta_k(u_j)}(\lambda) }{\partial \theta_{j,t}}\frac{\partial^2 f_{\theta_k(u_j)}^{-1}(\lambda)}{ \partial \theta_{j,s} \partial \theta_{j,r}} \bigg]_{r,s,t=1,\ldots, k+1}\,d\lambda
\notag\\
&&+\frac{1}{4\pi} \int_{-\pi}^{\pi}\bigg[ \frac{\partial }{\partial \theta_{j,t}}\bigg(\frac{\partial f_{\theta_k(u_j)}(\lambda)}{\partial \theta_{j,s}}\frac{1}{f_{\theta_k(u_j)}^2(\lambda)} \frac{\partial f_{\theta_k(u_j)}(\lambda)}{\partial \theta_{j,r}} \bigg)\bigg]_{r,s,t=1,\ldots, k+1}\,d\lambda ,
\label{Log-Likeli 3. Ableitung}
\end{eqnarray}
where we used the notation $(\theta_{j,1},\theta_{j,2}, \ldots, \theta_{j,k+1} ):=(d(u_j), a_{1}(u_j), \ldots, a_{k}(u_j))$. For the term $R_{3,T}$ we obtain with the well-known inequality $\| Ax\|_2 \leq \| A\|_{sp} \| x\|_2$
\begin{eqnarray*}
\|R_{3,T} \|_2& \leq &   \max_{\theta_k \in \Theta_{R,k}}\big\|\Gamma_k^{-1}(\theta_k)\big\|_{sp}\frac{1}{M}\sum_{j=1}^{M} \big\|\nabla^2 \mathcal{L}_{N,k}^\mu(\zeta_{j}^{(k)},u_j) - \nabla^2 \mathcal{L}_{N,k}^\mu(\theta_{0,k}(u_j),u_j)\big\|_{sp}\big\|  \hat{\theta}_{N,k}(u_j)-  \theta_{0,k}(u_j)\big\|_2 .
\end{eqnarray*}
By the mean value theorem there exist vectors $\tilde{\zeta}_{j}^{(k)} \in \mathbb{R}^k$ such that
\begin{eqnarray*}
  &&\big\|\nabla^2 \mathcal{L}_{N,k}^\mu(\zeta_{j}^k,u_j) - \nabla^2 \mathcal{L}_{N,k}^\mu(\theta_0(u_j),u_j)\big\|_{sp} \leq      k   \max_{r,s=1,\ldots,k }\big |  \big [ \nabla^2 \mathcal{L}_{N,k}^\mu(\zeta_{j}^{(k)},u_j) - \nabla^2 \mathcal{L}_{N,k}^\mu(\theta_{0,k}(u_j),u_j)\big]_{r,s}\big| \\
& = &     k    \max_{r,s=1,\ldots,k }\big |  \nabla\big [ \nabla^2 \mathcal{L}_{N,k}^\mu(\tilde{\zeta}_{j}^{(k)},u_j) \big]_{r,s}  \big(\zeta_{j}^{(k)}-\theta_{0,k}(u_j) \big)\big|  \leq     k   \max_{r,s=1,\ldots,k }\big\| \nabla\big [ \nabla^2 \mathcal{L}_{N,k}^\mu(\tilde{\zeta}_{j}^{(k)},u_j) \big]_{r,s}  \big\|_2 \big\|\zeta_{j}^{(k)}-\theta_{0,k}(u_j) \big\|_2  \\
& \leq &    k  \big\|  \hat{\theta}_{N,k}(u_j)-  \theta_{0,k}(u_j)\big\|_2   \sup_{\substack{\theta_k \in \Theta_{R,k} \\r,s=1,\ldots,k} }\big\| \nabla\big [ \nabla^2 \mathcal{L}_{N,k}^\mu(\theta_k,u_j) \big]_{r,s}  \big\|_2 ,
\end{eqnarray*}
where $\| \tilde{\zeta}_{j}^{(k)}  - \theta_{0,k}(u_j)\|_2 \leq \| \zeta_{j}^{(k)}-  \theta_{0,k}(u_j) \|_2$ for every $j=1,...,M$. Therefore, we obtain
\begin{eqnarray*}
\|R_{3,T} \|_2&\leq& k   \max_{j=1,\ldots,M}\big\|  \hat{\theta}_{N,k}(u_j)-  \theta_{0,k}(u_j)\big\|^2_2  \sup_{\theta_k \in \Theta_{R,k}}\big\|\Gamma_k^{-1}(\theta_k)\big\|_{sp}       \sup_{\substack{\theta_k \in \Theta_{R,k}; j=1, \ldots, M\\r,s=1,\ldots k} }\big\| \nabla\big [ \nabla^2 \mathcal{L}_N^\mu(\theta_k,u_j) \big]_{r,s}  \big\|_2 \\
& \leq &  k C \max_{j=1,\ldots,M}\big\|  \hat{\theta}_{N,k}(u_j)-  \theta_{0,k}(u_j)\big\|^2_2   \sup_{\theta_k \in \Theta_{R,k}}\big\|\Gamma_k^{-1}(\theta_k)\big\|_{sp}\\
      && \Big( k \cdot\sup_{\substack{\theta_k\in \Theta_{R,k};j=1,...,M\\r,s,t=1,\ldots,k} }\Big|\frac{1}{4\pi} \int_{-\pi}^{\pi} \big[ I^{\mu}_N(u_j,\lambda) - f_{\theta_k}(\lambda)\big]\frac{\partial^3 f^{-1}_{\theta_k}(\lambda)}{\partial \theta_{j,t} \partial \theta_{j,s} \partial \theta_{j,r}} \,d\lambda \Big| +  k \Big),
\end{eqnarray*}
where, in the last inequality, we have used the fact that the second and third term in (\ref{Log-Likeli 3. Ableitung}) are bounded by a constant [this follows directly from Assumption \ref{ass2}]. Before we investigate the order of this expression, we derive a similar bound for the term  $R_{4,T}$. Observing (\ref{Log-Likeli 2. Ableitung}) we obtain
\begin{eqnarray*}
\| R_{4,T}\|_2&\leq&  \max_{j=1,\ldots,M} \big\|  \hat{\theta}_{N,k}(u_{j})-  \theta_{0,k}(u_j)\big\|_2 \sup_{\theta_k \in \Theta_{R,k}}\big\|\Gamma_k^{-1}(\theta_k)\big\|_{sp} \max_{j=1,\ldots,M}\big\|\nabla^2 \mathcal{L}_{N,k}^\mu(\theta_{0,k}(u_j), u_j)-\Gamma_k(\theta_{0,k}(u_j))\big\|_{sp} \\
&=& \max_{j=1,\ldots,M} \big\|  \hat{\theta}_{N,k}(u_{j})-  \theta_{0,k}(u_j)\big\|_2  \sup_{\theta_k \in \Theta_{R,k}}\big\|\Gamma_k^{-1}(\theta_k)\big\|_{sp} \\ && \phantom{......} \times \max_{j=1,\ldots,M}\big\|\frac{1}{4\pi} \int_{-\pi}^{\pi}
\big[I^{\mu}_N(u_j,\lambda)-f_{\theta_{0,k}(u_j)}(\lambda) \big] \nabla^2 f^{-1}_{\theta_{0,k}(u_j)}(\lambda) \,d\lambda \big\|_{sp} \\
& \leq & k \max_{j=1,\ldots,M} \big\|  \hat{\theta}_{N,k}(u_{j})-  \theta_{0,k}(u_j)\big\|_2  \sup_{\theta_k \in \Theta_{R,k}}\big\|\Gamma_k^{-1}(\theta_k)\big\|_{sp} \\
&& \phantom{......} \times \max_{j=1,\ldots,M}\max_{r,s=1,...,k}\Big |\frac{1}{4\pi} \int_{-\pi}^{\pi}
\big[I^{\mu}_N(u_j,\lambda)-f_{\theta_{0,k}(u_j)}(\lambda) \big] \frac{\partial^2 f^{-1}_{\theta_{0,k}(u_j)}(\lambda)}{ \partial \theta_{j,s} \partial \theta_{j,r}} \,d\lambda \Big |.
\end{eqnarray*}
If we show
\begin{eqnarray*}
\max_{j=1,\ldots,M} \sup_{\substack{\theta_k \in \Theta_{R,k}\\r,s,t=1,\ldots,k} }\Big|\frac{1}{4\pi} \int_{-\pi}^{\pi} \big[ I^{\mu}_N(u_j,\lambda) - f_{\theta_k}(\lambda)\big] \frac{\partial^3 f^{-1}_{\theta_k}(\lambda)}{\partial \theta_{j,t} \partial \theta_{j,s} \partial \theta_{j,r}} \,d\lambda \Big| &=& O_p(1), \\
\max_{j=1,\ldots,M}\max_{r,s=1,...,k}\Big |\frac{1}{4\pi} \int_{-\pi}^{\pi}
\Big[I^{\mu}_N(u_j,\lambda)-f_{\theta_{0,k}(u_j)}(\lambda) \Big] \frac{\partial^2 f^{-1}_{\theta_{0,k}(u_j)}(\lambda)}{ \partial \theta_{j,s} \partial \theta_{j,r}} \,d\lambda \Big | &=& O_p(N^{-1/2+\epsilon/2}),
\end{eqnarray*}

it follows with Assumption \ref{ass2} (iv) in combination with  (\ref{null 1}) (under $H_0$) and \eqref{alt 1} (under $H_1$)  that the terms
$R_{3,T}$ and $R_{4,T}$ are of order $o_p(T^{-1/2})$ (under $H_0$) and $o_p(1)$ (under $H_1$). These two claims, however are a direct consequence of Theorem \ref{uniformtheorem} and  \eqref{null 4 implication}. $\hfill \Box$

\subsubsection{Proof of (\ref{null 1})} \label{sec33}

With the   same arguments as in the proof of Theorem 3.6 in \cite{dahlhaus1997} we obtain
\bea
\max_{j=1,\ldots,M} \big | \mathcal{L}_{N,k}^{\hat \mu}(\theta_k, u_j) - \mathcal{L}_{N,k}^{\mu}(\theta_k, u_j) \big | \leq \Pi_{1,T} + \Pi_{2,T},
\eea
where
\bea
\Pi_{1,T}&=& C \max_{t=1, \ldots, T} \max_{q=1,...,N} \sup_{\theta_k \in \Theta_{R,k}} \Big| \int_{-\pi}^{\pi}
d^{X-\mu}_N(t/T,\lambda)  f^{-1}_{\theta_{k}}(\lambda) \sum_{s=0}^{q-1}e^{is\lambda} \,d\lambda \Big| \\
&& \times \Big(
 \max_{t=1, \ldots, T}\Big|\Big\{\mu\Big(\frac{t-1}{T}\Big)-\hat{\mu}\Big(\frac{t-1}{T}\Big)\Big\}-\Big\{{\mu}\Big(\frac{t}{T}\Big)-\hat{\mu}\Big(\frac{t}{T}\Big)\Big\}\Big|  +\max_{t=1, \ldots, T}\big|\mu(t/T)-\hat{\mu}(t/T)\big|/N \Big) \\
\Pi_{2,T} &=& C N^{\epsilon} \max_{t=1, \ldots, T}\big|\mu(t/T)-\hat{\mu}(t/T)\big|^2
\eea
and $C$ denotes a positive constant. By proceeding as in the proof of Theorem \ref{uniformtheorem} one obtains
\bea
\max_{t=1, \ldots, T} \max_{q=1,...,N} \sup_{\theta_k \in \Theta_{R,k}} \Big| \int_{-\pi}^{\pi}
d^{X-\mu}_N(t/T,\lambda)  f^{-1}_{\theta_{k}}(\lambda) \sum_{s=0}^{q-1}e^{is\lambda} \,d\lambda \Big| = o(N^{1/2+\epsilon/2}),
\eea
which implies (observing the assumptions \eqref{null 2} and \eqref{null 3})
\be \label{ratebeweisnull1}
\max_{j=1,\ldots,M} \sup_{\theta_k \in \Theta_{R,k}} \big | \mathcal{L}_{N,k}^{\hat \mu}(\theta_k, u_j) - \mathcal{L}_{N,k}^{\mu}(\theta_k, u_j) \big |&= O_p(N^{-\epsilon}T^{-1/2}+N^{\epsilon-1}) &= o_P(N^{-1/2+\epsilon/2}k^{1/2}) \phantom{.......}
\ee
under $H_0$. Analogously we obtain
\be
&& \max_{j=1,\ldots,M} \sup_{\theta_k \in \Theta_{R,k}} \big \| \nabla \mathcal{L}_{N,k}^{\hat \mu}(\theta_k, u_j) - \nabla \mathcal{L}_{N,k}^{\mu}(\theta_k, u_j) \big \|_2 \nonumber \\
&=& O_p(k^{1/2}N^{-\epsilon}T^{-1/2}+k^{1/2}N^{\epsilon-1}) = o_P(N^{-1/2+\epsilon/2}k^{1/2})   \label{ratebeweisnull2}
\ee
under the null hypothesis. By using \eqref{ratebeweisnull1} and \eqref{ratebeweisnull2} instead of \eqref{ratebeweisalt1} and \eqref{ratebeweisalt2}, assertion  (\ref{null 1}) follows by the same arguments as given in the proof of Theorem \ref{uniform}.  $\hfill \Box$

\subsection{Proof of Theorem \ref{mean}}

A second order Taylor expansion yields
\be
\E(\hat \mu_L (t/T))
&=& \mu(t/T)+\frac{\mu'(t/T)}{L}\sum_{p=0}^{L-1}(-L/2+1+p)/T+O(L^2/T^2)
 = \mu(t/T)+O(1/T+L^2/T^2). \label{biasorder}
\ee

For $t_i \in \{1,...,T\}$ the cumulants of order $l \geq 2$
\bea
\cum(\hat \mu_L(t_1/T), \hat \mu_L(t_2/T),...,\hat \mu_L(t_l/T))= \frac{1}{L^l} \sum_{p_1,...,p_l=0}^{L-1} \sum_{m_1,...,m_l=0}^\infty \psi_{t_1,T,m_1} \cdots \psi_{t_l,T,m_l} \cum(Z_{p_1-m_1},...,Z_{p_l-m_l})
\eea
are bounded by
\bea
\frac{C}{L^l} \sum_{p_1=0}^{L-1}\sum_{\substack{m_1,...,m_l=0 \\ |m_i-m_{i+1}| \leq L}}^\infty \frac{1}{(I(m_1 \cdots m_l))^{1-D}} \leq C^l L^{1-l(1-D)} ,
\eea
where we used
the independence of the innovations, \eqref{apprbed} and \eqref{apprpsi} and the last inequality follows by replacing the sums by its corresponding approximating integrals and holds for some positive constant $C$ (which is independent of $l$ and may vary in the following arguments). This yields that $\hat \mu_L(t/T)$ estimates its true counterpart at a pointwise rate of $L^{1/2-D}$ and  we now continue by showing stochastic equicontinuity. The expansion \eqref{biasorder} and the bound $C^l L^{1-l(1-D)} $ for the $l$-th cumulant ($l \geq 2$) of $\hat \mu_L$ yield $\cum_l(L^{1/2-D-\alpha/2}(\hat \mu_L(t_1/T)-\hat \mu_L(t_2/T))) \leq (2C)^l L^{-l \alpha /2}$ for all $t_i \in \{1,...,T\}$ and every $\alpha>0$, from which we get
\bea
\E(L^{l(1/2-D-\alpha)}(\hat \mu_L(t_1/T)-\hat \mu_L(t_2/T))^l) \leq  (2l)! C^l L^{-l \alpha /2}  \quad \text{ for all even } l \in \mathbb{N} \text{ and } t_i \in \{1,...,T\}
\eea
[see the proof of Lemma 2.3 in \cite{dahlhaus1988} for more details]. By considering the order of the bias \eqref{biasorder} this yields
\bea
L^{1/2-D-\alpha}\max_{t=1, \ldots, T} \big|\mu(t/T)-\hat{\mu}_L(t/T)\big|=o_p(1), \quad \text{ for every } \alpha > 0,
\eea
as in the proof of Theorem \ref{uniformtheorem}. Consequently \eqref{null 2} [under the conditions of part a)] and \eqref{alt 2} [under the conditions of part b)] follow. So it remains to show \eqref{null 3} in the case $D=0$. For this purpose we define
\bea
\Delta (t/T)=\Big\{\mu\Big(\frac{t-1}{T}\Big)-\hat{\mu}_L\Big(\frac{t-1}{T}\Big)\Big\}-\Big\{{\mu}\Big(\frac{t}{T}\Big)-\hat{\mu}_L\Big(\frac{t}{T}\Big)\Big\},
\eea
and from (\ref{biasorder}) we obtain $\E(\Delta (t/T))=O(T^{-1}+L^2/T^2)$. A simple calculation reveals $\cum(\Delta (t_1/T),\Delta (t_2/T))=O(L^{-1}T^{-1})$ (where the estimate is independent of $t_i$) and with the Gaussianity of the innovations we get $\cum(\Delta (t_1/T),...,\Delta (t_l/T))=0$ for $l \geq 3$. This yields, as above,
\bea
L^{1/2-\alpha}T^{1/2}\max_{t=1, \ldots, T} | \Delta(t/T) | =o_p(1)
\eea

for every $\alpha>0$, and completes the proof of Theorem  \ref{mean}. $\hfill \Box$

\newpage

\section{Online supplement: Auxiliary results}
\label{aux results}
\def\theequation{7.\arabic{equation}}
\setcounter{equation}{0}

Finally, we state some lemmas which were employed  in the above proofs.

\begin{lem}
\label{Lemma 7.1}
Suppose it is $\mu, \nu, a, b \in \mathbb{R}$. Then there exists a constant $C \in \mathbb{R}$ such that the following holds:
\begin{itemize}
\item[a)] If $\mu , \nu >0$ and $b > a$, then
\begin{eqnarray}
\label{Gleichung 3.196(3)}
\sum_{\substack{p=0 \\p-a \geq 1 \\ -p + b \geq 1 }}^{N-1} \frac{1}{(p-a)^{1-\mu}} \frac{1}{(b-p)^{1-\nu}}
\leq \sum_{p=1+a}^{b-1} \frac{1}{(p-a)^{1-\mu}} \frac{1}{(b-p)^{1-\nu}}\leq  \frac{C}{(b-a)^{1-\mu-\nu}}.
\end{eqnarray}
\item[b)] If $0 <\mu,\nu$ and $0 < 1-\mu-\nu $, then it follows for $|a+b| > 0$
\begin{eqnarray}
\sum_{\substack{p=1 \\ p+b \geq 1 \\ p-a \geq 1}}^{N-1} \frac{1}{(p+b)^{1-\mu}} \frac{1}{(p-a)^{1-\nu}} \leq
\sum_{\substack{p=1 \\ p+b \geq 1 \\ p-a \geq 1}}^{\infty} \frac{1}{(p+b)^{1-\mu}} \frac{1}{(p-a)^{1-\nu}}
\leq \frac{C}{|a+b|^{1-\mu-\nu}}.
\label{eq1}
\end{eqnarray}
\item[c)] If $0 < \nu < 1-\mu$ and $y,z \geq 1$, then
\begin{eqnarray*}
&&\sum_{p=1+y}^{\infty} \frac{\log (p)}{p^{1-\mu}} \frac{1}{(p-y)^{1-\nu}}\leq  \frac{C\log(y)}{y^{1-\mu-\nu}}, \\
&&\sum_{p=1}^{\infty} \frac{\log(p+z)}{(p+z)^{1-\mu}} \frac{1}{p^{1-\nu}} \leq  \frac{C\log(z)}{z^{1-\mu-\nu}}.
\end{eqnarray*}
\end{itemize}
\end{lem}
\textbf{Proof:} The proof can be found in \cite{senpreudet2013}. $\hfill \Box$

\begin{lem}
\label{fourier}
For every $T \in \mathbb{N}$, let  $\eta_{T}: [-\pi, \pi] \mapsto \mathbb{R}$ be a symmetric and twice continuously differentiable function such that
$\eta_{T} = O( |\lambda|^{\alpha})$  for some $\alpha \in (-1,1)$ as $|\lambda| \rightarrow 0$ (where the constant in the $O(\cdot)$ term is independent of $T$). Then,   for $|h| \rightarrow \infty$, we have
\begin{eqnarray*}
\int_{-\pi}^{\pi}\eta_T( \lambda)e^{ih\lambda}\,d\lambda & = &O\Big( \frac{1}{|h|^{1+\alpha} }\Big)
\end{eqnarray*}
uniformly in $T$.
\end{lem}
\textbf{Proof: }
The assertion  follows from Lemma 4 and Lemma 5 in \cite{foxtaqqu1986}. $\hfill \Box$

\begin{lem}
\label{Fehler Varianz}
If Assumption \ref{ass1} holds, then
\begin{itemize}
\item [(i)]
\noindent
\vspace{-1.1cm}
\begin{eqnarray*}
&&  \frac{1}{ N^2}\frac{1}{M^2} \sum_{j_1=1}^M \sum_{\substack{q,r=0  }}^{N-1}\sum_{\substack{l,m,n,o=0  }}^{\infty} \sum_{\substack{j_2=1 \\ N \leq |r+l-n+(j_2-j_1)N| \\ 0 \leq q+o-m+(j_1-j_2)N \leq N-1 \\ |l-n+(j_2-j_1)N| \leq N-1 \\ |o-m+(j_1-j_2)N| \leq N-1 }}^M \psi_l(u_{j_1}) \psi_m(u_{j_1})\psi_n(u_{j_2})  \psi_o(u_{j_2})\\
&&\int_{-\pi}^{\pi} \phi_T(u_{j_1},\lambda_1)e^{-i(r-q+l-n+(j_2-j_1)N)\lambda_1} \,d\lambda_1
\int_{-\pi}^{\pi} \phi_T(u_{j_2},\lambda_2)e^{-i(r-q+m-o+(j_2-j_1)N)\lambda_2} \,d\lambda_2   =O\Big(\frac{1}{T}\frac{ g^2(k)}{N^{1-2D-2\epsilon}}\Big)
\end{eqnarray*}
\item [(ii)]
\noindent
\vspace{-1.1cm}
\begin{eqnarray*}
&&  \frac{1}{ N^2}\frac{1}{M^2} \sum_{j_1=1}^M \sum_{\substack{q,r=0  }}^{N-1}\sum_{\substack{l,m,n,o=0  }}^{\infty} \sum_{\substack{j_2=1  \\ -(N-1) \leq r+l-n+(j_2-j_1)N \leq -1 \\ 0 \leq q+o-m+(j_1-j_2)N \leq N-1 \\ |l-n+(j_2-j_1)N| \leq N-1 \\ |o-m+(j_1-j_2)N| \leq N-1 }}^M \psi_l(u_{j_1}) \psi_m(u_{j_1})\psi_n(u_{j_2})  \psi_o(u_{j_2})\\
&&\int_{-\pi}^{\pi} \phi_T(u_{j_1},\lambda_1)e^{-i(r-q+l-n+(j_2-j_1)N)\lambda_1} \,d\lambda_1
\int_{-\pi}^{\pi} \phi_T(u_{j_2},\lambda_2)e^{-i(r-q+m-o+(j_2-j_1)N)\lambda_2} \,d\lambda_2   =O\Big(\frac{1}{T}\frac{ g^2(k)}{N^{1-2D-2\epsilon}}\Big)
\end{eqnarray*}
\item [(iii)]
\noindent
\vspace{-1.1cm}
\begin{eqnarray*}
&&  \frac{1}{ N^2}\frac{1}{M^2} \sum_{j_2=1}^M \sum_{\substack{q=0  }}^{N-1}
 \sum_{\substack{r=-(N-1)  }}^{-1}
\sum_{\substack{l,m,n,o=0  }}^{\infty} \sum_{\substack{j_1=1    \\ 0 \leq q+o-m+(j_1-j_2)N \leq N-1 \\ |l-n+(j_2-j_1)N| \leq N-1 \\ |o-m+(j_1-j_2)N| \leq N-1 }}^M \psi_l(u_{j_1}) \psi_m(u_{j_1})\psi_n(u_{j_2})  \psi_o(u_{j_2})\\
&&\int_{-\pi}^{\pi} \phi_T(u_{j_1},\lambda_1)e^{-i(r-q+l-n+(j_2-j_1)N)\lambda_1} \,d\lambda_1
\int_{-\pi}^{\pi} \phi_T(u_{j_2},\lambda_2)e^{-i(r-q+m-o+(j_2-j_1)N)\lambda_2} \,d\lambda_2   =O\Big(\frac{1}{T}\frac{ g^2(k)}{N^{1-2D-2\epsilon}}\Big)
\end{eqnarray*}

\item [(iv)]
\noindent
\vspace{-1.1cm}
\begin{eqnarray*}
&&  \frac{1}{ N^2}\frac{1}{M^2} \sum_{j_1=1}^M \sum_{\substack{r,q=0  }}^{N-1}
\sum_{\substack{l,m,n,o=0  }}^{\infty} \sum_{\substack{j_2=1    \\ 0 \leq q+o-m+(j_1-j_2)N \leq N-1 \\ N \leq  |l-n+(j_2-j_1)N|  \\ |o-m+(j_1-j_2)N| \leq N-1 }}^M \psi_l(u_{j_1}) \psi_m(u_{j_1})\psi_n(u_{j_2})  \psi_o(u_{j_2})\\
&&\int_{-\pi}^{\pi} \phi_T(u_{j_1},\lambda_1)\phi_T(u_{j_2},\lambda_1)e^{-i(r-q+l-n+(j_2-j_1)N)\lambda_1} \,d\lambda_1
\int_{-\pi}^{\pi} e^{-i(r-q+m-o+(j_2-j_1)N)\lambda_2} \,d\lambda_2   =O\Big(\frac{1}{T}\frac{ g^2(k)}{N^{1-2D-2\epsilon}}\Big)
\end{eqnarray*}
\item [(v)]
\noindent
\vspace{-1.1cm}
\begin{eqnarray*}
&& \frac{1}{ N^2}\frac{1}{M^2} \sum_{j_1=1}^M \sum_{\substack{r,q=0  }}^{N-1}  \sum_{\substack{l,m,n,o=0 }}^{\infty} \sum_{\substack{j_2=1\\ |j_1-j_2| \geq 1  \\ 0 \leq q+o-m+(j_1-j_2)N \leq N-1 \\ |o-m+(j_1-j_2)N| \leq N-1  }}^M \psi_l(u_{j_1}) \psi_m(u_{j_1})\psi_n(u_{j_2})  \psi_o(u_{j_2})\\
&&\int_{-\pi}^{\pi} \phi_T(u_{j_1},\lambda_1)\phi_T(u_{j_2},\lambda_1)e^{-i(r-q+l-n+(j_2-j_1)N)\lambda_1} \,d\lambda_1
\int_{-\pi}^{\pi} e^{-i(r-q+m-o+(j_2-j_1)N)\lambda_2} \,d\lambda_2  =O\Big(\frac{1}{T}\frac{ g^2(k)}{N^{1-2D-2\epsilon}}\Big)
\end{eqnarray*}
\item [(vi)]
\noindent
\vspace{-1.1cm}
\begin{eqnarray*}
&& \frac{1}{ N^2}\frac{1}{M^2}\sum_{j_1=1}^M \sum_{\substack{q=0  }}^{N-1} \sum_{\substack{r \in \mathbb{Z}\\|r| \geq N }}^{}\sum_{\substack{l,m,n,o=0  }}^{\infty}
\sum_{\substack{j_2=1  \\ 0 \leq q+o-m+(j_1-j_2)N \leq N-1 \\ |l-n+(j_2-j_1)N| \leq N-1 \\ |o-m+(j_1-j_2)N| \leq N-1}}^{\infty}
\psi_l(u_{j_1}) \psi_m(u_{j_1})\psi_n(u_{j_2})  \psi_o(u_{j_2})\\
&&\int_{-\pi}^{\pi} \phi_T(u_{j_1},\lambda_1)e^{-i(r-q+l-n+(j_2-j_1)N)\lambda_1} \int_{-\pi}^{\pi}\Big[ \phi_T(u_{j_2},\lambda_2)-\phi_T(u_{j_2},\lambda_1) \Big]e^{-i(r-q+m-o+(j_2-j_1)N)\lambda_2} \,d\lambda_2  \,d\lambda_1\\
&=&O\Big(\frac{1}{T}\frac{ g^2(k)}{N^{1-2D-2\epsilon}}\Big)
\end{eqnarray*}
\end{itemize}
\end{lem}
\textbf{Proof:} Without loss of generality we restrict ourselves to a proof of part (i) and (v) and note that all other claims are proven by using the same arguments. \\

\textit{Proof of (i):}
We use (\ref{apprpsi}), (\ref{Ordnung theta}) and Lemma \ref{fourier} to bound the term in (i) (up to a constant) through
\begin{eqnarray*}
&& \frac{g^2(k)}{ N^2}\frac{1}{M^2} \sum_{j_1=1}^M \sum_{\substack{q,r=0  }}^{N-1}\sum_{\substack{l,m,n,o=1 }}^{\infty} \sum_{\substack{j_2=1 \\ N \leq |r+l-n+(j_2-j_1)N| \\ 0 \leq q+o-m+(j_1-j_2)N \leq N-1 \\ |l-n+(j_2-j_1)N| \leq N-1 \\ |o-m+(j_1-j_2)N| \leq N-1 \\ 1 \leq |r-q+m-o+(j_2-j_1)N|}}^{M} \frac{1}{l^{1-d_0(u_{j_1})}}\frac{1}{m^{1-d_0(u_{j_1})}}\frac{1}{n^{1-D}}\frac{1}{o^{1-D}}\\
&&\frac{1}{|r-q+l-n+(j_2-j_1)N|^{1+2d_0(u_{j_1})-\epsilon}}
\frac{1}{|r-q+m-o+(j_2-j_1)N|^{1-\epsilon}}.
\end{eqnarray*}
If the variables $j_1, o$ and $m$ are fixed, it follows with the constraint $0 \leq q+o-m+(j_1-j_2)N \leq N-1$ that there are at most two possible values for $j_2$ such that the resulting term is non vanishing. We now discuss for which combinations of $j_1$ and $j_2$ the above expression is maximized and then restrict ourselves to the resulting pair $(j_1,j_2)$.

If $j_1$ and $j_2$ are given, the variables $l,m,n,o$ can only be chosen such that $ |l-n+(j_2-j_1)N| \leq N-1$ and $ |o-m+(j_1-j_2)N| \leq N-1$ are fulfilled. Therefore, the possible values of the fractions
$(|r-q+l-n+(j_2-j_1)N|)^{-1}(|r-q+m-o+(j_2-j_1)N|)^{-1}$ are the same for any combination of $j_1$ and $j_2$. Consequently, in order to maximize the term above we need to maximize $l^{-1d_0(u_{j_1})} m^{-1+d_0(u_{j_1})}  n^{-1+D} o^{-1+D}$, which is achieved by the choice $j_1=j_2$ [since then $l,m,n,o$ can be jointly taken as small as possible due to the constraints $|l-n+(j_2-j_1)N| \leq N-1$ and $|o-m+(j_1-j_2)N| \leq N-1$]. Hence we can bound that above expression (up to a constant) by
\begin{eqnarray*}
& & \frac{g^2(k)}{ N^2}\frac{1}{M^2} \sum_{j_1=1}^M \sum_{\substack{q,r=0  }}^{N-1}\sum_{\substack{l,m,n,o=1 \\ N \leq |r+l-n|  \\ |l-n| \leq N-1 \\ |o-m| \leq N-1 \\ 1 \leq |r-q+m-o|}}^{\infty} \frac{1}{l^{1-d_0(u_{j_1})}}\frac{1}{m^{1-d_0(u_{j_1})}}\frac{1}{n^{1-D}}\frac{1}{o^{1-D}}\frac{1}{|r-q+l-n|^{1+2d_0(u_{j_1})-\epsilon}}
\frac{1}{|r-q+m-o|^{1-\epsilon}}.
\end{eqnarray*}
By setting $g:=r+l-n$ and $s:=q+o-m$ this term can be written as
\begin{eqnarray*}
&& \frac{ g^2(k)}{N^2}\frac{1}{M^2}  \sum_{j_1=1}^M  \sum_{\substack{q,r,s=0  \\ 1 \leq |r-s| }}^{N-1}\sum_{\substack{g \in \mathbb{Z} \\ |g| \geq N }}^{}\sum_{\substack{m,n=1\\  1 \leq g-r+n\\ 1 \leq s-q+m\\ |g-r| \leq N-1 }}^{\infty}\frac{1}{(g-r+n)^{1-d_0(u_{j_1})}}\frac{1}{m^{1-d_0(u_{j_1})}}\frac{1}{n^{1-D}}\\
&&\times\frac{1}{(s-q+m)^{1-D}}\frac{1}{|g-q|^{1+2d_0(u_{j_1})-\epsilon}}\frac{1}{|r-s|^{1-\epsilon}}
\end{eqnarray*}

Through an repeated application of (\ref{Gleichung 3.196(3)}) and (\ref{eq1}) the claim now follows. \\

\textit{Proof of (v):} By setting
\begin{eqnarray*}
f(u_{j_1},u_{j_2}, \lambda):= \frac{1}{2\pi} \sum_{\substack{l,n=0}}^{\infty}
\psi_l(u_{j_1}) \psi_n(u_{j_2}) e^{-i(l-n)\lambda}.
\end{eqnarray*}
we can write the term in (v) as
\begin{eqnarray*}
&& \frac{2\pi}{ N^2}\frac{1}{M^2} \sum_{j_1=1}^M \sum_{\substack{r,q=0  }}^{N-1}  \sum_{\substack{m,o=0 }}^{\infty} \sum_{\substack{j_2=1\\ |j_1-j_2| \geq 1  \\ 0 \leq q+o-m+(j_1-j_2)N \leq N-1 \\ |o-m+(j_1-j_2)N| \leq N-1  }}^M  \psi_m(u_{j_1})\psi_o(u_{j_2})\\
&&\times\int_{-\pi}^{\pi} \phi_T(u_{j_1},\lambda_1)\phi_T(u_{j_2},\lambda_1)f(u_{j_1},u_{j_2}, \lambda_1)e^{-i(r-q+(j_2-j_1)N)\lambda_1} \,d\lambda_1
\int_{-\pi}^{\pi} e^{-i(r-q+m-o+(j_2-j_1)N)\lambda_2} \,d\lambda_2.
\end{eqnarray*}
and by integrating over $\lambda_2$ this is the same as
\begin{eqnarray*}
&& \frac{4\pi^2}{ N^2}\frac{1}{M^2} \sum_{j_1=1}^M \sum_{\substack{q=0  }}^{N-1}  \sum_{\substack{m,o=0 }}^{\infty} \sum_{\substack{j_2=1\\ |j_1-j_2| \geq 1  \\ 0 \leq q+o-m+(j_1-j_2)N \leq N-1 \\ |o-m+(j_1-j_2)N| \leq N-1  }}^M  \psi_m(u_{j_1})\psi_o(u_{j_2})\int_{-\pi}^{\pi} \phi_T(u_{j_1},\lambda_1)\phi_T(u_{j_2},\lambda_1)f(u_{j_1},u_{j_2}, \lambda_1)e^{-i(o-m)\lambda_1} \,d\lambda_1 .
\end{eqnarray*}
By (\ref{Ordnung theta}) and Lemma \ref{fourier} this sum can be  bounded by
\begin{eqnarray*}
&& \frac{Cg^2(k)}{ N^2}\frac{1}{M^2} \sum_{j_1=1}^M \sum_{\substack{q=0  }}^{N-1}  \sum_{\substack{m,o=1  }}^{\infty} \sum_{\substack{j_2=1\\ |j_1-j_2| \geq 1   \\ 0 \leq q+o-m+(j_1-j_2)N \leq N-1 \\ |o-m+(j_1-j_2)N| \leq N-1  }}^M  \frac{1}{m^{1-d_0(u_{j_1})}} \frac{1}{o^{1-d_0(u_{j_1})}}\frac{1}{|o-m|^{1+d_0(u_{j_1})+d_0(u_{j_2})-2\epsilon}}\\
&\leq & \frac{C g^2(k)}{ N^2}\frac{1}{M^2} \sum_{j_1=1}^M \sum_{\substack{q=0  }}^{N-1}  \sum_{\substack{m,o=1  }}^{\infty} \sum_{\substack{j_2=1\\ |j_1-j_2| \geq 1   \\ 0 \leq q+o-m+(j_1-j_2)N \leq N-1 \\ |o-m+(j_1-j_2)N| \leq N-1  }}^M  \frac{1}{m^{1-D}} \frac{1}{o^{1-D}}\frac{1}{|o-m|^{1-2\epsilon}} .
\end{eqnarray*}
As in the proof of (i) we can argue  that there are at most two possible values for $j_2$ if $o,m$ and $j_1$ are chosen and that the expression is maximized for $|j_1-j_2|=1$. Therefore we can bound the above expression up to a constant through
\begin{eqnarray*}
&& \frac{g^2(k)}{ N^2}\frac{1}{M} \sum_{\kappa \in \{ -1,1\}} \sum_{\substack{q=0  }}^{N-1}  \sum_{\substack{m,o=1 \\ 0 \leq q+o-m+\kappa N \leq N-1 \\ |o-m+\kappa N| \leq N-1  }}^{\infty}  \frac{1}{m^{1-D}} \frac{1}{o^{1-D}}\frac{1}{|o-m|^{1-2\epsilon}}.
\end{eqnarray*}
By setting $p:=o-m+\kappa N$ the claim follows with (\ref{eq1}).  $\hfill \Box$

\end{document}